\documentclass[a4paper]{amsart}

\usepackage{epsfig}
\usepackage{graphicx}
\usepackage{graphicx}
\usepackage{amsmath}
\usepackage{amssymb}
\input xy
\xyoption{all}

\renewcommand{\mod}{\operatorname{mod}\nolimits}
\newcommand{\ann}{\operatorname{ann}\nolimits}
\newcommand{\add}{\operatorname{add}\nolimits}
\newcommand{\rad}{\operatorname{rad}\nolimits}
\newcommand{\Hom}{\operatorname{Hom}\nolimits}
\newcommand{\op}{\operatorname{op}\nolimits}
\newcommand{\bo}{\operatorname{b}\nolimits}

\newcommand{\Sub}{\operatorname{Sub}\nolimits}
\newcommand{\Ext}{\operatorname{Ext}\nolimits}
\newcommand{\id}{\operatorname{id}\nolimits}
\newcommand{\Mod}{\operatorname{Mod}\nolimits}
\newcommand{\End}{\operatorname{End}\nolimits}

\newcommand{\Tor}{\operatorname{Tor}\nolimits}
\newcommand{\fl}{\operatorname{f.l.}\nolimits}
\newcommand{\CM}{\operatorname{CM}\nolimits}
\newcommand{\pd}{\operatorname{pd}\nolimits}
\newcommand{\Pfaff}{\operatorname{Pfaff}\nolimits}
\newcommand{\gl}{\operatorname{gl.dim}\nolimits}
\newcommand{\pr}{\operatorname{pr}\nolimits}
\newcommand{\la}{\Lambda}
\newcommand{\ga}{\Gamma}
\newcommand{\ul}{\underline}

\newcommand{\SO}{\operatorname{SO}\nolimits}
\newcommand{\SL}{\operatorname{SL}\nolimits}

\renewcommand{\Im}{\operatorname{Im}\nolimits}
\newcommand{\Ker}{\operatorname{Ker}\nolimits}
\newcommand{\xto}{\xrightarrow}

\newcommand{\Fac}{\operatorname{Fac}\nolimits}
\newcommand{\A}{\operatorname{\bf A}\nolimits}
\newcommand{\B}{\operatorname{\mathcal B}\nolimits}
\newcommand{\C}{\operatorname{\mathcal C}\nolimits}
\newcommand{\M}{\operatorname{\mathcal M}\nolimits}
\newcommand{\X}{\operatorname{\bf X}\nolimits}
\newcommand{\F}{\operatorname{\bf F}\nolimits}
\newcommand{\D}{\operatorname{\mathcal D}\nolimits}
\newcommand{\E}{\operatorname{\mathcal E}\nolimits}

\newcommand{\Q}{\operatorname{\mathbb Q}\nolimits}
\newcommand{\Sc}{\operatorname{\bf S}\nolimits}
\newcommand{\ZZ}{\operatorname{\mathbb Z}\nolimits}

\newcommand{\rev}{\operatorname{rev}}

\newcommand{\s}{\scriptscriptstyle}
\newcommand{\p}{\scriptstyle}

\newcommand{\Res}{\operatorname{Res}}

\newcommand{\GL}{\operatorname{GL}}

\newcommand{\iso}{\operatorname{iso}}
\newcommand{\Gr}{\operatorname{Gr}}

\newcommand{\Di}{\displaystyle}

\newtheorem{theorem}{Theorem}[subsection]

\newtheorem{question}[theorem]{Question}
\newtheorem{conjecture}[theorem]{Conjecture}
\newtheorem{corollary}[theorem]{Corollary}
\newtheorem{lemma}[theorem]{Lemma}
\newtheorem{proposition}[theorem]{Proposition}

\begin{document}
\title{Cluster structures for 2-Calabi-Yau categories and unipotent groups}

\author[Buan]{A. B. Buan}
\address{Institutt for matematiske fag\\
Norges teknisk-naturvitenskapelige universitet\\
N-7491 Trondheim\\
Norway}
\email{aslakb@math.ntnu.no}

\author[Iyama]{O. Iyama}
\address{Graduate school of Mathematics \\
Nagoya University, Chikusa-ku\\
Nagoya\\
464-8602 Nagoya \\
Japan}
\email{iyama@math.nagoya-u.ac.jp}

\author[Reiten]{I. Reiten}
\address{Institutt for matematiske fag\\
Norges teknisk-naturvitenskapelige universitet\\
N-7491 Trondheim\\
Norway}
\email{idunr@math.ntnu.no}

\author[Scott]{J. Scott} \thanks{All authors were supported by a STORFORSK-grant 167130 from the Norwegian Research Council}
\address{Dept. of Pure Mathematics \\
University of Leeds \\
Leeds LS2 9JT \\
United Kingdom \\
England 
}
\email{jscott@maths.leeds.ac.uk}

\maketitle

\begin{abstract}
We investigate cluster tilting objects (and subcategories) in triangulated 2-Calabi-Yau 
categories and related categories. In particular we construct a new class of such categories 
related to preprojective algebras of non-Dynkin quivers 
associated with elements in the Coxeter group. This class of 2-Calabi-Yau categories
contains the cluster categories and the stable categories of preprojective algebras of Dynkin 
graphs as special cases. For these 2-Calabi-Yau categories we construct cluster tilting objects 
associated with each reduced expression. The associated quiver is described in terms of 
the reduced expression.
Motivated by the theory of cluster 
algebras, we formulate the notions of (weak) cluster structure and substructure, and give 
several illustrations of these concepts. We give applications to cluster algebras and 
subcluster algebras related to unipotent groups, both in the Dynkin and non-Dynkin case.
\end{abstract}

\tableofcontents


\section*{Introduction}
The theory of cluster algebras, initiated by Fomin-Zelevinsky in \cite{fz1}, and further
developed in a series of papers, including \cite{fz2,fz3,fz4}, has turned out to have 
interesting connections with many parts of algebra and other branches of mathematics. One of the 
links is with the representation theory of algebras, where a first connection was discovered in \cite{mrz}. 
A philosophy has been to model 
the main ingredients in the definition of a cluster algebra in a categorical/module theoretical setting. 
The cluster categories associated with  finite dimensional hereditary algebras were introduced for this 
purpose in \cite{bmrrt}, and shown to be triangulated in \cite{k1} (see also \cite{ccs} for the $A_n$ case), 
and the module categories $\mod\Lambda$ for $\la$ a preprojective algebra of a Dynkin quiver have been used for a  
similar purpose \cite{gls1}. This development has both inspired new directions of investigations on the categorical 
side, as well as interesting feedback on the theory of cluster algebras, see for example
\cite{abs,bmr1,bmr2,bm,bmrt,cc,ck1,ck2,gls1,gls3,hub,i1,i2,ir,iy,it,kr1,kr2,ringel,t} for material related to this paper.

Both the cluster categories and the stable categories $\ul{\mod}\la$ of preprojective algebras 
are triangulated Calabi-Yau categories of dimension 2 
(2-CY  for short). They both have what is called cluster tilting objects/subcategories \cite{bmrrt,kr1,iy} 
(called maximal 1-orthogonal in \cite{i1}), which are important since
they are the analogs of 
clusters. The investigation of cluster tilting objects/subcategories in 2-CY categories and related 
categories is interesting both from the point 
of view of cluster algebras and in itself. Hence it is of interest to develop methods for constructing 
2-CY categories together 
with the special objects/subcategories, and this is the main purpose of
the first two chapters.

The properties of cluster tilting objects in ($\Hom$-finite) 2-CY
categories which have been important for applications to cluster
algebras are (a) the unique exchange property for indecomposable
summands of cluster tilting objects, (b) the existence of associated
exchange triangles, (c) having no loops or 2-cycles (in the quiver of the
endomorphism algebra of a cluster tilting object) and (d) when passing
from the endomorphism algebra of a cluster tilting object $T$ to
the endomorphism algebra of another one $T^*$ via an exchange, the change in quivers is given by
Fomin-Zelevinsky mutation. The properties (a) and (b) are known to hold for any
2-CY triangulated category \cite{iy}, proved for cluster categories
in \cite{bmrrt} and for stable categories of preprojective algebras of
Dynkin type in \cite{gls1}. The property (c) does not always hold (see
\cite{bikr}), and hence it is of interest to establish criteria for
this to be the case, which is one of the topics in this paper. We then
show for any 2-CY categories that if (c) holds, then also (d) follows, as
previously shown by Palu for algebraic triangulated categories \cite{p}. We
construct new 2-CY categories with cluster tilting objects from old
ones via some subfactor construction, extending results from \cite{iy},
with a main focus on how condition (c) behaves under this
construction. Associated with this we introduce the notions of
cluster structures and cluster substructures.

Important examples, investigated in \cite{gls1}, are the categories
$\mod\la$ of finitely generated modules over the preprojective algebra
$\la$ of a Dynkin quiver. We deal with appropriate subcategories of
$\mod\la$. The main focus
in this paper is on the more general case of subcategories of the
category $\fl\la$ of finite length modules over the completion of the
preprojective algebra of a non-Dynkin quiver with no loops. Our main
tool is to extend the tilting theory developed for $\la$ in the
noetherian case in \cite{ir}.
%
This turns out to give a large class of 2-CY
categories associated with elements in the corresponding Coxeter
groups. For these categories we construct cluster tilting objects associated
with each reduced expression, and we describe the associated quiver
directly in
terms of the reduced expression. We prove that this class of 2-CY categories
contains all the cluster categories of finite dimensional
hereditary algebras and the stable categories $\ul{\mod}\la$ for a
preprojective algebra $\Lambda$ of Dynkin type. This also allows us to
get more information on the latter case.

We illustrate with applications to constructing subcluster algebras of cluster algebras, a notion 
which we define here, 
and which is already implicit in the literature. 
For this we define, inspired by maps from \cite{cc,ck1} and \cite{gls1},
(strong) cluster maps. These maps have the property that we can pass from cluster structures
and substructures to cluster algebras and subcluster algebras.

Associated with substructures for preprojective 
algebras of Dynkin type, 
we discuss examples from $\SO_8(\mathbb{C})$-isotropic Grassmanians, for the $G_{2,5}$ Schubert variety  
and for a unipotent 
cell of the unipotent subgroup of $\SL_4(\mathbb{C})$. For preprojective algebras of extended 
Dynkin type we use our results to 
investigate cluster structures for affine unipotent cells, in some special cases for the loop 
group $\SL_2(\mathcal{L})$. 

For a (non-Dynkin) quiver $Q$ with associated Coxeter group $W$ we can
for each $w\in W$ consider the coordinate ring $\mathbb{C}[U^w]$ of
the unipotent cell associated with $w$ in the corresponding Kac-Moody
group. We conjecture that this ring has a cluster algebra structure,
and that it is modelled by our (stably) 2-CY category associated
with the same $w$. As a support for this we prove the conjecture for
the case $\widehat{A_1}$ for the word $w$ of length at most 4.

The first chapter is devoted to introducing and investigating the notions of  
cluster structures and substructures, 
and giving sufficient conditions for such structures to occur. Also the three concrete examples 
mentioned above
are investigated, and used to illustrate the
connection with cluster algebras and subcluster algebras in Section \ref{c3_sec2}. In Chapter \ref{chap2} 
we use tilting theory to construct categories whose stable categories are 2-CY, 
along with natural cluster tilting objects in these categories. In
Section \ref{c3_sec2} we illustrate with examples for
preprojective algebras of Dynkin type, 
and in Section \ref{c3_sec3} we illustrate with examples from the extended Dynkin case.

Part of this work was done independently by Geiss-Leclerc-Schr\"oer  in \cite{gls3}. 
For Chapter \ref{chap1}, this concerns the 
development of 2-CY categories (in a different language) in the case of 
subcategories of the form $\Sub P$ (or $\Fac P$) for $P$ 
projective, over a preprojective algebra of Dynkin type, with a somewhat 
different approach. Concerning Section \ref{c3_sec2}, 
examples arising from $\Sub P$ were done independently in \cite{gls3}. For 
this connection, the last author was inspired 
by a lecture of Leclerc in 2005, where cluster algebras associated with $\Sub P$ 
in the $A_n$-case were discussed. There is also recent work of Geiss-Leclerc-Schr\"oer \cite{gls5} 
related to Chapter \ref{chap2}, where completely different methods are used.

For general background on representation theory of algebras, we refer
to \cite{ars, ass, rin, h1, ahk}, and for Lie theory we refer to \cite{bl}.

Our modules are usually left modules and composition of maps $fg$ means
first $f$, then $g$.

The second author would like to thank William Crawley-Boevey and Christof Geiss for 
answering a question on references about 2-CY property of preprojective algebras.
He also would like to thank Bernald Leclerc for valuable comments.

\section{2-CY categories and substructures}\label{chap1}

The cluster algebras of Fomin and Zelevinsky have motivated work on
trying to model the essential ingredients in the definition of a
cluster algebra in a categorical /module theoretical way. In particular, 
this led to the theory of cluster categories and the investigation of new aspects 
of the module theory of preprojective algebras of Dynkin type. In Section \ref{c1_sec1} 
we give some of the main categorical requirements needed for the modelling, for the cases
with and without coefficients, leading to the notions of weak cluster structure 
and cluster structure. The main examples do, like the above mentioned
examples, have 2-Calabi-Yau type properties.

We introduce substructures of (weak) cluster structures in Section \ref{c1_sec2}. 
For this it is natural to deal with (weak) cluster structures with what will be called
coefficients, at least for the substructures. Of particular interest
for our applications to cluster algebras is the case of completions of 
preprojective algebras $\Lambda$ of a finite connected quiver with
no loops over an algebraically closed field $k$, where the
interesting larger category is the stable category
$\underline{\fl}\Lambda$ of the finite length $\Lambda$-modules. For
Dynkin quivers this is the stable category $\underline{\mod}\Lambda$
of the finitely generated $\Lambda$-modules, and in the non-Dynkin
case $\underline{\fl}\Lambda =\fl\Lambda$. The first case is discussed in Section \ref{c1_sec3},
while Chapter \ref{chap2} is devoted to the second case.\\
\hspace{7mm}

\subsection{Cluster structures}\label{c1_sec1}
${}$ \\
In this section we introduce the concepts of weak cluster structure and cluster 
structure for extension closed subcategories of 
triangulated categories or for exact categories. We illustrate with 2-CY categories 
and closely related categories, and the main objects we investigate are the cluster 
tilting ones. These cases are particularly nice when the quivers of the cluster 
tilting subcategories 
have no loops or 2-cycles. Also the closely related maximal rigid objects (see \cite{gls1}) provide interesting examples.

We start with introducing the notions of weak cluster structure and cluster structure. 
Throughout this chapter all categories are Krull-Schmidt categories 
over an algebraically closed field $k$, that is, each object is isomorphic to a finite
direct sum of indecomposable objects with local endomorphism ring. The categories we consider 
are either exact (for example abelian) categories or extension closed subcategories of triangulated
categories. Note that an extension closed subcategory of an exact category is again exact.
We refer to \cite{k2,k3} for definition and basic properties 
of exact categories, which behave very much like abelian categories, also with respect to 
derived categories and Ext-functors. 

We often identify a set of indecomposable objects with the additive
subcategory consisting of all summands of direct sums of these
indecomposable objects. We also identify an object with the set of 
indecomposable objects appearing in a direct sum decomposition,
and with the subcategory obtained in the above way.

Assume that we have a collection of sets $\underline{x}$ 
(which may be infinite) of non-isomorphic indecomposable objects, 
called \emph{clusters}. The union of all indecomposable objects in clusters are called 
\emph{cluster variables}. Assume also that there is a subset $\ul{p}$
(which may be infinite) of indecomposable objects which are not cluster
variables, called  {\em coefficients}. We denote by $T$
the union of the indecomposable objects in $\ul{x}$ and $\ul{p}$, sometimes viewed as a category with these objects, 
and call it an \emph{extended cluster}. 

We say that the clusters, together with the prescribed set of coefficients $\ul{p}$, give a 
\emph{weak cluster structure} on $\C$ if the following hold:

\begin{itemize}
\item[(a)]{For each extended cluster $T$ and each cluster variable $M$ which is a summand in $T$, there is  a unique 
indecomposable object $M^{\ast}\not\simeq
M$ such  that we get a new extended cluster $T^{\ast}$ by replacing $M$ by $M^{\ast}$. 
We denote this operation, called {\em exchange}, by $\mu_{M}(T)=T^{\ast}$,
and we call $(M, M^{\ast})$ an \emph{exchange pair}.}
\item[(b)]{There are triangles/short exact sequences
$M^{\ast}\xrightarrow{f}B\xrightarrow{g}M$ and
$M\xrightarrow{s}B^{'}\xrightarrow{t}M^{\ast}$, where the maps $g$ and $t$ are
minimal right $\add(T\char92 \{M\})$-approximations and $f$ and $s$
are minimal left $\add(T\char92 \{M\})$-approximations.
These are called {\em exchange triangles/sequences}.}
\end{itemize}

Denote by $Q_{T}$  the quiver of 
$T$, where the vertices correspond to the indecomposable objects in $T$ and the number of arrows $T_i \to T_j$ 
between two indecomposable objects $T_i$ and $T_j$ is
given by the dimension of the space of irreducible maps $\rad(T_i,T_j)/\rad^2(T_i,T_j)$.
Here $\rad(\ ,\ )$ denotes the radical in $\add T$, where the objects are finite direct sums
of objects in $T$. For an algebra $\la$ (where $\la$ has a unique decomposition as a direct sum of indecomposable objects),
the quiver of $\la$ is then the opposite of the quiver of $\add \la$.

We say that a quiver $Q= (Q_0,Q_1)$ is an {\it extended quiver} with respect to a
subset of vertices $Q_0'$ if there are no
arrows between two vertices in $Q_0\backslash Q_0'$. We regard the quiver $Q_T$ of
an extended cluster as an extended quiver by neglecting all arrows
between two vertices corresponding to coefficients.

We say that $\C$, with a fixed set of clusters and coefficients, has
{\it no loops} (respectively, {\em no 2-cycles}) if in the extended
quiver of each extended cluster there are no loops (respectively, no 2-cycles).
When $\underline{x}$ is finite, this is the opposite quiver of the factor algebra $\underline{\End}(T)$ of $\End(T)$ 
by the maps factoring through direct sums of objects from $\underline{p}$. 

We say that we have a \emph{cluster structure} if the following additional conditions hold:
\begin{itemize}
\item[(c)]There are no loops or 2-cycles. (In other words, for a cluster variable $M$, any 
non-isomorphism $u\colon M\to M$ factors through $g\colon B\to M$ and through $s\colon M\to B'$, 
and any non-isomorphism $v\colon M^{\ast}\to M^{\ast}$ factors through $f\colon M^{\ast}\to B$ and 
through $t\colon B'\to M^{\ast}$, and $B$ and $B'$ have no common indecomposable summand.)
\item[(d)] For an extended cluster $T$, passing from $Q_{T}$ to $Q_{T^{\ast}}$ is given by Fomin-Zelevinsky mutation 
at the vertex of $Q_{T}$ given by the cluster variable $M$.
\end{itemize}

Note that (c) is needed for (d) to make sense, but it is still convenient to write two separate statements.

We recall that for an extended quiver $Q$ without loops or 2-cycles
and a vertex $i$ in $Q_0'$, the Fomin-Zelevinsky
mutation $\mu_i(Q)$ of $Q$ at $i$ is the quiver obtained from $Q$
making the following changes \cite{fz1}:
\begin{itemize}
\item[-]{Reverse all the arrows starting or ending at $i$.}
\item[-]{Let $s\neq i$ and $t\neq i$ be vertices in $Q_0$ such that at
least one vertex belongs to $Q_0'$.
If we have $n>0$ arrows from $t$ to $i$ and $m>0$ arrows from $i$
to $s$ in $Q$ and $r$ arrows from $s$ to $t$ in $Q$ (interpreted as $-r$ arrows from $t$ to $s$ if $r<0$), 
then we have $nm-r$
arrows from $t$ to $s$ in the new quiver $\mu_i(Q)$ (interpreted as $r-nm$
arrows from $s$ to $t$ if $nm-r<0$).}
\end{itemize}

The main known examples of triangulated $k$-categories with finite
dimensional homomorphism spaces (Hom-finite for short) which have a weak cluster structure, and usually 
cluster structure, are 2-CY categories. These are triangulated $k$-categories with  functorial 
isomorphisms $D\Ext^1(A,B)\simeq \Ext^1(B,A)$ for  all $A,B$ in $\C$, where $D=\Hom_k(\ ,k)$.
A $\Hom$-finite triangulated category is 2-CY if and only if it
has almost split triangles with translation $\tau$ and $\tau \colon \C
\to \C$ is a functor isomorphic to the shift functor $[1]$ (see
\cite{rv}). 

We have the following examples of 2-CY categories.

\noindent 
(1) The cluster category $\C_{H}$ associated with a finite
dimensional hereditary $k$-algebra $H$ is by definition the orbit
category ${\bf D}^{\bo}(H)/\tau^{-1}[1]$, where ${\bf D}^{\bo}(H)$ is
the bounded derived category of finitely generated $H$-modules, and
$\tau$ is the AR-translation of ${\bf D}^{\bo}(H)$ \cite{bmrrt}. It is
a $\Hom$-finite triangulated category \cite{k2}, and it is 2-CY since $\tau=[1]$.

\bigskip
\noindent
(2) The stable category of maximal Cohen-Macaulay modules $\underline{\CM}(R)$ over a 3-dimensional
complete local commutative noetherian Gorenstein isolated singularity
$R$ containing the residue field $k$ \cite{a} (see \cite{yo}).

\bigskip
\noindent 
(3) The preprojective algebra $\Lambda$ associated to a finite connected quiver $Q$ 
without loops is defined  as follows: 
Let $\widetilde{Q}$ be the quiver constructed from $Q$ by adding an arrow $\alpha^{\ast}\colon i\to j$ 
for each arrow $\alpha\colon j\to i$ in $Q$. 
Then $\la=k\widetilde{Q}/I$, where $I$
is the ideal generated by the sum of commutators
$\sum_{\beta\in Q_1}[\beta,\beta^{\ast}]$.
Note that $\Lambda$ is uniquely determined up to isomorphism by the
underlying graph of $Q$.

When $\Lambda$ is the preprojective algebra of a Dynkin quiver over $k$,
the stable category $\ul{\mod}\la$ is 2-CY (see \cite[3.1,1.2]{ar2}\cite{cb}\cite[8.5]{k2}). 

When $\Lambda$ is the completion of the preprojective algebra of a
finite connected quiver without loops which is not Dynkin,
the bounded derived category ${\bf D}^{\bo}(\fl\Lambda)$
of the category $\fl\la$ of the modules of finite length is 2-CY
(see \cite{b,cb,bbk}\cite[section 8]{gls2}).

\bigskip
We shall also use the terminology 2-CY in more general situations.
Note that from now on we will 
usually write just ``category'' instead of ``$k$-category''.
 
We say that an exact Hom-finite category $\C$ is {\em derived 2-CY} if the triangulated category 
${\bf D}^{\bo}(\C)$ is 2-CY, i.e. if
$D\Ext^i(A,B)\simeq\Ext^{2-i}(B,A)$ for all $A$, $B$ in ${\bf D}^{\bo}(\C)$ and all $i$.  
Note that 
when $\C$ is derived 2-CY, then $\C$ has no non-zero projective or injective objects. 
The category $\fl\Lambda$ where $\la$ is the completion of the preprojective algebra of a non-Dynkin 
connected quiver without loops is an important example of a derived 2-CY category.

We say that an exact category $\C$ is {\em stably 2-CY} if
it is Frobenius, that is, $\C$ has enough projectives and injectives, which coincide, 
and the stable category $\underline{\C}$, which is triangulated \cite{h1}, is Hom-finite 2-CY.
Recall that $\C$ is said to have enough projectives if for each $X$ in $\C$ there is an
exact sequence $0 \to Y \to P \to X \to 0$ in $\C$ with $P$ projective. Having enough injectives
is defined in a dual way. 

We have the following characterization of stably 2-CY categories.
\begin{proposition}\label{propI1.1}
  Let $\C$ be an exact Frobenius category. 
Then $\C$ is stably 2-CY if and only if $\Ext^1_{\C}(A,B)$ 
is finite dimensional and we have functorial isomorphisms 
$D\Ext^1_{\C}(A,B)\simeq\Ext^1_{\C}(B,A)$ for all $A$, $B$ in $\C$. 
\end{proposition}
\begin{proof}
  Let $A$ and $B$ be in $\C$, and let $0\to A\to P\to \Omega^{-1}A\to 0$ 
be an exact sequence in $\C$ where $P$ is projective injective. Apply 
$\Hom_{\C}(B,\text{ })$ to get the exact sequence $0\to 
\Hom_{\C}(B,A)\to \Hom_{\C}(B,P)\to \Hom_{\C}(B,\Omega^{-1}A)\to 
\Ext^1_{\C}(B,A)\to 0$. Then we get $\Ext^1_{\C}(B,A)\simeq
\Hom_{\underline{\C}}(B,\Omega^{-1}A)= \Ext^1_{\underline{\C}}(B,A)$. 

Assume that $\C$ is stably 2-CY, that is, the stable category $\ul{\C}$ 
is a Hom-finite triangulated 2-CY category. Then $\Ext^1_{\ul{\C}}(B,A)$ 
is finite dimensional for all $A$, $B$ in $\C$, and hence $\Ext^1_{\C}(B,A)$ 
is finite dimensional, and we have functorial isomorphisms $D\Ext^1_{\C}(A,B)\simeq\Ext^1_ {\C}(B,A)$.

The converse also follows directly.
\end{proof}

Examples \sloppy of exact stably 2-CY categories are categories of maximal Cohen-Macaulay modules $\CM(R)$ for a 
3-dimensional complete local commutative isolated Gorenstein singularity $R$ (containing the residue field $k$)
and $\mod\Lambda$ 
for $\Lambda$ being the preprojective algebra of a Dynkin quiver. We shall see several further examples later.

We are especially interested in pairs of 2-CY categories
$(\C,\underline{\C})$ where $\C$  is an exact stably 2-CY category. The only difference in
indecomposable objects between $\C$ and $\underline{\C}$ is the indecomposable projective objects in
$\C$. Also note that given an exact sequence $0\to A\to B\to C\to 0$ in $\C$,
there is an associated triangle $A\to B\to C\to A[1]$ in
$\underline{\C}$. Conversely, given a triangle $A\to B\stackrel{\ul{g}}{\to} C\to A[1]$
in $\underline{\C}$, we lift $\ul{g}\in\Hom_{\ul{\C}}(B,C)$ to $g\in\Hom_{\C}(B,C)$, and obtain an
exact sequence $0\to A\to B\oplus P\to C\to 0$ in $\C$, where $P$ is projective. We then have the 
following useful fact.

\begin{proposition}\label{prop1.1}
  Let $\C$ be an exact stably 2-CY category with a set of clusters
$\ul{x}$ and a set of coefficients $\ul{p}$ which are the
indecomposable projective objects. For the stable 2-CY category
$\ul{\C}$ consider the same set of clusters $\ul{x}$ and with no
coefficients. Then we have the following. 
\begin{itemize}
 \item[(a)]The $(\ul{x},\ul{p})$ give a weak cluster structure on $\C$ if
   and only if the $(\ul{x},\emptyset)$ give a weak cluster structure on
   $\ul{\C}$. 
\item[(b)]$\C$ has no loops if and only if $\ul{\C}$ has no loops.
  \item[(c)]If $\C$ has a cluster structure, then $\ul{\C}$ has a cluster structure.
  \end{itemize}
\end{proposition}

\begin{proof}
(a) Assume we have a collection of extended clusters for $\C$. 
Assume also that $\mu_{M}(T)$ is defined for each indecomposable non-projective object $M$ 
in the cluster $T$, and that we have the required exchange exact sequences. Then the 
induced clusters 
for $\underline{\C}$ determine a weak cluster structure for $\underline{\C}$. 
The converse also follows directly. 

\noindent
(b) For an extended cluster $T$ in $\C$, the quiver $\bar{Q}_{T}$ is obtained from the quiver $Q_{T}$  
by removing the vertices corresponding to the indecomposable projective objects. The claim is then obvious.

\noindent
(c) It is clear that if there is no 2-cycle for $\C$, then there is no
2-cycle for $\ul{\C}$, and the claim follows from this. 
\end{proof}

In the examples of 2-CY categories with cluster structure
which have been investigated the extended clusters have been the
subcategories $T$ where $\Ext^1(M,M)=0$ for all $M\in T$,
and whenever $X\in\C$ satisfies $\Ext^1(M,X)=0$ for all $M\in T$,
then $X\in T$. Such $T$ has been called {\it cluster tilting
  subcategory} in \cite{bmrrt,kr1} if it is 
in addition functorially finite in the sense of \cite{as}, which is
automatically true when $T$ is finite. Such $T$
has been called maximal 1-orthogonal subcategory in \cite{i1,i2}, and
$\Ext$-configuration in \cite{bmrrt}, without the assumption on
functorially finiteness. 

We have the following nice connections between $\C$ and 
$\underline{\C}$ for an exact stably 2-CY category $\C$ when using the 
cluster tilting subcategories.

\begin{lemma}\label{lem1.2}
Let $\C$ be an exact stably 2-CY category, and let $T$ be a
subcategory of $\C$ containing all indecomposable projective objects. 
Then $T$ is a cluster tilting subcategory in $\C$ if and only if 
it is the same in 
$\underline{\C}$.
\end{lemma}

\begin{proof}
We have $\Ext^1_{\C}(C,A) \simeq \Ext^1_{\ul{\C}}(C,A)$ from the proof of 
Proposition \ref{prop1.1}. It is easy to see that $T$ is functorially finite in $\C$ if and only 
if it is functorially finite in $\ul{\C}$ \cite{as}. Hence $T$ is cluster tilting in $\C$ if and 
only if it is cluster tilting in $\underline{\C}$.
\end{proof}

\begin{lemma}\label{lem1.3}
Let $\C$ be an exact stably 2-CY category and $T$ a cluster tilting object in $\C$, 
with an indecomposable non-projective summand $M$. Then there is no loop at $M$ for 
$\End_{\C}(T)$ if and only if  there is no loop at $M$ for 
$\End_{\underline{\C}}(T)$. If $\C$ has no 2-cycles, there are none for $\ul{\C}$.
\end{lemma}

\begin{proof}
This is a direct consequence of Proposition \ref{prop1.1} and the definitions.
\end{proof}

Note that an exact stably 2-CY category $\C$, with the cluster tilting 
subcategories, gives a situation where we have a natural set of coefficients, 
namely the indecomposable projective objects which clearly belong to all cluster tilting 
subcategories, whereas $\underline{\C}$ with the cluster tilting subcategories 
gives a case where it is natural to choose no coefficients. We have
the following useful observation, which follows from Proposition
\ref{prop1.1}. 

\begin{proposition}\label{prop1.4}
Let $\C$ be a $\Hom$-finite exact stably 2-CY category. Then the cluster tilting subcategories 
in $\C$, with the indecomposable projectives as coefficients,  determine a weak cluster structure 
on $\C$ if and only if the cluster tilting subcategories in $\underline{\C}$ 
determine a weak cluster structure on $\underline{\C}$.
\end{proposition}
When $\C$ is Hom-finite triangulated 2-CY, then $\C$ has a weak cluster structure, 
with the extended clusters being the cluster tilting subcategories and the indecomposable
projectives being the coefficients \cite{iy}. Properties 
(c) and (d) hold for cluster categories and the stable category $\ul{\mod}\la$ of a 
preprojective algebra of Dynkin type \cite{bmrrt,bmr2,gls1}, but (c) does not hold in 
general \cite{bikr}. However, we show that when we have some cluster
tilting object in the 2-CY category $\C$, then (d) 
holds under the assumption that (c) holds. This was first proved in \cite{p} when $\C$ is {\em algebraic}, that is
by definition the 
stable category of a Frobenius category, as a special case of a more general result. Our proof 
is inspired by \cite[7.1]{ir}.

\begin{theorem}\label{teoI1.6}
Let  $\C$ be $\Hom$-finite triangulated (or exact stably) 2-CY
category with some cluster tilting subcategory. If $\C$ has no loops or
2-cycles, then the cluster tilting subcategories determine a cluster structure for $\C$.
\end{theorem}

\begin{proof}
We give a proof for the triangulated 2-CY case. Using exact sequences instead of triangles, a similar argument works for
the stably 2-CY case.
Note that in the stably 2-CY case we do not have to consider arrows
between projective vertices.

Let $T=\oplus_{i=1}^nT_i$ be a cluster tilting subcategory in $\C$. Fix
a vertex $k\in\{1,\cdots ,n\}$,
and let $T^{\ast}=\oplus_{i\ne k}T_i\oplus T_k^{\ast}=\mu_k(T)$. 
We have exchange triangles $T_k^{\ast}\to B_k\to T_k$ and $T_k\to B_k'\to T_k^{\ast}$, showing that 
when passing from $\End(T)$ to $\End(T^{\ast})$ we reverse all arrows in the quiver of $\End(T)$ 
starting or ending at $k$.

We need to consider the situation where we have arrows $j \to k \to i$. Since there
are no 2-cycles, there is no arrow $i \to k$. Consider the exchange triangles
$T_i^{\ast} \to B_i \to T_i$ and $T_i \to B_i' \to T_i^{\ast}$. Then $T_k$ is not a 
direct summand of $B_i'$, and we write $B_i = D_i \oplus T_k^m$ for some $m>0$, where
$T_k$ is not a direct summand of $D_i$.

Starting with the maps in the upper square and the triangles they
induce, we get by the octahedral axiom the diagram below, where the
third row is a triangle.
$$\xymatrix{            
&(T_k^{\ast})^m[1] \ar@{=}[r]&   (T_k^{\ast})^m[1]& \\
T_i[-1] \ar[r]\ar@{=}[d] &   T_i^{\ast}\ar[r]\ar[u] & D_i\oplus T_k^m \ar[u]\ar[r]& T_i\ar@{=}[d]\\
T_i[-1]\ar[r] &  X \ar[r]\ar[u] &D_i\oplus B_k^m \ar[r]\ar[u]& T_i\\            
& (T_k^{\ast})^m \ar@{=}[r] \ar[u] &    (T_k^{\ast})^m \ar[u] &
}$$

Using again the octahedral axiom, we get the following commutative diagram of triangles, 
where the second row is an exchange triangle and the third column is the second column of the 
previous diagram.
$$\xymatrix{ 
& (T_k^{\ast})^m[1] \ar@{=}[r] &    (T_k^{\ast})^m[1] & \\
T_i\ar@{=}[d]\ar[r]   &    B_i' \ar[r]\ar[u] &   T_i^{\ast} \ar[r]\ar[u]& T_i[1]\ar@{=}[d]\\
T_i\ar[r] & Y \ar[r]\ar[u]&      X\ar[r]\ar[u]& T_i[1]\\
    &        (T_k^{\ast})^m\ar@{=}[r] \ar[u]& (T_k^{\ast})^m \ar[u]&
}$$
Since $T_k$ is not in $\add B_i'$, we have $(B_i', (T_k^{\ast})^m[1])=0$, and hence $Y=B_i'\oplus(T_k^{\ast})^m$.

Consider the triangle $X\to D_i\oplus B_k^{m}\xto{a}T_i \to X[1]$.
Let $\overline{T}^{\ast} = (\oplus_{t \neq i,k} T_t) \oplus T_k^{\ast}$.
We observe that $D_i \oplus B_k^m$ is in $\add \overline{T}^{\ast}$.
For $D_i \oplus B_k^m = B_i$ is in $\add T$.
Since there is no loop at $i$, then $T_i$ is not a direct summand of $D_i$,
and $T_i$ is not a direct summand of $B_k$ since there is no arrow from $i$ to $k$.
Further $T_k$ is not a direct summand of $D_i$ by the choice of $D_i$,
and $T_k$ is not a direct summand of $B_k$ since there is no loop at $k$. Hence
we see that $B_i$ is in $\add  \overline{T}^{\ast}$.

We next want to show that $a$ is a right $\add \overline{T}^{\ast}$-approximation.
It follows from the first commutative diagram that any map $g \colon T_t \to T_i$,
where $T_t$ is an indecomposable direct summand of $ \overline{T}^{\ast}$
not isomorphic to $T_k^{\ast}$, factors through $a$. Let then $f \colon T_k^{\ast} \to T_i$
be a map, and $h \colon T_k^{\ast} \to B_k$ the minimal left $\add \overline{T}$-approximation, where 
$\overline{T} = \oplus_{t \neq k} T_t$. Then
there is some $s \colon B_k \to T_i$ such that $hs=f$. Then $s$ factors through $a$ by
the above, since $B_k$ is in $\add \overline{T}^{\ast}$ (using that $T_i$ is not a direct summand in $B_k$), and
$T_k^{\ast}$ is not a direct summand of $B_k$. It follows that $a$ is  
a right $\add \overline{T}^{\ast}$-approximation.

Consider now the triangle $T_i \to B_i' \oplus (T_k^{\ast})^m \xto{b} X \to T_i[1]$.
Then $ B_i' \oplus (T_k^{\ast})^m$ is clearly in $\add \overline{T}^{\ast}$, since
$T_k$ is not a direct summand of $B_i'$. Since $T_i$ is in both $T$ and $T^{\ast}$,
we have that $\Hom(\overline{T}^{\ast},T_i[1]) = 0$, and hence $b$ is a right 
$\add \overline{T}^{\ast}$-approximation. Note that the approximations $a$ and $b$
need not be minimal.

Recall that we are interested in paths of length two $j \to k \to i$
passing through $k$. By the above, the number of arrows from $j$ to $i$
in the quiver $Q_{T^{\ast}}$ is   
$$u = \alpha_{D_i \oplus B_k^m}(T_j) - \alpha_{B_i' \oplus (T_k^{\ast})^m} (T_j)$$
where $\alpha_X(T_j)$ denotes  the multiplicity of $T_j$ in $X$. 
We have $$u =\alpha_{D_i}(T_j) + m \alpha_{B_k}(T_j) - \alpha_{B_i'}(T_j) = 
\alpha_{B_i}(T_j) + m \alpha_{B_k}(T_j) - \alpha_{B_i'}(T_j),$$
since $B_i = D_i \oplus T_k^m$.
The last expression says that $u$ is equal to the number of arrows from $j$ to $i$
in $Q_T$, minus the number of arrows from $i$ to $j$, plus the product of the number
of arrows from $j$ to $k$ and from $k$ to $i$. This is what is required for having the 
Fomin-Zelevinsky mutation, and we are done.
\end{proof}

We shall also use the terminology stably 2-CY for certain subcategories of triangulated 
categories. Let $\B$ be a functorially finite extension closed subcategory of a Hom-finite triangulated 
2-CY category $\C$. We say that $X\in \B$ is {\it projective} in $\B$ if $\Hom(X,\B[1])=0$.
In this setting we shall prove in \ref{teoI2.1} that the category 
$\B$ modulo projectives in $\B$ has a 2-CY triangulated structure.
Note that $\B$ does not necessarily have enough projectives or injectives, for example if $\B = \C$.
We then say that $\B$ is \emph{stably 2-CY}. 

We illustrate the concept with the following.

\bigskip
\noindent
\textbf{Example.}
Let $\C_Q$ be the cluster category of the path algebra $kQ$, where $Q$ is the quiver 
$\stackrel{1}{\cdot}\to\stackrel{2}{\cdot}\to\stackrel{3}{\cdot}$. We have the following 
AR-quiver for $\C_Q$, where $S_i$ and $P_i$ denote the simple and projective modules 
associated with vertex $i$ respectively.
$$
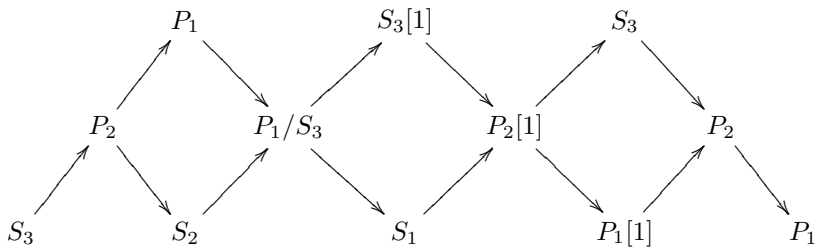
\xymatrix@C0.5cm{
&&P_1\ar[dr]&& S_3[1]\ar[dr]&&S_3\ar[dr]&&\\
&P_2\ar[ur]\ar[dr]&& P_1/S_3\ar[ur]\ar[dr]&& P_2[1]\ar[ur]\ar[dr]&& P_2\ar[dr]& \\
S_3\ar[ur]&& S_2\ar[ur]&& S_1 \ar[ur]&& P_1[1]\ar[ur] && P_1
}$$

Then $\B=\mod kQ$  is an extension closed subcategory of $\C_Q$ and it is easy to 
see that $P_1$ is the only indecomposable projective object in $\B$. Then $\B/P_1$ is clearly 
equivalent to the cluster category $\C_{Q'}$ where $Q'$ is a quiver of type $A_2$, which is a  
triangulated 2-CY category. Hence $\B$ is stably 2-CY.

\bigskip
In addition to the cluster tilting objects, also the maximal rigid objects have played an important 
role in the investigation of 2-CY categories. We now investigate the concepts of cluster structure 
and weak cluster structure with respect to these objects.

Recall that a subcategory $T$ of a category $\C$ is said to
be {\em rigid} if $\Ext^1(M,M)=0$ for all $M$ in $T$, 
and {\em maximal rigid} if $T$ is maximal among rigid subcategories
\cite{gls1}. It is clear that any cluster tilting subcategory
is maximal rigid, but the converse is not the case \cite{bikr}. 
There always exists a maximal rigid subcategory in $\C$ if the
category $\C$ is skeletally small, while the existence of a cluster
tilting subcategory is rather restrictive. It is of interest
to have sufficient conditions for the two concepts to coincide.
For this the following is useful (see \cite{bmr1,i1,kr1} for (a) and
the argument in \cite[5.2]{gls1} for (b)).

\begin{proposition}\label{pro-extra}
Let $\C$ be a triangulated (or exact stably) 2-CY category.
\begin{itemize}
\item[(a)] Let $T$ be a cluster tilting subcategory. Then for any $X$
  in $\C$, there exist triangles (or short exact sequences) $T_1 \to
  T_0 \to X$ and $X \to T_0' \to T_1'$ with $T_i, T_i'$ in $T$.
\item[(b)] Let $T$ be a functorially finite maximal rigid subcategory.
Then for any $X$ in $\C$ which is rigid, the same conclusion as in (a) holds.
\end{itemize}
\end{proposition}

Then we have the following.

\begin{theorem}\label{prop2.1}
  Let $\C$ be an exact stably 2-CY category, with some cluster tilting object. 
  \begin{itemize}
    \item[(a)]Then any maximal rigid object in $\C$ (respectively, $\ul{\C}$) is a cluster tilting object.
    \item[(b)]Any rigid subcategory in $\C$ (respectively, $\ul{\C}$) has an additive generator which 
      is a direct summand of a cluster tilting object.
    \item[(c)]All cluster tilting objects in $\C$ (respectively,
      $\ul{\C}$) have the same number of non-isomorphic indecomposable
      summands.
  \end{itemize}
\end{theorem}

\begin{proof}
(a) Let $N$ be maximal rigid in $\C$.
We only have to show that any $X\in\C$ satisfying $\Ext^1(N,X)=0$ is contained in $\add N$.

(i) Let $M$ be a cluster tilting object in $\C$. Since $N$ is maximal rigid and $M$ is rigid, 
there exists an exact sequence $0\to N_1\to N_0\to M\to 0$ with $N_i \in \add N$
by Proposition \ref{pro-extra}(b). In particular, we have ${\rm pd }_{{\rm End}(N)}{\rm
    Hom}(N,M)\le1$.

(ii) Since $M$ is cluster tilting, there is, by Proposition 
\ref{pro-extra}(a), an exact sequence $0 \to X \to M_0 \to M_1\to 0$ for $X$ as above, with
$M_i \in \add M$, obtained by taking the minimal left $\add M$-approximation $X \to M_0$.
Applying $(N,\ )$, 
we have an exact sequence $0\to
(N,X) \to (N,M_0) \to (N,M_1) \to \Ext^1(N,X)=0$. By (i), we have $\pd_{\End(N)}\Hom (N,X) \le1$. 
Take a projective resolution
$0 \to (N,N_1) \to (N,N_0) \to (N,X) \to 0$. 
Then we have a complex
\begin{equation}\label{N approximation of X}
0\to N_1\to N_0\to X\to0
\end{equation}
in $\C$. Since $0 \to (P,N_1) \to (P,N_0) \to (P,X)
\to 0$ is exact for any projective $P$ in $\C$, it follows from the axioms
of Frobenius categories that the complex \eqref{N approximation of X} is an exact sequence in
$\C$. Since $\Ext^1(X,N)=0$, we have
$X\in \add N$, and hence $N$ is cluster tilting. 

(b)  Let $M$ be a cluster tilting object in $\C$ and $N$ a
rigid object in $\C$. By \cite[5.3.1]{i2}, $\Hom(M,N)$ is a partial
tilting $\End(M)$-module. In particular, the number of non-isomorphic
indecomposable direct summands of $N$ is not greater than that of $M$.
Consequently, any rigid object in $\C$ is a direct summand of some
maximal rigid object in $\C$, which is cluster tilting by (a).

(c) See \cite[5.3.3]{i2}.
\end{proof}

For a Hom-finite triangulated  2-CY category we also get a weak cluster structure, 
and sometimes a cluster structure, determined by the maximal rigid objects, if there are any. 
Note that there are cases where the maximal rigid objects are not cluster tilting \cite{bikr}. 
But we suspect the following.

\begin{conjecture}
Let $\C$ be a connected $\Hom$-finite triangulated 2-CY category. Then any maximal rigid object without 
loops or 2-cycles in its quiver is a cluster tilting object.
\end{conjecture}

Furthermore, we have the following.

\begin{theorem}\label{theoI1.8}
  Let $\C$ be a $\Hom$-finite triangulated 2-CY category (or exact stably
  2-CY category) having some functorially finite maximal rigid subcategory.
  \begin{itemize}
    \item[(a)]{The functorially finite maximal rigid subcategories determine a weak cluster structure on $\C$.}
    \item[(b)]{If there are no loops or 2-cycles for the functorially finite maximal rigid subcategories, 
then they determine a cluster structure on $\C$.}
\end{itemize}
\end{theorem}
\begin{proof}
(a) This follows from \cite[5.1,5.3]{iy}. Note that the arguments there are stated only
for cluster tilting subcategories, but work also for functorially finite maximal rigid subcategories.

\noindent
(b) The proof of Theorem \ref{teoI1.6} works also in this setting.
\end{proof}

There exist triangulated or exact categories with cluster tilting objects also when the 
categories are not 2-CY or stably 2-CY (see \cite{i1,kz,eh}). But we do not necessarily 
have even a weak cluster structure in this case. For let $\la$ be a Nakayama algebra with  
two simple modules $S_1$ and $S_2$, with associated projective covers $P_1$ and $P_2$. 
Assume first that $P_1$ and $P_2$ have length 3. Then in $\mod\la$ we have that $S_1\oplus P_1\oplus P_2$, 
$\begin{smallmatrix} S_1\\ S_2 \end{smallmatrix}\oplus P_1\oplus P_2$, 
$\begin{smallmatrix} S_2\\ S_1 \end{smallmatrix}\oplus P_1\oplus P_2$ are  the cluster 
tilting  objects, so we do not have the unique exchange property.

If $P_1$ and $P_2$ have length 4, then the cluster tilting objects are 
$S_1\oplus \begin{smallmatrix} S_1\\ S_2 \\S_1 \end{smallmatrix}\oplus P_1\oplus P_2$ and 
$S_2\oplus \begin{smallmatrix} S_2\\ S_1\\S_2 \end{smallmatrix}\oplus P_1\oplus P_2$, 
and so there is no way of exchanging $S_1$ in the first object to obtain a new cluster tilting object.

\bigskip
We end ths section with some information on the endomorphism algebras of cluster tilting
objects in stably 2-CY categories. Such algebras are studied as
analogs of Auslander algebras in \cite{gls1,i1,i2,kr1}.
We denote by $\mod\C$ the category of finitely presented $\C$-modules.
If $\C$ has pseudokernels, then $\mod\C$ forms an abelian category \cite{aus1}.

\begin{proposition}\label{prop2.5}
Let ${\C}$ be an exact stably 2-CY category.
Assume that $\C$ has pseudokernels and the global dimension of ${\rm mod}\C$ is finite.
Let $\Gamma=\End(T)$ for a cluster tilting object $T$ in $\C$.
\begin{itemize}
\item[(a)]{$\Gamma$ has finite global dimension.} 
\item[(b)]{If ${\C}$ is $\Hom$-finite, then the quiver of $\Gamma$ has no loops.
If moreover ${\C}$ is an extension closed subcategory of an abelian
category closed under subobjects, then the quiver of $\Gamma$ has no
2-cycles.} 
\end{itemize}
\end{proposition}

\begin{proof}
(a) Let $m={\rm gl.dim} (\mod\C)$. 
For any $X\in \mod \Gamma$, take a projective
presentation $(T,T_1)\to(T,T_0)\to X\to0$. By our assumptions, there
exists a complex $0\to F_m\to\cdots\to F_2\to T_1\to T_0$ in
$\C$ such that
$0\to(\ ,F_m)\to\cdots\to(\ ,F_2)\to(\ ,T_1)\to(\ ,T_0)$ is exact in
$\mod \C$. Since $T$ is cluster tilting, 
we have an exact sequence $0\to T_1\to T_0\to F_i\to 0$, with 
$T_1$ and $T_0$ in $\add T$ by Proposition \ref{pro-extra}. Hence we have ${\pd}_{{\Gamma}}(T,F_i)\le 1$ and 
consequently ${\pd}_{\Gamma}X\le m+1$. It follows that $\Gamma$ has finite global dimension.

(b) By (a), $\Gamma$ is a finite dimensional algebra of finite global dimension.
By \cite{l,ig}, the quiver of $\Gamma$ has no loops.

We shall show the second assertion.
Our proof is based on \cite[6.4]{gls1}. 
We start with showing that $\Ext^2_\Gamma(S,S)=0$ for any simple
$\Gamma$-module $S$, assumed to be the top of the projective
$\Gamma$-module $(T,M)$ for an indecomposable summand $M$ of $T$.

First, we assume that $M$ is not projective in $\C$.
Take exact exchange sequences $0\to M^{\ast} \stackrel{f}{\to}B\stackrel{g}{\to}M\to0$
and $0\to M\stackrel{s}{\to}B'\stackrel{t}{\to}M^{\ast}\to0$.
Since $\Gamma$ has no loops, we have a projective presentation
$0\to(T,M)\stackrel{\cdot s}{\to}(T,B')\stackrel{\cdot tf}{\to}(T,B)
\stackrel{\cdot g}{\to}(T,M)\to S\to0$
of the $\Gamma$-module $S$. 
Since $M$ is not a summand of $B'$, we have
${\Ext}^2_\Gamma(S,S) = 0$.

Next, we assume that $M$ is projective in $\C$.
Take a minimal projective presentation $(T,B)\stackrel{\cdot g}{\to}(T,M)\to S\to0$
of the $\Gamma$-module $S$. By assumption, $\Im g$ in the abelian category belongs to $\C$. 
Then $g \colon B\to{\Im} g$ is a minimal right ${\add} T$-approximation.
By Proposition \ref{pro-extra}(a), we have that $B'={\Ker} g$ belongs to ${\add} T$.
Thus we have a projective resolution
$0\to(T,B')\stackrel{}{\to}(T,B)\stackrel{\cdot g}{\to}(T,M)\to S\to0$
of the $\Gamma$-module $S$.
Since $g$ is right minimal, $B'$ does not have an injective summand.
Thus we have ${\Ext}^2_\Gamma(S,S) = 0$.

Since in both cases $\Ext^2_\Gamma(S,S)=0$, we can not have a 2-cycle
by \cite[3.11]{gls1}.
\end{proof}
\hspace{7mm}

\subsection{Substructures}\label{c1_sec2} 
${}$ \\
For extension closed subcategories of triangulated or exact categories both having a weak 
cluster structure, we introduce the notion of substructure. Using heavily \cite{iy}, 
we give sufficient conditions for having a substructure, when starting with a Hom-finite triangulated 
2-CY category or an exact stably 2-CY category, and using the cluster tilting subcategories,
with the indecomposable projectives as coefficients.

Let $\C$ be an exact or triangulated $k$-category,
and $\B$ a subcategory of $\C$
closed under extensions. Assume that both $\C$ and
$\B$  have a weak cluster structure. We say that we have
a \emph{substructure} of $\C$ induced by an extended cluster $T$ in
$\B$ if we have the following:

There is
a set $A$ of indecomposable objects in $\C$ such that
$\widetilde{T'} = T'\cup A$ is an extended cluster in $\C$
for any extended cluster $T'$ in $\B$ which is obtained
by a finite number of exchanges from $T$.

Note that for each  sequence of cluster variables $M_1,\cdots, M_t$, with 
$M_{i+1}$ in $\mu_{M_i}(T)$, we have
$\mu_{M_t}(\cdots\mu_{M_1}(T))\cup A=\tilde{\mu}_{M_t}(\cdots
\tilde{\mu}_{M_1}(\widetilde{T}))$, where $\mu$ denotes
the exchange for $\B$ and $\tilde{\mu}$ the exchange for $\C$.

We shall investigate substructures arising from certain extension closed subcategories of 
triangulated 2-CY categories and of exact stably 2-CY categories. We start with the 
triangulated case, and here we first recall some results from \cite{iy} specialized to the 
setting of 2-CY categories. 

For a triangulated category $\C$ and full subcategories $\B$ and $\B'$,
let ${\B}^\perp =\{X \in \C  \mid \Hom({\B},X)=0\}$ and
${}^\perp{\B}=\{ X\in \C \mid \Hom(X,\B)=0\}$.
We denote by $\B \ast {\B}'$ the full subcategory of $\C$
consisting of all $X \in \C$ such that there exists a triangle
$B\to X\to B'\to B[1]$ with $B\in \B$ and $B'\in {\B}'$.

We get the following sufficient 
conditions for constructing 2-CY categories, and hence categories with weak cluster structures.

\begin{theorem}\label{teoI2.1}
Let ${\C}$ be a $\Hom$-finite triangulated 2-CY category
and ${\B}$ a functorially finite extension closed subcategory of $\C$.
\begin{itemize}
\item[(a)]{${\B}^\perp$ and ${}^\perp{\B}$ are functorially finite
extension closed subcategories of $\C$.
Moreover, ${\B}\ast {\B}^\perp={\C}={}^\perp{\B}\ast{\B}$
and ${}^\perp({\B}^\perp)={\B}=({}^\perp{\B})^\perp$ hold.}

\item[(b)]{Let $\D= \B \cap {}^\perp{\B}[1]$. 
Then $\B/\D$ is a $\Hom$-finite triangulated 2-CY 
category, so that $\B$ is a stably 2-CY category.
Moreover, ${\B}\subseteq({\D}\ast{\B}[1]) \cap ({\B}[-1]\ast{\D})$ holds,
and $\D$ is a functorially finite rigid subcategory of $\C$.}

\item[(c)]{Let $\D$ be a functorially finite rigid subcategory
of $\C$ and ${\B'}={}^\perp{\D}[1]$. Then ${\B'}$ is a
functorially finite extension closed subcategory of ${\C}$
and ${\B'}/{\D}$ is a triangulated 2-CY category.
Moreover, there exists a one-one correspondence between cluster tilting
(respectively, maximal rigid, rigid) subcategories of ${\C}$ containing ${\D}$ and
cluster tilting (respectively, maximal rigid, rigid) subcategories of
${\B'}/{\D}$. It is given by $T\mapsto T/\D$.}
\end{itemize}
\end{theorem}

\begin{proof}
(a) Since $\B^\perp={}^\perp\B[2]$ holds by the 2-CY property, the
assertion follows from \cite[2.3]{iy}.

(b) Clearly $\B/\D$ is $\Hom$-finite, since $\C$ is.
To show that $\B/\D$ is a triangulated 2-CY category, we only need to check
${\B}\subseteq({\D}\ast{\B}[1]) \cap ({\B}[-1]\ast{\D})$ by \cite[4.2]{iy}.
Let $Z$ be in $\B$. Since $\B$, and hence $\B[1]$, is functorially finite in 
$\C$, it follows from (a) that we have a triangle $X\to Y\to Z\to X[1]$ with $Y$ in 
$^{\bot}\B[1]$ and $X[1]$ in $\B[1]$. Since $\B$ is extension closed, 
$Y$ is in $\B$, and consequently $Y$ is in $\B\cap{}^{\bot}\B[1]=\D$. 
It follows that $Z$ is in $\D{\ast}\B[1]$, and similarly in $\B[-1] \ast \D$.

To see that $\D$ is functorially finite in $\C$,
we only have to show that $\D$ is functorially finite in $\B$.
For any $Z \in \B$, take the above triangle $X\to Y\stackrel{f}{\to}Z\to X[1]$
with $Y$ in $\D$ and $X[1]$ in $\B[1]$.
Since $(\D,X[1])=0$, we have that $f$ is a right $\D$-approximation.
Thus $\D$ is contravariantly finite in $\B$,
and similarly covariantly finite in $\B$.

(c) See \cite[4.9]{iy}.
\end{proof}

The example of the cluster category $\C$ of the path algebra $kQ$ where $Q$ is of type 
$A_3$ from the previous section illustrates part of this theorem. For let $\D=\add P_1$. 
Then $\B'={^{\bot}\D[1]}=\mod kQ$, and $\B'/\D=\C_{kQ'}$, where $Q'$ is a quiver of type $A_2$. 
The cluster tilting objects in $\C$ containing $P_1$ are $P_1\oplus S_3\oplus P_2$, $P_1
\oplus P_2\oplus S_2$, $P_1\oplus S_2\oplus P_1/S_3$, $P_1\oplus P_1/S_3\oplus S_1$, 
$P_1\oplus S_1\oplus S_3$, which are in one-one correspondence with the cluster tilting objects in $\B'/\D$.

In order to get sufficient conditions for having a substructure we investigate cluster tilting 
subcategories in $\B$. For this the following lemma is useful.

\begin{lemma}\label{lemI2.3}
Let $\C$ be a $\Hom$-finite triangulated 2-CY category. For any
functorially finite and thick subcategory $\C_1$ of
$\C$, there exists a functorially finite and thick
subcategory $\C_2$ of $\C$ such that
$\C=\C_1\times\C_2$.
\end{lemma}

\begin{proof}
Let $\C_2 =\C_1^\perp$.
Then we have
$\C_2=\C_1^\perp={}^\perp\C_1[2]={}^\perp\C_1$
by Serre duality, using that $\C_1$ is triangulated. We only have to show that any object in
$\C$ is a direct sum of objects in $\C_1$ and
$\C_2$. For any $X\in\C$, there exists a triangle
$A_1\to X \to A_2\stackrel{f}{\to}A_1[1]$ in $\C$ with $A_1$
in $\C_1$ and $A_2$ in $\C_2=\C_1^{\bot}$ by 
Theorem \ref{teoI2.1}(a). Since $f=0$, we have $X\simeq A_1\oplus A_2$. 
Thus we have $\C=\C_1\times\C_2$.
\end{proof}

Using Lemma \ref{lemI2.3}, we get the following decomposition of
triangulated categories.

\begin{proposition}\label{propIadded}
Let $\C$ be a $\Hom$-finite triangulated 2-CY category and $\B$ a
functorially finite extension closed subcategory of $\C$. Let $\D= \B
\cap {}^\perp{\B}[1]$ and $\B'={}^\perp\D[1]$. 
\begin{itemize}
\item[(a)]{There exists a
functorially finite and extension closed subcategory $\B''$ of $\C$
such that $\D\subseteq\B''\subseteq\B'$ and $\B'/\D=\B/\D\times\B''/\D$ as
a triangulated category.}

\item[(b)]{There exists a one-one correspondence between pairs
    consisting of cluster tilting (respectively, maximal rigid,
    rigid) subcategories of $\B$ and of $\B''$, and cluster tilting
    (respectively, maximal rigid, rigid) subcategories of $\B'$.
It is given by $(T,T'')\mapsto T\oplus T''$.}
\end{itemize}
\end{proposition}

\begin{proof}
(a) We know by Theorem \ref{teoI2.1}(b)(c) that 
$\D$ is functorially finite rigid, and that $\B/\D$ 
and $\B'/\D$ are both triangulated 2-CY categories. The inclusion functor 
$\B/\D\to \B'/\D$ is a triangle functor by the construction 
of their triangulated structures in \cite[4.2]{iy}. In particular $\B/\D$ is a thick subcategory 
of $\B'/\D$, and hence we have a decomposition by Lemma \ref{lemI2.3}. 

(b) This follows by Theorem \ref{teoI2.1}(c).
\end{proof}

Then we get the following.
\begin{corollary}\label{corI2.4}
Let $\C$ be a $\Hom$-finite 2-CY algebraic triangulated category with a cluster tilting object, 
and $\B$ a functorially finite extension closed subcategory of $\C$. Then we have the following.
\begin{itemize}
\item[(a)]{The stably 2-CY category $\B$ also has some cluster tilting
    object. Any maximal rigid object in $\B$ is a cluster tilting
    object in $\B$.}

\item[(b)]{There is some rigid object $A$ in $\C$ such that $T\oplus A$ is a cluster tilting 
object in $\C$ for any cluster tilting object $T$ in $\B$.}

\item[(c)]{Any cluster tilting object $T$ in $\B$ determines a substructure for the  
weak cluster structures on $\B$ and $\C$ given by cluster tilting objects.}
\end{itemize}
\end{corollary}

\begin{proof}
(a) Let $\D ={}^\perp\B[1]$ and $\B'={}^{\bot}\D[1]$.
Since $\C$ is algebraic by Theorem \ref{prop2.1}, we have
a cluster tilting object $T$ in $\C$ containing $\D$.
By Proposition \ref{propIadded}, we have decompositions 
$\B'/\D=\B/\D\times\B''/\D$ for some subcategory $\B''$ of $\B'$ 
and $T=T_1\oplus T_2$ with a cluster
tilting object $T_1$ (respectively, $T_2$) in $\B$ (respectively, $\B''$).
Thus $\B$ has a cluster tilting object.

Now we show the second assertion. 
Let $M$ be maximal rigid in $\B$. By Proposition \ref{propIadded}(b),
we have that $M\oplus T_2$ is maximal rigid in $\C$. By Theorem  \ref{prop2.1},
it follows that $M\oplus T_2$ is cluster tilting in $\C$ and by Proposition
\ref{propIadded}(b), we have that $M$ is cluster tilting in $\B$.

\noindent
(b) We only have to let $A = T_2$.

\noindent
(c) This follows from (b).
\end{proof}

It is curious to note that combining Proposition \ref{propIadded} with Theorem \ref{teoI2.1} we obtain a 
kind of classification of functorially finite extension closed subcategories of a  
triangulated 2-CY category in terms of functorially finite rigid
subcategories, analogous to results from \cite{ar}.

\begin{theorem}\label{theoI2.5}
Let $\C$ be a 2-CY triangulated category. Then the functorially finite  
extension closed subcategories $\B$ of $\C$ are all obtained as preimages under the 
functor $\pi\colon \C\to \C/\D$ of the direct summands of $^{\bot}\D[1]/\D$ as a triangulated category, 
for functorially finite rigid subcategories $\D$ of $\C$.
\end{theorem}

\begin{proof}
Let $\D$ be functorially finite rigid in $\C$. Then 
$\B'={^{\bot}\D}[1]$ is functorially finite extension closed in $\C$ by Theorem \ref{teoI2.1}(a). 
Then the preimage under $\pi\colon\C\to\C/\D$ of any direct summand of $\B'/\D$ as a triangulated category
is functorially 
finite and extension closed in $\C$.

Conversely, let $\B$ be a functorially finite extension closed subcategory of
$\C$ and $\D=\B\cap{}^{\bot}\B[1]$. By Proposition \ref{propIadded}, we have that $\B/\D$ is a direct
summand of ${^{\bot}}\D[1]/\D$.
\end{proof}

We now investigate substructures also for exact categories which are stably 2-CY. 
We have the following main result.

\begin{theorem}\label{teoI2.7}
Let $\C$ be an exact stably 2-CY category, and
$\B$  a functorially finite extension closed subcategory of $\C$. 
Then $\B$ has enough projectives and injectives, and is a stably 2-CY category.
\end{theorem}

\begin{proof}
We know that
$\B$ is an exact category and $\D=\B\cap{}^\perp\B[1]$ is the
subcategory of projective injective objects. Since $\B\subseteq\B[-1]\ast\D$
holds by Theorem \ref{teoI2.1}(b), then for any $X\in\B$, there exists
a triangle $X\to Y\to Z\to X[1]$ with $Y\in\D$
and $Z\in\B$. This is induced from an exact sequence $0\to X\to Y\to
Z\to0$ in $\C$. Thus $\B$ has enough injectives. Dually, $\B$ has
enough projectives, which coincide with the injectives. Hence $\B$ is a Frobenius category, and
consequently, $\B$ is stably 2-CY.
\end{proof}

Alternatively we give a direct approach, where the essential information is given by the following lemma and its dual.

\begin{lemma}\label{lemI2.3new}
Let $\C$ be an exact category with enough injectives, and
$\B$ a contravariantly finite extension closed subcategory of $\C$. Then $\B$ is an exact
category with enough injectives.
\end{lemma}

\begin{proof}
It is clear that $\B$ is also an exact category. Let $X$ be in $\C$ and
take an exact sequence $0 \to X \to I \to X' \to 0$ with $I$ injective in $\C$. Then we have an exact 
sequence of functors $(\ ,X') \to \Ext ^1(\ ,X)\to 0$. Since $\B$ is Krull-Schmidt and
contravariantly finite in $\C$, we can take a projective cover 
$\phi:(\ ,Y)\to{\Ext}^1(\ ,X)|_{\B}\to0$ of $\B$-modules.
This is induced by an exact sequence $0\to X\to Z\stackrel{}{\to}
Y\to0$ with terms in $\B$.

We will show that $Z$ is injective. Take any exact sequence $0\to
Z\stackrel{}{\to}Z'\to Z''\to0$ with terms in $\B$. We will
show that this splits. Consider the following exact commutative diagram:
\begin{equation}\label{commutative diagram}
\begin{array}{ccccccccc}
&&&&0&&0&&\\
&&&&\downarrow&&\downarrow&&\\
0&\to&X&\to&Z&\to&Y&\to&0\\
&&\parallel&&\downarrow&&\downarrow^{a}&&\\
0&\to&X&\to&Z'&\to&Y'&\to&0\\
&&&&\downarrow&&\downarrow&&\\
&&&&Z''&=&Z''&&\\
&&&&\downarrow&&\downarrow&&\\
&&&&0&&0&&
\end{array}
\end{equation}
Then $Y'\in\B$, and we have the commutative diagram
\begin{equation}\label{commutative diagram2}
\begin{array}{ccccccccccc}
0&\to&(\ ,X)&\to&(\ ,Z)&\to&(\ ,Y)&\stackrel{\phi}{\to}&{\rm Ext}^1(\
,X)|_{\B}&\to&0\\
&&\parallel&&\downarrow^{}&&\downarrow^{\cdot a}&&\parallel&&\\
0&\to&(\ ,X)&\to&(\ ,Z')&\to&(\ ,Y')&\to&{\rm Ext}^1(\ ,X)|_{\B}&&
\end{array}
\end{equation}
of exact sequences of $\B$-modules. Since $\phi$ is a projective cover, we have that
$(\cdot a)$ is a split monomorphism. Thus $a$ is a split monomorphism.
We see that the sequence $0 \to \Ext^1(Z'',Z) \to \Ext^1(Z'',Y)$ is
exact by evaluating the upper sequence in \eqref{commutative diagram2}
at $Z''$. Since the right vertical sequence 
in \eqref{commutative diagram} splits, it follows that the middle
vertical sequence in \eqref{commutative diagram} splits. Hence $Z$ is
injective, and consequently $\B$ has enough injectives.
\end{proof}

It follows from $\C$ being 2-CY that the projectives and injectives in $\B$ coincide,
and hence $\B$ is Frobenius by Lemma \ref{lemI2.3new} and its dual. It follows as before
that $\B$ is stably 2-CY, and the alternative
proof of Theorem \ref{teoI2.7} is completed.

\bigskip
We have the following interesting special case, as a consequence of Theorem \ref{teoI2.7} and Corollary \ref{corI2.4}.
For $X$ in $\C$ we denote by $\Sub X$ the subcategory of $\C$ whose objects are subobjects of finite direct sums of copies of $X$.
\begin{corollary}\label{cor2.6}
Let $\C$ be a $\Hom$-finite abelian stably 2-CY category, and let $X$ be an object in
$\C$  with $\Ext^1(X,X)=0$, and $\id X\leq 1$.
\begin{itemize}
\item[(a)]{Then $\Sub X$ is a functorially finite extension closed subcategory of $\C$ and is exact stably 2-CY.}
\item[(b)]{If $\C$ has a cluster tilting object, then so does $\Sub X$, and any cluster 
tilting object in $\Sub X$ determines a substructure of the cluster structure for $\C$.}
\item[(c)]{If $\C$ is abelian, then $\Sub X$ has no loops or 2-cycles.}
\end{itemize}
\end{corollary}

\begin{proof}
(a) We include the proof for the convenience of the reader.
We first want to show that $\Sub X$ is extension closed. Let $0\to A\to B\to C\to 0$ be an 
exact sequence with $A$ and $C$ in $\Sub X$, and consider the diagram
$$\xymatrix{
    & 0\ar[d]&& 0\ar[d] &\\
    0\ar[r]& A\ar[r]^i\ar[d]^f& B\ar[r]^j& C\ar[r]\ar[d]^q& 0\\
    0\ar[r]& X_0\ar[r]& X_0\oplus X_1\ar[r]& X_1\ar[r] & 0
}$$
with $X_0, X_1$ in $\add X$. Since $\id X_0\le 1$, we have the exact sequence 
$\Ext^1(X_1,X_0)\to \Ext^1(C,X_0)\to 0$, which shows that $\Ext^1(C,X_0)=0$. 
Then the exact sequence $(B,X_0)\to (A,X_0)\to \Ext^1(C,X_0)$ shows that there 
is some $t\colon B\to X_0$ such that $it=f$. This shows that $B$ is in $\Sub X$, 
which is then closed under extensions.
It is also functorially finite \cite{as}, and clearly
  Krull-Schmidt. So $\Sub X$ has enough projectives and injectives,
  with the projectives coinciding with the injectives. Hence $\Sub X $
  is Frobenius, and as we have seen before, it follows that the stable category
  $\underline{\Sub} X$ is 2-CY.\\
(b) This follows directly using Corollary \ref{corI2.4}.\\
(c) This follows from Proposition \ref{prop2.5}.
\end{proof}

In order to see when we have cluster structures we next want to give sufficient conditions 
for algebraic triangulated (or stably) 2-CY categories not to have loops or 2-cycles.

\begin{proposition}\label{propI2.11}
Let $\C$ be a $\Hom$-finite algebraic triangulated (or exact stably) 2-CY
category with a cluster tilting object, and $\B$ a functorially finite
extension closed subcategory. 
\begin{itemize}
\item[(a)]{If $\C$ has no 2-cycles, then also $\B$ has no 2-cycles.}
\item[(b)]{If $\C$ has no loops, then $\B$ has no loops.}
\end{itemize}
\end{proposition}

\begin{proof}
We give a proof for the algebraic triangulated 2-CY case. 
A similar argument works for the stably 2-CY case.

(a) Let $\D = \B \cap{}^\perp \B[1]$ and $\B' ={}^\perp \D[1]$.
Since cluster tilting objects in $\B'$ are exactly
cluster tilting objects in $\C$ which contain $\D$,
our assumption implies that $\B'$ has no 2-cycles.

We shall show that $\B$ has no 2-cycles.
Let $T$ be a cluster tilting object in $\B$.
By Corollary \ref{corI2.4}(b), there exists $T'\in \B'$ such that $T\oplus T'$ is
a cluster tilting object in $\C$.
We already observed that $T\oplus T'$ has no 2-cycles.
If $T$ has a 2-cycle, then at least one arrow in the 2-cycle
represents a morphism $f \colon X\to Y$ which factors through an object in $T'$.
We write $f$ as a composition of $f_1 \colon X\to Z$ and $f_2 \colon Z\to Y$ with $Z\in T'$.
Since $\B/\D$ is a direct summand of $\B'/\D$ by Proposition \ref{propIadded},
any morphism between $T$ and $T'$ factors through $\D$.
Thus we can write $f_1$ (respectively $f_2$) as a composition
of $g_1 \colon X\to W_1$ and $h_1 \colon W_1\to Z$ (respectively, $g_2 \colon Z\to W_2$ and 
$h_2 \colon W_2\to Y$)
with $W_1\in \D$ (respectively, $W_2\in \D$).
We have $f=f_1f_2=g_1(h_1g_2)h_2$, where $h_1g_2$ is in $\rad \B$ and
at least one of $h_2$ and $g_1$ is in $\rad \B$, since
at least one of $X$ and $Y$ is not in $\D$.
So $f$ can not be irreducible in $\add T$, a contradiction.

\noindent 
(b) This follows in a similar way.
\end{proof}

Note that the quiver $Q_T$ may have 2-cycles between coefficients.
For example, let $\C = \mod \la$ for the preprojective algebra of a Dynkin quiver and let $\B$ be the
subcategory $\add \la$. Then there are no 2-cycles for $\C$, but there are 2-cycles for $\B$,
since $\la$ is the only cluster tilting object in $\B$. \\
\hspace{7mm}

\subsection{Preprojective algebras of Dynkin type}\label{c1_sec3}
${}$ \\
In this section we specialize our general results from Section \ref{c1_sec2} to the case of the finitely generated 
modules over a preprojective algebra of Dynkin type. We also illustrate with three concrete examples. 
The same examples will be used in the Chapter \ref{chap3} to illustrate how 
to use this theory to construct subcluster algebras of cluster algebras.

The category $\C =\mod\Lambda$ for $\Lambda$ preprojective of Dynkin 
type is a Hom-finite Frobenius category. By \cite{gls1}, see also Section \ref{c2_sec2}, a rigid
$\la$-module is cluster tilting if and only if the number of
non-isomorphic indecomposable summands is the number of positive
roots, so $\frac{n(n+1)}{2}$ for $A_n$, $n(n-1)$ for $D_n$, $36$ for
$E_6$, $63$ for $E_7$ and $120$ for $E_8$.

Let $\B$ be an extension closed functorially finite 
subcategory of $\C$. We  know that $\B$ is stably 2-CY by Theorem \ref{teoI2.7}. It is  
known that $\C$ and $\ul{\C}$ have no loops or 2-cycles for the cluster 
tilting objects \cite{gls1}, and this also follows from Proposition
\ref{prop2.5}. Then it follows from Proposition \ref{propI2.11} that 
there are also no loops or 2-cycles for $\B$ and the subcategory $\ul{\B}$ of $\ul{\C}$.
Note that $\ul{\B}$ is not the stable category of $\B$ since $\B$ may
have more projectives than $\C$.

We then have the following.

\begin{theorem}\label{teoI3.1}
Let $\B$ be an extension closed functorially finite subcategory of the category $\C=\mod\la$ for the 
preprojective algebra $\la$ of a Dynkin quiver. Then we have the following.
\begin{itemize}
\item[(a)]The exact stably 2-CY category $\B$ has a cluster tilting object, and any 
maximal rigid object in $\B$ is a cluster tilting object,
which can be extended to a cluster tilting object for $\C$, and which gives rise to a substructure.
\item[(b)]{The category $\underline{\B}$ is a stably 2-CY Frobenius
    category with no loops or 2-cycles for the cluster tilting objects, and hence has a cluster structure}.
\end{itemize}
\end{theorem}

\begin{proof}
\noindent (a) This follows from Theorem \ref{prop2.1} and Corollary \ref{cor2.6}.

\noindent (b) This follows from the above comments and Theorem \ref{teoI1.6}.
\end{proof}

We now give some concrete examples of weak cluster structures and substructures. 
In Chapter \ref{chap3} these examples will be revisited, and used to model
cluster algebras and subcluster algebras.

We denote by $P_i$ the indecomposable projective module associated to vertex $i$, 
by $J$ the radical of a ring, and by
$S_i$ the simple top of $P_i$. 
Usually, we represent a module $M$ by its radical filtration, 
the numbers in the first row represent the indices of the simples in $M/J M$, and
the numbers in the $i$'th row represent the indices of the simples in $J^{i-1} M / J^i M$. 
e.g. 
$$\begin{matrix} & 2 & \\ 1 & & 3 \\ & 2 & 
\end{matrix}$$
represents the indecomposable projective module $P_2$ for the preprojective algebra of type
$A_3$, which has quiver
$$
\xymatrix{
  1 \ar@<0.5ex>[r] & 2 \ar@<0.5ex>[l] \ar@<0.5ex>[r] & 3 \ar@<0.5ex>[l] \\
}
$$

\noindent
\textbf{Example 1.}
Let $\la$ be the preprojective algebra of a Dynkin quiver $A_4$.
This algebra has quiver
$$
\xymatrix{
  1 \ar@<0.5ex>[r] & 2 \ar@<0.5ex>[l] \ar@<0.5ex>[r] & 3 \ar@<0.5ex>[l] \ar@<0.5ex>[r] & 4 \ar@<0.5ex>[l] \\
}
$$
We consider the modules $P_3$ and $M = J P_3$. 
These modules are represented by their radical filtrations:
$$
\begin{smallmatrix} & & 3 & \\ & 2 & & 4 \\ 1 &  & 3 &  \\ & 2 & &  
\end{smallmatrix} \text{ and }
\begin{smallmatrix} & 2 & & 4 \\ 1 &  & 3 &  \\ & 2 & &  
\end{smallmatrix}
$$
We let $\C'=\Sub P_3$ and $\B = \{X \in \C' | \Ext^1(M,X) = 0 \}$.
The AR-quiver of $\C'$ is given below, where
we name the indecomposables in $\C'$ to ease notation.
The indexing will be explained in Section \ref{c3_sec2}.
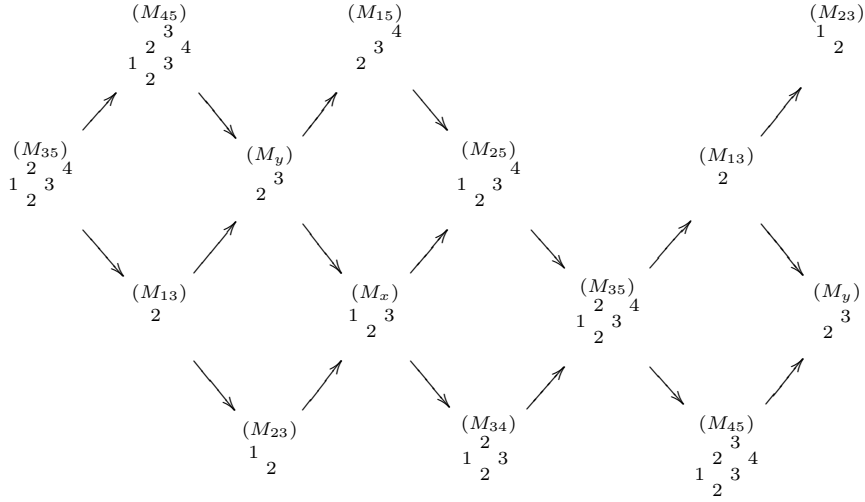
$$
\xymatrix@C0.4cm@R0.5cm{
& 
{\stackrel{(M_{45})}{\begin{smallmatrix} & & 3 & \\ & 2 & & 4  \\ 1 & & 3 & \\ & 2 & & \end{smallmatrix}}} \ar[dr] &
&
{\stackrel{(M_{15})}{\begin{smallmatrix} & & & 4 \\ &  & 3 &   \\  & 2 & & \\ & & & \end{smallmatrix}}}  \ar[dr] &
&
&
&
{\stackrel{(M_{23})}{\begin{smallmatrix} 1 & &  & \\ & 2 & &   \\  & & & \\ & & & \end{smallmatrix}}}  \\ 
{\stackrel{(M_{35})}{\begin{smallmatrix} & 2 &  & 4 \\ 1 &  & 3 &   \\  & 2 & & \\ & & & \end{smallmatrix}}}  \ar[dr] \ar[ur]  &
&
{\stackrel{(M_{y})}{\begin{smallmatrix} & &  & \\ &  & 3 &   \\  & 2 & & \\ & & & \end{smallmatrix}}}  \ar[dr] \ar[ur]  &
&
{\stackrel{(M_{25})}{\begin{smallmatrix} & &  & 4 \\ 1 &  & 3 &   \\  & 2 & & \\ & & & \end{smallmatrix}}}  \ar[dr] &
&
{\stackrel{(M_{13})}{\begin{smallmatrix} & &  & \\ & 2 & &   \\  & & & \\ & & & \end{smallmatrix}}}  \ar[dr] \ar[ur] & \\
&
{\stackrel{(M_{13})}{\begin{smallmatrix} & &  &  \\ &  2 &  &   \\  & & & \\ & & & \end{smallmatrix}}}  \ar[dr] \ar[ur] &
&
{\stackrel{(M_{x})}{\begin{smallmatrix} & &  & \\ 1 &  & 3 &   \\  & 2 & & \\ & & & \end{smallmatrix}}}  \ar[dr] \ar[ur] &
&
{\stackrel{(M_{35})}{\begin{smallmatrix} & 2 &  & 4 \\ 1 &  & 3 &   \\  & 2 & & \\ & & & \end{smallmatrix}}}  \ar[dr] \ar[ur]  &
&
{\stackrel{(M_{y})}{\begin{smallmatrix} & &  & \\ &  & 3 &   \\  & 2 & & \\ & & & \end{smallmatrix}}}    \\
&
&
{\stackrel{(M_{23})}{\begin{smallmatrix} & &  & \\ 1 &  &  &   \\  & 2 & & \\ & & & \end{smallmatrix}}}  \ar[ur] &
&
{\stackrel{(M_{34})}{\begin{smallmatrix} & 2 &  &  \\ 1 &  & 3 &   \\  & 2 & & \\ & & & \end{smallmatrix}}}  \ar[ur] &
&
{\stackrel{(M_{45})}{\begin{smallmatrix} & & 3 & \\ & 2 &  &4   \\ 1 &  & 3 & \\ & 2 & & \end{smallmatrix}}}  \ar[ur] &
}
$$
>From the AR-quiver we see that 
the indecomposable projectives in $\C'$ are $M_{45}, M_{34}, M_{23}$ and $M_{15}$.
The indecomposables of the subcategory $\B$
are obtained from $\C'$ by deleting the indecomposable objects
$P_3 = M_{45}, M_{x}$ and $M_y$.
The category $\B$ is extension closed by definition, and the
indecomposable projectives in $\B$ are $M_{35}, M_{34}, M_{23},
M_{15}$.

Let $T =  M_{34} \oplus M_{23} \oplus M_{13} \oplus M_{15} \oplus M_{35}$, then clearly $\Ext^1(T,T) = 0$, 
and the unique indecomposable
in $\B$ which is not a summand in $T$ is $M_{25}$, which has a non-zero extension with $T$.
Hence $T$ is a cluster tilting object in $\B$, and $\B$ has a cluster structure,
with coefficients $M_{35}, M_{34}, M_{23}, M_{15}$. 

Now, since $M_{45}$ is projective in $\C'$, and $M_x$ as well as $M_y$ has non-zero extensions with $T$,
it is clear that
$T' = T \oplus M_{45}$ is a cluster tilting object in $\C'$, and hence
$\C'$ has a cluster structure, such that we have a substructure for $\B$ induced by $T$.

We claim that the cluster tilting object $T'$ in $\C'$ can be extended to a cluster tilting
object $\widetilde{T} = T' \oplus P_1 \oplus P_2 \oplus P_4 \oplus Z$ of $\mod \la$,
where $Z$ is the $\la$-module with radical filtration 
${\begin{smallmatrix} &  & & \\ 1& &  & \\ &  2 & &  \\  & & 3 & \\ & & & \end{smallmatrix}}$

To see that $\Ext^1(\widetilde{T},\widetilde{T}) = 0$, it is sufficient to
show that $\Ext^1(Z, X \oplus Z) = 0$ for all $X$ in $\C'$. 
There is an exact sequence $0 \to S_4 \to P_1 \to Z \to 0$, and hence 
for every $X$ in $\C'$, there is an exact sequence
$$\Hom(S_4,X \oplus Z) \to \Ext^1(Z,X \oplus Z) \to \Ext^1(P_1,X \oplus Z).$$
Note that $\Hom(S_4, P_3) = 0$
and hence $\Hom(S_4, X) = 0$ for all $X$ in $\C'$, and that $\Hom(S_4,Z) = 0$.
It follows that $\Ext^1(Z,X \oplus Z) = 0$.
Thus $\widetilde{T}$ is a cluster tilting object since it has the correct number
$10=\frac{4\cdot 5}{2}$ of indecomposable direct summands.

\bigskip
\noindent
\textbf{Example 2.}
For our next example, let $\la$ be the preprojective algebra
of type $A_3$. It has the quiver
$$
\xymatrix{
  1 \ar@<0.5ex>[r] & 2 \ar@<0.5ex>[l] \ar@<0.5ex>[r] & 3 \ar@<0.5ex>[l] \\
}
$$
The AR-quiver of $\C=\mod \la$ is given by the following.
We name the indecomposables in $\C$ according to the following table.
The indexing will be explained in Section \ref{c3_sec2}.
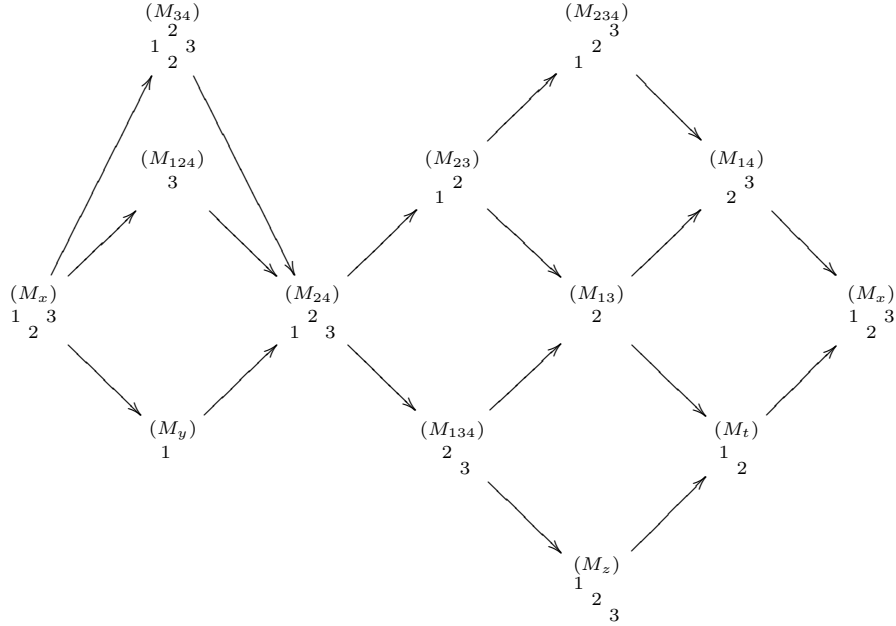
$$
\xymatrix{
& 
{\stackrel{(M_{34})}{\begin{smallmatrix} & 2 & \\ 1 & & 3 \\ & 2 &  \end{smallmatrix}}} \ar[ddr] & & & 
{\stackrel{(M_{234})}{\begin{smallmatrix} & & 3 \\ & 2 & \\ 1 & &  \end{smallmatrix}}} \ar[dr] & & \\
&
{\stackrel{(M_{124})}{\begin{smallmatrix} & & \\ & 3 & \\ & &  \end{smallmatrix}}} \ar[dr] &
&
{\stackrel{(M_{23})}{\begin{smallmatrix} & & \\ & 2 & \\ 1 & &   \end{smallmatrix}}} \ar[dr] \ar[ur] &
&
{\stackrel{(M_{14})}{\begin{smallmatrix} & & \\ & & 3 \\ & 2 &  \end{smallmatrix}}} \ar[dr] & \\
{\stackrel{(M_{x})}{\begin{smallmatrix} & & \\ 1 & & 3 \\ & 2 &  \end{smallmatrix}}} \ar[uur] \ar[ur] \ar[dr] &
&
{\stackrel{(M_{24})}{\begin{smallmatrix} & & \\ & 2 & \\ 1 & & 3 \end{smallmatrix}}} \ar[dr] \ar[ur] &
&
{\stackrel{(M_{13})}{\begin{smallmatrix} & & \\ & 2 & \\ & &  \end{smallmatrix}}} \ar[dr] \ar[ur] & 
&
{\stackrel{(M_{x})}{\begin{smallmatrix} & & \\ 1 & & 3 \\ & 2 &  \end{smallmatrix}}} \\
&
{\stackrel{(M_{y})}{\begin{smallmatrix} & & \\ 1 & & \\ & &  \end{smallmatrix}}} \ar[ur] &
&
{\stackrel{(M_{134})}{\begin{smallmatrix} & & \\ & 2 & \\ & & 3  \end{smallmatrix}}} \ar[dr] \ar[ur]  &
&
{\stackrel{(M_{t})}{\begin{smallmatrix} & & \\ 1 & & \\ & 2 &  \end{smallmatrix}}} \ar[ur] & \\
& & & & 
{\stackrel{(M_{z})}{\begin{smallmatrix} 1 & & \\ & 2 & \\ & & 3  \end{smallmatrix}}} \ar[ur]  & & 
}
$$
The indecomposable projectives in $\C$ are $M_{34}, M_z, M_{234}$.

Let $\B$ be the full subcategory of $\C$ generated by $P_2 \oplus P_3$.
Then $\B=\add (M_{34} \oplus M_{124} \oplus M_{24} \oplus M_{23} \oplus M_{134} \oplus M_{234} 
\oplus M_{13} \oplus M_{14})$.
In addition to $M_{34}, M_{234}$, also $M_{134}$ becomes projective in $\B$. 
It is straightforward to see that $M_{23} \oplus M_{13}$ is extension-free,
so $T = M_{34} \oplus M_{234} \oplus M_{134} \oplus M_{23} \oplus M_{13}$ has $\Ext^1(T, T)= 0$.
Let $\widetilde{T} =  T \oplus M_{z}$. Then also
$\Ext^1(\widetilde{T}, \widetilde{T})= 0$.
Since $\widetilde{T}$ has the correct number of indecomposable direct 
summands $6=\frac{3\cdot 4}{2}$, it is a cluster tilting object.
Hence $T$ is a cluster tilting object in $\B$.

\bigskip
\noindent
\textbf{Example 3.}
In this example we let $\la$ be the preprojective algebra of a Dynkin quiver $D_4$.
$$
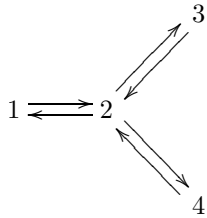
\xymatrix{
& & 3 \ar@<0.5ex>[dl] \\
1 \ar@<0.5ex>[r] & 2 \ar@<0.5ex>[l] \ar@<0.5ex>[ur] \ar@<0.5ex>[dr] & \\
& & 4 \ar@<0.5ex>[ul]
}
$$
We consider the subcategory $\B=\Sub P_2$.
Using Corollary \ref{cor2.6} we have that $\B$ is extension closed.
We know by Theorem \ref{teoI3.1} that $\B$ has a cluster tilting object that
can be extended to a cluster tilting object for $\C=\mod \la$.

The following gives $P_2$ as a representation of the quiver with relations
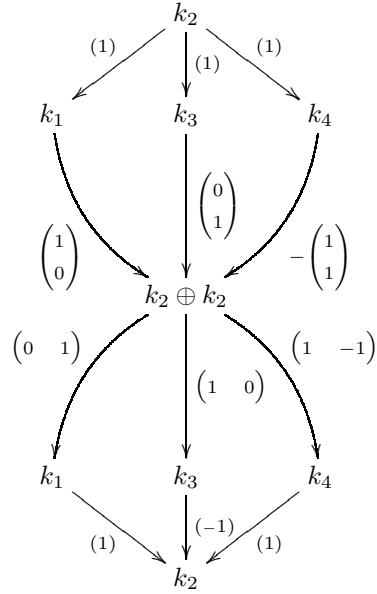
\[ \xymatrix{ & k_2 \ar[ld]_{(1)} \ar[d]^{(1)} \ar[rd]^{(1)} & \\ 
k_1\ar@/_1.1pc/[rdd]_{\begin{pmatrix} \p 1 \\ \p 0 \end{pmatrix}} 
& k_3 \ar[dd]^{\begin{pmatrix} \p 0 \\ \p 1 \end{pmatrix}} 
& k_4 \ar@/^1.1pc/[ldd]^{{\p -}\begin{pmatrix}  \p 1 \\ \p 1 \end{pmatrix}}  \\ \\
&k_2 \oplus k_2 
\ar@/_1.1pc/[ldd]_{\begin{pmatrix} \p 0 & \p 1 \end{pmatrix}} 
\ar[dd]^{\begin{pmatrix} \p 1 & \p 0 \end{pmatrix}}
\ar@/^1.1pc/[rdd]^{\begin{pmatrix} \p 1 & \p -1 \end{pmatrix}} & \\ \\
k_1\ar[rd]_{(1)} & k_3 \ar[d]^{({\s -}1)} 
& k_4 \ar[ld]^{(1)} \\
&k_2 & } \]
The modules in $\B$ do not necessarily have a simple socle, and in fact
the subcategory is not of finite type. As noted earlier, it is functorially finite.
However, the indecomposable direct summands in the 
cluster tilting object we will construct all have simple socle.

The indecomposable submodules of $P_2$ we will need to construct a cluster tilting object 
have the following radical filtrations. The indexing will be explained in Chapter \ref{chap3}.
\begin{center}
\begin{tabular}{|rr|rr|rr|rr|}
\hline
$M_{16}$ & ${\begin{smallmatrix} & & & & \\ & & & & \\  & 3 & & 4&   &  \\ & & 2 & & \\ & & & & \\ & & & &
\end{smallmatrix}}$ &
$M_{24}$ &  ${\begin{smallmatrix}  & & & & \\ & & &  & \\ & 1 &  & 3 &   \\ & & 2&  &  \\ & & & &  \\ & & & & 
\end{smallmatrix}}$ &
$M_{25}$ &  ${\begin{smallmatrix}  & & & & \\ & & &  & \\ & 1 &  & 4  & \\ & & 2 &  &  \\ & & & & \\ & &  & & 
\end{smallmatrix}}$ 
&
$M_{26}$ &  ${\begin{smallmatrix}  & & & & \\ & &  & &   \\ & 1 & 3 & 4  & \\ & & 2 & & \\  & & & & \\ & & & & 
\end{smallmatrix}}$ \\
\hline
$M_{68}$ & ${\begin{smallmatrix}  & & & & \\ 1 & & 3 & & 4   \\  & 2 & & 2 & \\ 1 & & 3 & & 4  \\ & & 2 & &  \\ & & & & 
\end{smallmatrix}}$ 
&
$M_{18}$ & ${\begin{smallmatrix}  & & & & \\ & & 1 & &   \\  & & 2 &  & \\ & 3 && 4 & \\ && 2 & & \\  & & & & 
\end{smallmatrix}}$ 
&
$M_{-}$ & ${\begin{smallmatrix} & & & & \\ & & 4 &  &  \\ & & 2 & &   \\  & 1 & & 3 & \\ & & 2 & & \\ & & & &
\end{smallmatrix}}$ 
&
$M_{+}$ & ${\begin{smallmatrix} & & & & \\ & & 3 &  & \\ & & 2 & &  \\ & 1 & & 4 &  \\ & & 2 & & \\ & & & &
\end{smallmatrix}}$ \\
\hline
\end{tabular}
\end{center}
The indecomposable projectives in $\B$ are $P_2, M_{18}, M_{+}$ and $M_{-}$. 
This follows from the following.

\begin{lemma}
Let $\la$ be a finite dimensional algebra with $X$ in $\mod\la$ such that $\Sub X$ is extension 
closed in $\mod\la$. Then the indecomposable projective objects in $\Sub X$ are of the form 
$P/\alpha(P)$ where $P$ is an indecomposable projective $\la$-module and $\alpha(P)$ is the smallest 
submodule of $P$ such that $P/\alpha(P)$ is in $\Sub X$.
\end{lemma}

\begin{proof}
  For convenience of the reader, we include a proof.
Let $P$ be indecomposable projective in $\mod\la$. It is clear that there is a smallest 
submodule $\alpha(P)$ of $P$ such that $P/\alpha(P)$ is in $\Sub X$. For, if $A$ and $B$ are 
submodules of $P$ with $P/A$ and $P/B$ in $\Sub X$, then clearly $P/{A\cap B}\subseteq P/A\oplus P/B$ 
is in  $\Sub X$. It is clear that the natural map $f\colon P\to P/\alpha(P)$ is a minimal left 
$\Sub X$-approximation and that every module in $\Sub X$ is a factor of a direct sum of $\la$-modules 
of the form  $P/\alpha(P)$ for $P$ indecomposable projective. To see that each $P/\alpha(P)$ is 
projective in $\Sub X$, consider the exact sequence $0\to A\xto{s} B\xto{t}P/\alpha(P)\to 0$ in $\Sub X$. 
Then there is some $ u\colon P\to B$ such that $ut=f$. Since $B$ is in $\Sub X$, then $u(\alpha(P))=0$, 
so the sequence splits. Clearly there are no other indecomposable
projectives in $\Sub X$ since all modules in 
$\Sub X$ are factors of direct sums of those of the form $P/\alpha(P)$.
\end{proof}

In addition we need the following, where we leave the details to the reader.

\begin{lemma}
Let $M = M_{16} \oplus M_{24} \oplus M_{25} \oplus M_{26} \oplus M_{68}$.
Then we have $\Ext^1(M,M)= 0$.
\end{lemma}

Hence, the module $T = M \oplus P_2 \oplus M_{18} \oplus M_{+} \oplus M_{-}$ in $\B$ is rigid. 
If we add the other projectives we obtain the module $\widetilde{T} = T \oplus P_1 \oplus P_3 \oplus P_4$,
which also satisfies $\Ext^1(\widetilde{T},\widetilde{T}) = 0$.

Since $\widetilde{T}$ has the correct number $12=4\cdot 3$ of indecomposable
summands, it is a cluster tilting object in $\C=\mod \la$.
It is also clear from this that $T$ is a cluster tilting object in
$\B=\Sub P_2$, since
we added only projectives/injectives to $T$ to obtain $\widetilde{T}$. 
Note that $\B$ has a substructure of the cluster structure of $\C$.


\section{Preprojective algebras for non-Dynkin quivers}\label{chap2}

In this chapter we deal with completions of preprojective algebras of a finite  
connected quiver $Q$ with no oriented cycles, and mainly those which
are not Dynkin. In this case the modules of finite length coincide
with the nilpotent modules over the preprojective algebra.
These algebras $\la$ are known to be derived 2-CY (see
\cite{b,cb,bbk,gls2}).
Tilting $\la$-modules of projective dimension at most one were
investigated in \cite{ir} when the quiver $Q$ is a (generalized)
extended Dynkin quiver.
It was shown that such tilting modules are exactly
the ideals in $\la$ which are finite products of two-sided
ideals $I_i=\la(1-e_i)\la$, where $e_1,\cdots ,e_n$ correspond to
the vertices of the quiver, and that they are in
one-one correspondence with the elements of the corresponding Weyl
group, where $w=s_{i_1}\cdots s_{i_k}$ corresponds to
$I_w=I_{i_1}\cdots I_{i_k}$. Here we generalize some of the results
from \cite{ir} beyond the 
noetherian case. In particular, we show that any finite product of
ideals of the form $I_i$ is a tilting module, and show that there is a
bijection between cofinite tilting ideals and elements of the
associated Coxeter group $W$.

For any descending chain of tilting ideals of the form $\la \supseteq
I_{i_1}\supseteq I_{i_1}I_{i_2}\supseteq I_{i_1}I_{i_2}\cdots
I_{i_k}\supseteq\cdots$ we show that for
$\la_m=\la/{I_{i_1}\cdots I_{i_m}}$, the categories $\Sub\la_m$ and
$\underline{\Sub}\la_m$ are respectively stably 2-CY and 2-CY with nice cluster tilting objects. 
In this way we get, for any $w\in W$, a stably 2-CY category
$\C_w=\Sub(\Lambda/I_w)$, and for any reduced expression
$w=s_{i_1}\cdots s_{i_k}$, a cluster tilting object
$\bigoplus_{j=1}^k\Lambda/I_{s_{i_1}\cdots s_{i_j}}$ in $\C_w$.
We also construct cluster tilting subcategories of the derived 2-CY category $\fl\la$. 
This way we get many examples of weak cluster
structures without loops or 2-cycles which are then cluster structures by Theorem \ref{teoI1.6}. 
We also get many examples of substructures. In particular, any
cluster category and the stable category $\ul{\mod}\la$ of a preprojective algebra of Dynkin  
type occur amongst this class.
We give a description of the quivers of the cluster tilting
objects/subcategories in terms of the associated reduced expressions.
For example, the quiver of the preprojective component of the hereditary
algebra with additional arrows from $X$ to $\tau X$ occur this way.
In Section \ref{c3_sec3} results in this chapter are used
to show that coordinate rings of some unipotent cells of $\SL_2(
\mathbb{C}[t,t^{-1}] )$ have a cluster algebra structure.

We refer to \cite{i3} for corresponding results for $d$-CY algebras.\\
\hspace{7mm}

\subsection{Tilting modules over 2-CY algebras}\label{c2_sec1}
${}$ \\
Let $Q$ be a finite connected quiver without oriented cycles which is
not Dynkin, $k$ an algebraically closed field and $\la$ the completion
of the associated preprojective algebra. In \cite{ir} the
tilting $\la$-modules of projective dimension at most one were investigated
in the noetherian case, that is, when $Q$ is extended Dynkin \cite{bgl} (and
also the generalized ones having loops). In this section we
generalize some of these results to the non-noetherian case,
concentrating on the aspects that will be needed for our
construction of new 2-CY categories with cluster tilting
objects/subcategories in the next sections.
Note that since $\Lambda$ is complete, the Krull-Schmidt theorem holds
for finitely generated projective $\Lambda$-modules.

We say that a finitely presented $\la$-module $T$ is a {\em tilting module} 
if 
(i) there exists an exact sequence $0\to P_n\to\cdots\to
P_0\to\Lambda\to 0$ with finitely generated projective
$\Lambda$-modules $P_i$, (ii) $\Ext_{\la}^i(T,T)=0$ for any $i>0$, (iii) 
there exists an exact sequence $0\to\Lambda\to T_0\to\cdots\to T_n\to 0$ with
$T_i$ in $\add T$.

We say that $T\in{\bf D}(\Mod\la)$ a {\em tilting complex} \cite{rick} if 
(i$'$) $T$ is quasi-isomorphic to an object in the
category ${\bf K}^{\bo}(\pr\la)$ of bounded complexes of finitely generated
projective $\la$-modules $\pr\la$,
(ii$'$) $\Hom_{{\bf D}(\mod\la)}(T,T[i])=0$ for any $i\neq0$,
(iii$'$) $T$ generates ${\bf K}^{\bo}(\pr\la)$.

A tilting module is nothing but a module which is a tilting complex
since the condition (iii) can be replaced by (iii$'$).
A {\em partial tilting complex} is a direct summand of a tilting
complex. A {\em partial tilting module} is a module which is a partial
tilting complex.

Let $1,\cdots,n$ denote the vertices in $Q$, and let $e_1,\cdots,e_n$
be the corresponding idempotents. For each $i$ we denote by $I_i$ the
ideal $\la(1-e_i)\la$. Then $S_i=\la/I_i$ is a simple $\la$-module and 
$\la^{\op}$-module  since by assumption there are no loops in the quiver. 
We shall show that each $I_i$, and any finite
product of such ideals, is a tilting ideal in $\la$, and give
some information about how the different products are related. But first
we give several preliminary results, where new proofs are needed
compared to \cite{ir} since we do not assume $\la$ to be noetherian.

\begin{lemma}\label{lemII1.1}
Let $T$ be a partial tilting $\Lambda$-module of projective dimension at
most 1 and $S$ a simple $\Lambda^{\op}$-module. Then at least one
of the statements $S\otimes_\Lambda T=0$ and ${\Tor}^\Lambda_1(S,T)=0$ holds.
\end{lemma}

\begin{proof}
We only have to show that there is a projective resolution $0\to
P_1\to P_0\to T\to0$ such that $P_0$ and $P_1$ do not have a common
summand. This is shown as in \cite[1.2]{hu}.
\end{proof}

Recall that 
for rings $\la$ and $\Gamma$, we call an object $T$ in ${\bf D}(\Mod\Lambda\otimes_{\mathbb{Z}}\Gamma^{\op})$ a {\em two-sided tilting complex} if $T$ is a 
tilting complex in ${\bf D}(\Mod \Lambda)$ and $\End_{{\bf D}(\Mod \Lambda)}(T)\simeq\Gamma$ naturally.

The following result is useful (see \cite{rick}\cite[1.7]{ye}).

\begin{lemma}\label{lemII1.3}
Let $T\in{\bf D}(\Mod \Lambda\otimes_{\mathbb{Z}}\Gamma^{\op})$
be a two-sided tilting complex.
\begin{itemize}
\item[(a)] For any tilting complex (respectively, partial tilting
  complex) $U$ of $\Gamma$, we have a tilting complex (respectively,
  partial tilting complex) $T\stackrel{\bf L}{\otimes}_{\Gamma}U$ of $\Lambda$
such that $\End_{{\bf D}(\Mod\Lambda)}(T\stackrel{\bf L}{\otimes}_{\Gamma}U)\simeq\End_{{\bf D}(\Mod\Gamma)}(U)$.
\item[(b)] ${\bf R}\Hom_{\Lambda}(T,\Lambda)$ and ${\bf
    R}\Hom_{\Gamma^{\op}}(T,\Gamma)$ are two-sided tilting complexes and
  isomorphic in ${\bf D}(\Mod \Gamma\otimes_{\mathbb{Z}}\Lambda^{\op})$.
\end{itemize}
\end{lemma}

We collect some basic information on preprojective algebras. 

\begin{proposition}\label{2CY}
Let $\Lambda$ be the completion of the preprojective algebra of a finite
connected non-Dynkin diagram without loops. 
\begin{itemize}
\item[(a)] Let $\Gamma$ be the completion of $\la\otimes_{kQ_0}\la^{\op}$ with respect to the ideal
$J\otimes_{kQ_0}\la^{\op}+\la\otimes_{kQ_0}J^{\op}$ where $J$ is the radical of $\la$.
Then there exists a commutative diagram
$$\begin{array}{ccccccccccc}
0&\to&P_2&\stackrel{f_2}{\to}&P_1&\stackrel{f_1}{\to}&P_0&\stackrel{}{\to}&\Lambda&\to&0\\
&&\downarrow\wr&&\downarrow\wr&&\downarrow\wr&&\downarrow\wr&&\\
0&\to&\Hom_{\Gamma}(P_0,\Gamma)&\stackrel{f_1\cdot}{\to}&\Hom_{\Gamma}(P_1,\Gamma)&\stackrel{f_2\cdot}{\to}&
\Hom_{\Gamma}(P_2,\Gamma)&\stackrel{}{\to}&\Lambda&\to&0 
\end{array}$$
of exact sequences of $\Gamma$-modules such that each $P_i$ is
a finitely generated projective $\Gamma$-module
and $P_0\simeq P_2\simeq\Gamma$.
\item[(b)]There exists a functorial
  isomorphism $\Hom_{{\bf D}(\Mod\la)}(X,Y[1])\simeq 
  D\Hom_{{\bf D}(\Mod\la)}(Y,X[1])$ for any $X\in{\bf D}^{\bo}(\fl\la)$ and $Y\in{\bf K}^{\bo}(\pr\la)$.
\item[(c)]$\fl\la$ is derived 2-CY and $\gl\la=2$. In particular, any left ideal $I$ of $\la$ satisfies $\pd{}_\la I\le1$.
\item[(d)]$\Ext^i_\la(X,\la)=0$ for $i\neq 2$ and
  $\Ext^2_\la(X,\la)\simeq DX$ for any $X\in\fl\la$.
\end{itemize}
\end{proposition}

\begin{proof}
(a) See \cite[Section 8]{gls1} and \cite[Section 4.1]{bbk}.

(b) This follow from (a) and \cite[4.2]{b}.

(c)(d) Immediate from (a) and (b).\end{proof}

We are now ready to show that each $I_i$, and a finite product of such
ideals, is a tilting module.

\begin{proposition}\label{propII1.4}
$I_i$ is a tilting $\Lambda$-module of projective dimension at most
one and $\End_\Lambda(I_i)=\Lambda$.
\end{proposition}

\begin{proof}
We have $\Ext^n_\Lambda(S_i,\Lambda)\simeq
D\Ext^{2-n}_\Lambda(\Lambda,S_i)=0$ for $n=0,1$ by Proposition \ref{2CY}. 
Applying $\Hom_\Lambda(\text{ },\Lambda)$ to the exact sequence $0\to I_i\to\Lambda\to
S_i \to 0$, we get $\Hom_\Lambda(I_i,\Lambda)=\Lambda$. Applying
$\Hom_\Lambda(I_i,\text{ })$, we get an exact sequence $0\to 
\End_\Lambda(I_i) \to \Hom_\Lambda(I_i,\Lambda)\to 
\Hom_\Lambda(I_i,S_i)$. Since $\Hom_\Lambda(I_i,S_i)=0$, we
have $\End_\Lambda(I_i)= \Hom_\Lambda(I_i,\Lambda)=\Lambda$.

Applying $\ \otimes_\Lambda S_i$ to the exact sequence in Proposition
\ref{2CY}(a) we have a projective resolution
\begin{equation}\label{simple resolution}
0 \to \Lambda e_i\stackrel{g}{\to}P\stackrel{f}{\to}\Lambda e_i\to
S_i \to 0
\end{equation}
with $\Im f=I_ie_i$ and $P \in \add \Lambda(1-e_i)$.
In particular $I_i=\Im f\oplus\Lambda(1-e_i)$ is a finitely presented $\la$-module with $\pd I_i\le1$.

We have $\Ext^1_\la(I_i,I_i)\simeq\Ext^2_\la(S_i,I_i)\simeq D\Hom_\la(I_i,S_i)=0$. Using \eqref{simple resolution}, we have an exact sequence
\[0\to\la\to P\oplus\la(1-e_i)\to I_ie_i\to0\]
such that the middle and the right terms belong to $\add I_i$.
Thus $I_i$ is a tilting $\la$-module.
\end{proof}


\begin{proposition}\label{propII1.5}
Let $T$ be a tilting $\Lambda$-module of projective dimension at
most one.
\begin{itemize}
\item[(a)]{If $\Tor^\Lambda_1(S_i,T)=0$, then $I_i\stackrel{\bf
    L}{\otimes}_\Lambda T = I_i \otimes_{\la} T =I_iT$ is a tilting $\Lambda$-module of
    projective dimension at most one.}
\item[(b)]{$I_iT$ is always a tilting $\Lambda$-module of projective
    dimension at most one, and $\End_\la(I_iT)\simeq\End_\Lambda(T)$.}
\end{itemize}
\end{proposition}

\begin{proof}
(a) Since $\Tor^\Lambda_1(I_i,T)=\Tor^\Lambda_2(S_i,T)=0$
because $\pd T\le1$, we have $I_i\stackrel{\bf L}{\otimes}_\Lambda
T=I_i\otimes_\Lambda T$.
Since we have an exact sequence
$$
0= \Tor^\Lambda_1(S_i,T)\to I_i\otimes_\Lambda
T\to\Lambda\otimes_\Lambda T\to S_i\otimes_\Lambda T\to0,
$$
we have $I_i\otimes_\Lambda T=I_iT$.
Thus $I_i\stackrel{\bf L}{\otimes}_\Lambda T=I_iT$ is a tilting $\la$-module by Lemma \ref{lemII1.3} and Proposition \ref{propII1.4}.
Since $\pd T \leq 1$ and $\pd T/I_i T \leq 2$ by Proposition \ref{2CY}, we have $\pd I_iT\le1$. 

\noindent 
(b) By Lemma \ref{lemII1.1}, either $S_i \otimes_\Lambda T=0$ or $\Tor^\Lambda_1(S_i,T)=0$
holds. If $S_i\otimes_\Lambda T=0$, then $I_iT=T$ holds. If $\Tor^\Lambda_1(S_i,T)=0$, 
then we apply (a). For the rest we use Lemma \ref{lemII1.3}.
\end{proof}

A left ideal $I$ of $\la$ is called {\em
  cofinite} if $\Lambda/I\in {\fl} \Lambda$, and called {\em tilting}
(respectively, {\em partial tilting}) if it is a tilting
(respectively, partial tilting) $\la$-module.
Similarly, a {\em cofinite} (respectively, {\em (partial) tilting})
{\em right ideal} of $\la$ is defined. An ideal $I$ of $\la$ is called
{\em cofinite tilting} if it is cofinite tilting as a left and right ideal.
We denote by $\langle I_1,...,I_n\rangle$ the ideal semigroup
generated by $I_1,...,I_n$. Then we have the following result. 

\begin{theorem}\label{teoII1.6}
\begin{itemize}
\item[(a)]Any $T\in\langle I_1,...,I_n\rangle$ is a 
cofinite tilting ideal and satisfies ${\rm End}_\Lambda(T)=\Lambda$.
\item[(b)]Any cofinite tilting ideal of $\Lambda$ belongs to $\langle
  I_1,...,I_n\rangle$.
\item[(c)] Any cofinite partial tilting left (respectively, right)
  ideal of $\Lambda$ is a cofinite tilting ideal.
\end{itemize}
\end{theorem}

\begin{proof}
(a) This is a direct consequence of Propositions \ref{propII1.4} and
\ref{propII1.5}.

(b)(c) Let $T$ be a cofinite partial tilting left ideal of $\Lambda$.
If $T\neq\Lambda$, then there exists a simple submodule $S_i$ of $\Lambda/T$.
Since $\Hom_\la(S_i,\la)=0$, we have $\Ext^1_\Lambda(S_i,T)\neq0$. Thus we have
$\Tor_1^\Lambda(S_i,T)\simeq D\Ext^1_\Lambda(T,S_i)\simeq\Ext^1_\Lambda(S_i,T)\neq0$.
By Lemma \ref{lemII1.1}, we have $S_i\otimes_\Lambda T=0$.

Put $U={\bf R}\Hom_\Lambda(I_i,T)$. By Lemma \ref{lemII1.3}, we have that $U\simeq
{\bf R}\Hom_\Lambda(I_i,\Lambda)\stackrel{{\bf L}}{\otimes}_\Lambda T$
is a partial tilting complex of $\Lambda$. Since $\pd I_i\le 1$ and
$\Ext^1_\Lambda(I_i,T)\simeq\Ext^2_\Lambda(S_i,T)\simeq
D\Hom_\Lambda(T,S_i)\simeq S_i\otimes_\Lambda
T=0$, we have $U=\Hom_\Lambda(I_i,T)$, which
is a partial tilting $\la$-module.
Since we have a commutative diagram
$$\begin{array}{ccccccccc}
0=\Hom_\la(S_i,\la)&\to&\la&\to&\Hom_\la(I_i,\la)&\to&\Ext^1_\la(S_i,\la)=0&&\\
&&\cup&&\cup&&&&\\
&&T&\to&\Hom_\la(I_i,T)&\to&\Ext^1_\la(S_i,T)&\to&0\end{array}$$
of exact sequences, $U$ is a cofinite partial
tilting left ideal of $\Lambda$
containing $T$ properly, such that $U/T$ is a direct sum of copies of $S_i$.
By $S_i\otimes_\Lambda T=0$, we have $T=I_iU$.
Thus we have $T\in\langle I_1,\cdots,I_n\rangle$ by induction on the length of $\Lambda/T$.
\end{proof}

We here pose the following question, where there is a positive answer
in the extended Dynkin case \cite{ir}.

\begin{question}
For any tilting $\la$-module $T$ of projective dimension at most one,
does there exist some $U$ in $\langle I_1, \cdots, I_n \rangle$ such that $\add T = \add U$?
\end{question}

We have some stronger statements on products of the ideals $I_i$, generalizing results in the noetherian case 
from \cite{ir}.

\begin{proposition}\label{propII1.7}
  The following equalities hold for multiplication of
  ideals.
  \begin{itemize}
  \item[(a)] $I_i^2=I_i$,
  \item[(b)] $I_iI_j=I_jI_i$ if there is no arrow between $i$ and $j$ in $Q$,
  \item[(c)] $I_iI_jI_i=I_jI_iI_j$ if there is precisely one arrow between $i$ and $j$ in $Q$.
  \end{itemize}
\end{proposition}
\begin{proof}
(a) is obvious.

Parts (b) and (c) are proved in \cite[6.12]{ir} for module-finite
2-CY algebras. Here we give a direct proof for an arbitrary preprojective
algebra $\la$ associated with a finite quiver $Q$ without oriented cycles.
Let $I_{i,j} = \Lambda(1-e_i-e_j)\Lambda$.
Then any product of ideals $I_i$ and $I_j$ contains $I_{i,j}$. If there is no arrow from
$i$ to $j$, then
$\Lambda/I_{i,j}$ is semisimple. Thus $I_iI_j$ and $I_jI_i$
are contained in $I_{i,j}$, and we have $I_iI_j=I_{i,j}=I_jI_i$.

If there is precisely one arrow from $i$ to $j$, then $\la/I_{i,j}$ 
is the preprojective algebra of type $A_2$. Hence
there are two indecomposable projective $\Lambda/I_{i,j}$-modules,
whose Loewy series are ${i\choose j}$ and ${j\choose i}$.
Thus $I_iI_jI_i$ and $I_jI_iI_j$ are contained in $I_{i,j}$,
and we have $I_iI_jI_i=I_{i,j}=I_jI_iI_j$.
\end{proof}

Now let $W$ be the Coxeter group associated to the quiver $Q$,
so $W$ has generators $s_1,...,s_n$ with relations
$s_i^2=1$, $s_is_j=s_js_i$ if there is no arrow between $i$ and $j$ in $Q$,
and $s_is_js_i=s_js_is_j$ if there is a precisely one arrow between
$i$ and $j$ in $Q$.

\begin{theorem}\label{teoII1.8}
  There exists a bijection $W\to\langle I_1,...,I_n\rangle$.
  It is given by $w\mapsto I_w =I_{i_1}I_{i_2}...I_{i_k}$ for any reduced
  expression $w=s_{i_1}s_{i_2}...s_{i_k}$.
\end{theorem}
\begin{proof}
  The corresponding result was proved in \cite{ir} in the noetherian case, 
using a partial order of tilting modules. Here we use instead properties of Coxeter groups.

We first show that the map is well-defined.
Take two reduced expressions
$w=s_{i_1}s_{i_2}...s_{i_k}=s_{j_1}s_{j_2}...s_{j_k}$.
By \cite[3.3.1(ii)]{bb}, two words $s_{i_1}s_{i_2}...s_{i_k}$
and $s_{j_1}s_{j_2}...s_{j_k}$ can be connected by a sequence of
the following operations:
(i) replace $s_is_j$ by $s_js_i$ (there is no arrow from $i$ to $j$),
(ii) replace $s_is_js_i$ by $s_js_is_j$ (there is precisely one arrow from $i$ to $j$).
Consequently, by Proposition \ref{propII1.7}(b)(c), we have
$I_{i_1}I_{i_2}...I_{i_k}=I_{j_1}I_{j_2}...I_{j_k}$.
Thus the map is well-defined.

Next we show that the map is surjective.
For any $I\in\langle I_1,...,I_n\rangle$, take an expression
$I=I_{i_1}I_{i_2}...I_{i_k}$ with a minimal number $k$.
Let $w =s_{i_1}s_{i_2}...s_{i_k}$.
By \cite[3.3.1(i)]{bb}, a reduced expression of $w$ is obtained from
the word $s_{i_1}s_{i_2}...s_{i_k}$ by a sequence of
the operations (i)(ii) above and
(iii) remove $s_is_i$.
By Proposition \ref{propII1.7}, the operation (iii) can not appear since $k$ is minimal.
Thus $w=s_{i_1}s_{i_2}...s_{i_k}$ is a reduced expression, and we have $I=I_w$.

Finally we show that the map is injective in a similar way as in
\cite{ir}. Let $\E ={\bf K}^{\bo}(\text{pr } \Lambda)$. For any $i$,
we have an autoequivalence $I_i\stackrel{{\bf L}}{\otimes}_\Lambda\ $
of $\E$ and an automorphism $[I_i\stackrel{{\bf
L}}{\otimes}_\Lambda\ ]$ of the Grothendieck group $K_0(\E)$.
By \cite[proof of 6.6]{ir}, we have the action
$s_i\mapsto[I_i\stackrel{{\bf L}}{\otimes}_\Lambda\ ]$ of $W$ on
$K_0(\E)\otimes_{\bf Z}\mathbb{C}$, which is shown to be faithful \cite[4.2.7]{bb}.

For any reduced expression $w=s_{i_1}s_{i_2}...s_{i_k}$, we have
$I_w=I_{i_1}\stackrel{{\bf L}}{\otimes}_\Lambda...\stackrel{{\bf L}}{\otimes}_\Lambda I_{i_k}$,
by Proposition \ref{propII1.5}(a) and the minimality of $k$.
Thus the action of $w$ on $K_0(\E)\otimes_{\mathbb{Z}}\mathbb{C}$ coincides with
$[I_w\stackrel{{\bf L}}{\otimes}_\Lambda\ ]$.
In particular, if $w,w'\in W$ satisfy $I_w=I_{w'}$,
then the actions of $w$ and $w'$ on $K_0(\E)\otimes_{\mathbb{Z}}\mathbb{C}$ coincide,
so we have $w=w'$ by the faithfulness of the action.
\end{proof}

We denote by $l(w)$ the length of $w\in W$.
We say that an infinite expression $s_{i_1}s_{i_2}\cdots s_{i_k}\cdots$ is
{\it reduced} if the expression $s_{i_1}s_{i_2}\cdots s_{i_k}$ is
reduced for any $k$.

\begin{proposition}\label{length}
Let $w\in W$ and $i\in\{1,\cdots,n\}$.
If $l(ws_i)>l(w)$, then we have $I_wI_i=I_{ws_i}\subsetneq I_w$.
If $l(ws_i)<l(w)$, then we have $I_wI_i=I_w\subsetneq I_{ws_i}$.
\end{proposition}

\begin{proof}
Let $w=s_{i_1}\cdots s_{i_k}$ be a reduced expression.
If $l(ws_i)>l(w)$, then $ws_i=s_{i_1}\cdots s_{i_k}s_i$ is a reduced
expression, so the assertion follows from Theorem \ref{teoII1.8}.
If $l(ws_i)<l(w)$, then $u=ws_i$ satisfies $l(us_i)>l(u)$, so
$I_uI_i=I_{us_i}=I_w\subsetneq I_u$.
\end{proof}

Let $s_{i_1}s_{i_2}\cdots s_{i_k}\cdots$ be a (finite or infinite) expression
such that $i_k\in\{1,\cdots,n\}$. Let
$$w_k=s_{i_1}s_{i_2}\cdots s_{i_k},\ \ T_k=I_{w_k}=I_{i_1}I_{i_2}\cdots I_{i_k}\ \mbox{ and }\ \la_k=\la/T_k.$$
We have a descending chain
$$\Lambda=T_0\supseteq T_1\supseteq T_2\supseteq \cdots$$
of cofinite tilting ideals of $\Lambda$, and a chain
$$\la_1\leftarrow\la_2\leftarrow\la_3\leftarrow\cdots$$
of surjective ring homomorphisms.
We have the following properties of the chain.

\begin{proposition}\label{strict}
\begin{itemize}
\item[(a)] If $T_{m-1}\neq T_m$, then $\la_m$ differs from $\la_{m-1}$ in exactly
one indecomposable summand $\Lambda_me_{i_m}$.
\item[(b)] Let $k\le m$. Then $\la_ke_{i_k}$ is a projective $\la_m$-module
  if and only if $i_k\notin\{i_{k+1},i_{k+2},\cdots,i_m\}$.
\item[(c)] $T_1\supsetneq T_2\supsetneq T_3\supsetneq\cdots$ holds if and only if $s_{i_1}s_{i_2}\cdots$ is
reduced.
\end{itemize}
\end{proposition}

\begin{proof}
(a) This follows from $T_m(1-e_{i_m}) = T_{m-1} I_{i_m} (1- e_{i_m}) = T_{m-1}(1-e_{i_m})$.

(b) If $i_k\notin\{i_{k+1},\cdots,i_m\}$, then $\la_ke_{i_k}$ is a
summand of $\la_m$ by (a), so it is a projective $\la_m$-module.
Otherwise, take the smallest $k'$ with $k<k'\le m$ satisfying
$i_k=i_{k'}$. Then we have $\la_ke_{i_k}=\la_{k'-1}e_{i_k}$ and that
$\la_ke_{i_k}$ is a proper factor module of $\la_{k'}e_{i_k}$ by (a).
Hence $\la_ke_{i_k}$ is not a projective $\la_m$-module.

(c) This follows from Proposition \ref{length}.
\end{proof}

Our next goal is to show that $\Ext_{\la}^1(T_k,T_m)=0$  
for $k\le m$. For this the following result will be useful. 

\begin{lemma}\label{lemII1.9}
Let the notation and assumptions be as above. Then $^{\perp_{>0}}{T_{m-1}}\subseteq
{}^{\perp_{>0}}T_m$, where $^{\perp_{>0}}T=\{X\in\mod\la \mid \Ext^i_{\la}(X,T)=0 \text{ for all }i>0\}$.
\end{lemma}
\begin{proof}
We can assume $T_{m-1}\neq T_m$. Then we have that $T_{m-1}\otimes_\Lambda
S_{i_m}\neq0$. Hence $\Tor_1^\Lambda(T_{m-1},S_{i_m})=0$ by Lemma
\ref{lemII1.1}, so
$T_{m-1}\otimes_\Lambda I_{i_m}=T_{m-1}I_{i_m}=T_m$ by Proposition \ref{propII1.5}.
Let $0\to P_1\to P_0\to I_{i_m}\to 0$ be a projective resolution.
We have $\Tor^\Lambda_1(T_{m-1},I_{i_m})\simeq\Tor^\Lambda_2(T_{m-1},S_{i_m})=0$.
Applying $T_{m-1}\otimes_\Lambda\ $, we have an exact sequence
$0\to T_{m-1}\otimes_\Lambda P_1\to T_{m-1}\otimes_\Lambda P_0\to T_m\to 0$.
This immediately implies ${}^{\perp_{>0}} T_{m-1}\subseteq {}^{\perp_{>0}} T_m$.
\end{proof}

We now have the following consequence.

\begin{proposition}\label{propII1.10}
With the above notation and assumptions, we have $\Ext_{\la}^1(T_k,T_m)=0$ for $k\le m$.
\end{proposition}
\begin{proof}
By Lemma \ref{lemII1.9} we have $^{\perp_{>0}}T_k\subseteq
^{\perp_{>0}}T_{k+1}\subseteq\cdots\subseteq^{\perp_{>0}}T_m$. Since $T_k$ is in
$^{\perp_{>0}}T_k$, we then have that $T_k$ is in $^{\perp_{>0}}T_m$, and hence
$\Ext^1_{\la}(T_k,T_m)=0$ for $k\le m$.
\end{proof}



Later we shall use the following observation.

\begin{lemma}\label{hom-calculation}
Assume that the expression $s_{i_1}s_{i_2}\cdots$ is reduced.
Let $T_{k,m}=I_{i_k}\cdots I_{i_m}$ if $k\le m$ and $T_{k,m}=\Lambda$ otherwise.
Then we have $\Hom_\la(T_k,T_m)\simeq T_{k+1,m}=\{x\in\Lambda \mid
T_kx\subseteq T_m\}$ and $\Hom_\la(\la_k,\la_m)\simeq T_{k+1,m}/T_m$.
\end{lemma}

\begin{proof}
Let $U = \{ x \in \la \mid T_kx \subseteq T_m \}\supseteq T_{k+1,m}$.

If $k\ge m$, then clearly $U=\Lambda=T_{k+1,m}$ holds, and
$\Hom_\la(T_k,T_m)\subseteq\End_\la(T_k)\simeq\la$ by Theorem \ref{teoII1.6}.
Thus $\Hom_\la(T_k,T_m)\simeq\Lambda$.

We assume $k<m$.
Since $T_m = T_k \stackrel{\bf L}{\otimes}_\Lambda T_{k+1,m}$ holds by 
Proposition \ref{propII1.5}(a) and Lemma \ref{lemII1.3},
we have ${\bf R} \Hom_{\la}(T_k, T_m)=
{\bf R} \Hom_{\la}(T_k, (T_k \stackrel{\bf L}{\otimes}_\Lambda T_{k+1,m})) = 
{\bf R} \Hom_{\la}(T_k,T_k) \stackrel{\bf L}{\otimes}_\Lambda T_{k+1,m}
= \la \stackrel{\bf L}{\otimes}_\Lambda T_{k+1,m} = T_{k+1,m}$.
In particular, we have $\Hom_\la(T_k,T_m)=T_{k+1,m}$.
On the other hand, we have a commutative diagram
$$\begin{array}{ccccc}
&\la&\stackrel{}{\to}&\Hom_\la(T_k,\la)&\\
&\cup&&\cup&\\
&U&\to&\Hom_\la(T_k,T_m)&\simeq T_{k+1,m},
\end{array}$$
where the horizontal map is given by $x\mapsto(\cdot x)$ for any $x\in\la$,
which is injective. Thus we have $U\subseteq T_{k+1,m}$, and so $U=T_{k+1,m}$.

Now we show the second equality.
For any $f\in\Hom_\la(\la_k,\la_m)$, there exists a unique
element $x\in\la_m$ such that $f(y)=yx$ for any $y\in\la$.
Since $T_kx\subseteq T_m$ holds, we have $x\in U$.
Thus we have $\Hom_\la(\la_k,\la_m)\simeq U/T_m=T_{k+1,m}/T_m$.
\end{proof} 
\hspace{7mm}

\subsection{Cluster tilting objects for preprojective algebras}\label{c2_sec2}
${}$ \\
Let again $\Lambda$ be the completion of the preprojective algebra of
a finite connected non-Dynkin quiver without loops over the field
$k$. We show that for a large class of
cofinite tilting ideals $I$ in $\Lambda$ we have  that $\Lambda/I$ is a finite dimensional 
$k$-algebra which is Gorenstein of dimension at most one, and the categories $\Sub\Lambda/I$ 
and $\underline{\Sub}\Lambda/I$ are stably 2-CY and 2-CY respectively. We describe some cluster 
tilting objects in these categories, using tilting ideals. We also describe cluster tilting 
subcategories in the derived 2-CY abelian category $\fl\Lambda$, which have an infinite number of 
nonisomorphic indecomposable objects. Hence we get examples of cluster structures with infinite 
clusters (see \cite{kr2} for other examples). 

We start with investigating $\Lambda/T$ for  our special cofinite tilting
ideals $T$ as a module over $\Lambda$ and over the factor ring
$\Lambda/U$ for a cofinite tilting ideal $U$ contained in $T$.

\begin{lemma}\label{lemII2.1}
Let $T$ and $U=TU'$ be cofinite tilting ideals in $\Lambda$.
Then $\Ext_{\Lambda}^1(\Lambda/T,\Lambda/U)=0=\Ext_{\Lambda}^1
(\Lambda/U,\Lambda/T)$.
\end{lemma}

\begin{proof}
Consider the exact sequence $0\to U\to \Lambda\to \Lambda/U\to 0$.
Applying $\Hom_\la(\la/T,\ )$, we have an exact sequence
\[\Ext^1_\la(\la/T,\la)\to\Ext^1_\la(\la/T,\la/U)\to\Ext^2_\la(\la/T,U).\]
We have $\Ext^1_{\la}(\Lambda/T,\Lambda)=0$ by Proposition \ref{2CY}.
It follows from Corollary \ref{propII1.10} that $\Ext_{\Lambda}^1(T,U)=0$. 
Since $\Ext_{\Lambda}^2(\Lambda/T,U)\simeq\Ext^1_\la(T,U)=0$, it follows that
$\Ext_{\Lambda}^1(\Lambda/T,\Lambda/U)=0$.
Since $\Lambda$ is derived 2-CY it follows that also $\Ext_{\Lambda}^1(\Lambda/U,\Lambda/T)=0$.
\end{proof}

Using this lemma we obtain more information on $\Lambda/T$.

\begin{proposition}\label{propII2.2}
\begin{itemize}
\item[(a)]For a cofinite ideal $T$ in $\Lambda$ with
  $\Ext^1_\la(\la/T,\la/T)=0$, the algebra $\Lambda/T$ is Gorenstein
  of dimension at most one.
\item[(b)]For a cofinite tilting ideal $T$ in $\Lambda$, the factor
  algebra $\Lambda/T$ is Gorenstein of dimension at most one.
\end{itemize}
\end{proposition}

\begin{proof}
(a) Consider the exact sequence $0\to \Omega_\Lambda( D(\Lambda/T))\to P\to
{D}(\Lambda/T)\to 0$ with a projective $\la$-module $P$. Using Lemma
\ref{lemII2.1} and \cite{ce}, we have
$\Tor^{\Lambda}_1(\Lambda/T,{D}(\Lambda/T))\simeq {D}
\Ext^1_{\Lambda^{\op}}(\Lambda/T,\Lambda/T)=0$. Applying $\Lambda/T
\otimes_{\Lambda}\ $ to the above exact sequence, we get the
exact sequence
$0\to\Lambda/T\otimes_{\Lambda}\Omega_\Lambda({D}(\Lambda/T))\to
\Lambda/T\otimes_{\Lambda} P\to
\Lambda/T\otimes_{\Lambda}{D}(\Lambda/T)\to 0$. The $\Lambda/T$-module
$\Lambda/T\otimes_{\Lambda}P$ is projective. To see that also
$\Lambda/T\otimes_{\Lambda}\Omega_\Lambda({D}(\Lambda/T))$ is a projective
$\Lambda/T$-module, we show that the functor
$\Hom_{\Lambda/T}(\Lambda/T\otimes_{\Lambda}\Omega_\Lambda({D}(\Lambda/T)),\text{
})\simeq\Hom_{\Lambda}(\Omega_\Lambda({D}(\Lambda/T)),\text{ })$ is exact
on $\mod\Lambda/T$. This follows from the functorial isomorphisms
\begin{eqnarray*}
  \Ext_{\Lambda}^1(\Omega_\Lambda({D}(\Lambda/T)), \text{
  })&\simeq&\Ext^2_{\Lambda}({D}(\Lambda/T), \text{ })\\
  &\simeq& {D}\Hom_{\Lambda}(\text{ },{D}(\Lambda/T))\\
  &\simeq& {D}\Hom_{\Lambda/T}(\text{ },{D}(\Lambda/T))\simeq \id_{\mod\la/T}
\end{eqnarray*}
Hence we conclude that $\pd_{\Lambda/T}D(\Lambda/T)\leq 1$. Then it
is  well known and easy to see that $\pd_{(\la/T)^{\op}}D(\la/T)\le 1$, so that by definition $\Lambda/T$ is Gorenstein of
dimension at most one.

(b) This is a direct consequence of (a) and Lemma \ref{lemII2.1}.\end{proof}

When $\Lambda/T$ is Gorenstein of dimension at most one, the category
of Cohen-Macaulay modules is the category $\Sub(\Lambda/T)$ of first
syzygy modules (see \cite{ar, h2}). It is known that $\Sub(\Lambda/T)$ is a Frobenius
category, with $\add(\Lambda/T)$ being the category of projective and
injective objects, and the stable category
$\underline{\Sub}(\Lambda/T)$ is triangulated \cite{h1}.
Moreover $\Sub(\Lambda/T)$ is an
extension closed subcategory of $\mod\Lambda/T$ by Corollary \ref{cor2.6}, since $\id_{\la/T}\la/T\le1$ 
and $\Ext^1_{\la/T}(\la/T,\la/T)=0$. But to show that the
stably 2-CY property is deduced from $\fl\la$ being derived 2-CY we need that $\Sub\Lambda/T$ is also
extension closed in $\fl\Lambda$.

\begin{proposition}\label{propII2.3}
  Let $T$ be a cofinite ideal with $\Ext^1_\la(\la/T,\la/T)=0$ (for example
  a cofinite tilting ideal).
  \begin{itemize}
    \item[(a)]
      $\Ext_{\Lambda}^1(\Lambda/T, X)=0=\Ext_{\Lambda}^1(X,\Lambda/T)$
      for all $X$ in $\Sub\Lambda/T$.
    \item[(b)]
      $\Sub\Lambda/T$ is an extension closed subcategory of $\fl\Lambda$
    \item[(c)]
      $\Sub\Lambda/T$ and $\underline{\Sub}\Lambda/T$  are stably 2-CY and 2-CY respectively.
  \end{itemize}
\end{proposition}

\begin{proof}
(a) For $X$ in $\Sub\la/T$ we have an exact sequence $0\to X\to P\to Y\to 0$ with  $Y$ in
  $\Sub\Lambda/T$ and $P$ in $\add\Lambda/T$. Applying
  $\Hom_{\Lambda}(\Lambda/T,\text{
  })\simeq\Hom_{\Lambda/T}(\Lambda/T,\text{ })$, the sequence does not
  change. Since $\Ext_{\Lambda}^1(\Lambda/T, \Lambda/T)=0$, we conclude that
  $\Ext_{\Lambda}^1(\Lambda/T,X)=0$. Hence
  $\Ext_{\Lambda}^1(X,\Lambda/T)=0$ by the derived 2-CY property of
  $\fl\Lambda$.

\noindent 
(b) Let $0\to X\to Y\to Z\to 0$ be an exact sequence in
$\fl\Lambda$, with $X$ and $Z$ in $\Sub\la/T$. Then we have a
monomorphism $X\to P$, with $P$ in $\add\Lambda/T$. Since
$\Ext_{\Lambda}^1(Z,P)=0$ by (a), we have a commutative diagram of
exact sequences
$$\xymatrix@C0.5cm@R0.5cm{
  0\ar[r]& X\ar@{^{(}->}[d]\ar[r]& Y\ar@{^{(}->}[d]\ar[r]&
  Z\ar@{=}[d]\ar[r]& 0\\
  0\ar[r]& P\ar[r]& P\oplus Z\ar[r]& Z\ar[r]& 0
}$$
Thus $Y$ is a submodule of $P\oplus Z\in\Sub\Lambda/T$, and we have $Y\in\Sub\Lambda/T$.

\noindent
(c) Since $\Sub\Lambda/T$ is extension closed in $\fl \la$, we have $\Ext^1_{\Sub \la/T}(X,Y)= \Ext^1_{\la}(X,Y)$. 
Since $\Sub\Lambda/T$ is Frobenius, it follows from Proposition \ref{propI1.1}, that 
$\Sub \la/T$ is stably 2-CY, since $\fl \la$ is derived 2-CY,
and so $\ul{\Sub}\Lambda/T$ is 2-CY.
\end{proof}

We now want to investigate the cluster tilting objects in
$\Sub\Lambda/T$ and $\underline{\Sub}\Lambda/T$  for certain tilting
ideals $T$, and later also the cluster tilting subcategories of
$\fl\Lambda$. The following observation will be useful.

\begin{lemma}\label{lemII2.4}
  Let $\Delta$ be a finite dimensional algebra and $M$ a
  $\Delta$-module which is a generator. Let
  $\Gamma=\End_{\Delta}(M)$, and assume $\gl\Gamma\leq 3$ and
  $\pd_{\Gamma}D(M)\leq 1$.
  Then for any $X$ in $\mod\Delta$ there is an exact sequence $0\to M_1\to M_0\to X\to 0$, 
with $M_0$ and $M_1$ in $\add M$.
\end{lemma}

\begin{proof}
  Let $X$ be in $\mod\Delta$, and consider the exact sequence $0\to
  X\to  I_0\to I_1$ where $I_0$ and $I_1$  are injective. Apply
  $\Hom_{\Delta}(M,\text{ })$ to get an exact sequence $0\to
  \Hom_{\Delta}(M,X)\to \Hom_{\Delta}(M,I_0)\to
  \Hom_{\Delta}(M,I_1)$. Since by assumption
  $\pd_{\Gamma}\Hom_{\Delta}(M,I_i)\leq 1$ for $i=0,1$ and $\gl
  \Gamma\le 3$, we obtain $\pd_{\Gamma}\Hom_{\Delta}(M,X)\le
  1$. Hence we have an exact sequence $0\to P_1\to P_0\to
  \Hom_{\Delta}(M,X)\to 0$ in $\mod\ga$  with $P_0$ and $P_1$
  projective. This sequence is the image under the functor $\Hom_{\Delta}(M,\text{
  })$ of the complex $0\to M_1\to M_0\to X\to 0$ in $\mod\Delta$, with
  $M_0$ and $M_1$  in $\add M$. Since $M$ is assumed to be a
  generator, this complex must be exact, and we have our desired exact
  sequence.
\end{proof}

Let now $\la=T_0\supsetneq T_1\supsetneq T_2\supsetneq \cdots$ be a strict
descending chain of tilting ideals corresponding to a (finite or
infinite) reduced expression $s_{i_1}s_{i_2}s_{i_3}\cdots$.
We want to describe some
natural cluster tilting objects for the algebras $\la_m=\la/T_m$. Let
$$\la_k = \la/T_k\ \mbox{ and }\ M_m=\oplus_{k=0}^m\la_k,$$
and $\ga=\End_{\la_m}(M_m)$. The following will be essential.


\begin{proposition}\label{propII2.5}
  With the above notation we have the following.
  \begin{itemize}
    \item[(a)]
      For $X$ in $\mod\la_m$ there is an exact sequence $0\to N_1\to
      N_0\to X\to 0$  in $\mod\la_m $, with $N_i$ in $\add M_m$  for
      $i=1,2$.
    \item[(b)]
      $\gl\ga\le 3$.
  \end{itemize}
\end{proposition}

\begin{proof}
We prove (a) and (b) by induction on $m$. Assume first that $m=1$. Then $\la_1=\la/T_1$, 
which is a simple  $\la_1$-module. Since $M_1 =\la/T_1$, (a) and (b) are trivially  satisfied.

Assume now that $m> 1$ and that (a) and (b) have been proved for
$m-1$. Then we first prove (b) for $m$. Note that since there are no loops for $\la$, 
we have $T_{m-1}J\subseteq T_m$ where $J$ is the Jacobson radical of
$\la$, so that $J\la_m$ is a $\la_{m-1}$-module ($\ast$). 

For an indecomposable object $X$ in $\M_m=\add M_m$, let
$f\colon C_0\to X$ be a minimal right almost split map in
$\M_m$. We first assume that $X$ is not a projective
$\la_m$-module. Then $f$ must be surjective. An
indecomposable object which is in $\M_m$ but not
in $\M_{m-1}$ is a projective $\la_m$-module, so we can write
$C_0=C_0'\oplus P$ where $C_0'\in \M_{m-1}$ and $P$ is a
projective $\la_m$-module. Since $f$ is right minimal, we have $\Ker
f\subseteq C_0'\oplus JP$, so that $\Ker f$ is a $\la_{m-1}$-module by
($\ast$). It follows by the induction assumption that there is an
exact sequence $0\to C_2\to C_1\to \Ker f\to 0$ with $C_1$ and $C_2$
in $\M_{m-1}$. Hence we have an exact sequence $0\to C_2\to
C_1\to C_0\to X\to 0$. Applying $\Hom_{\la}(M_m, \text{ })$ gives an
exact sequence
\begin{multline*}
0\to {\Hom_{\la}(M_m,C_2)}\to {\Hom_{\la}(M_m,C_1)}\to
{\Hom_{\la}(M_m,C_0)} \to {\Hom_{\la}(M_m,X)}\to S\to 0.
\end{multline*}
Then the module $S$, which is a simple module in the top of $\Hom_{\la}(M_m,X)$ in
$\mod\ga$, has projective dimension at most 3.

Assume now that $X$ is a projective $\la_m$-module. Then by  ($\ast$) we have that $JX$ is in $\mod\la_{m-1}$. 
By the induction assumption there is then an exact sequence $0\to C_1\to C_0\to J X\to 0$, with $C_0$ and $C_1$ in
$\M_{m-1}$. Hence we have an exact sequence $0\to C_1\to
C_0\to X$. Applying $\Hom_{\la}(M_m,\text{ })$ gives the exact
sequence
$$0\to \Hom_{\la}(M_m,C_1)\to \Hom_{\la}(M_m,C_0)\to
\Hom_{\la}(M_m,X)\to S\to 0$$
where $S$ is the simple top of $\Hom_{\la}(M_m,X)$, and hence $\pd_{\Gamma}S\le2$. It now follows that 
$\gl\ga_m\le 3$.

We now want to show (a) for $m$. By Proposition \ref{propII2.2} we have an
exact sequence $0\to P_1\to P_0\to {D}(\la_m)\to 0$ in $\mod\la_m$, where $P_0$ and $P_1$ are projective 
$\la_m$-modules. By
Lemma \ref{lemII2.1} we have $\Ext^1_{\la}(M_m,\la_m)=0$. Applying
$\Hom_{\la}(M_m,\text{ })$ gives the exact sequence
$$0\to \Hom_{\la}(M_m,P_1)\to \Hom_{\la}(M_m,P_0)\to
\Hom_{\la}(M_m,{ D}(\la_m))\to 0$$
Since $\Hom_{\la}(M_m,{D}(\la_m))\simeq {D}(M_m)$, we have
$\pd_{\ga_m}{D}(M_m)\le 1$. Now our desired result follows from Lemma
\ref{lemII2.4}.
\end{proof}

We can now describe some cluster tilting objects in $\Sub\la_m$  and
$\underline{\Sub}\la_m$.

\begin{theorem}\label{teoII2.6}
With the above notation, $M_m$ is a cluster tilting object in
$\Sub\la_m$ and in $\ul{\Sub}\la_m$.
\end{theorem}

\begin{proof}
We already have that $\Ext_{\la}^1(M_m,M_m)=0$ by Lemma \ref{lemII2.1},
so $\Ext_{\la_m}^1(M_m,M_m)=0$. Note that 
$\Sub\la_m=\{X\in\mod\Lambda_m \mid \Ext_{\la_m}^1(X,\la_m)=0\}$ because $\la_m$ is a cotilting module with $\id \la_m \leq 1$. Since $\la_m$ is a summand
of $M_m$, we have that $M_m$ is in $\Sub\la_m$. Assume then that
$\Ext_{\la_m}^1(X,M_m)= 0 $ for $X$ in $\mod\la_m$. By Proposition \ref{propII2.5}(a) there is
an exact sequence $0\to C_1\to C_0\to X\to 0$ with $C_1$ and $C_0$ in
$\add M_m$, which must split by our assumption. Hence $X$ is in
$\add M_m$, and it follows that $M_m$ is a cluster tilting object in
$\Sub\la_m$. It then follows as usual that it is a cluster tilting
object also in $\underline{\Sub}\la_m$.
\end{proof}

We have now obtained a large class of 2-CY categories $\Sub\la/I_w$
and $\ul{\Sub}\la/I_w$ defined via elements $w$
of the associated Coxeter group $W$, along with cluster tilting objects
associated with reduced expressions of elements in $W$. We call these
{\em standard cluster tilting objects} for $\Sub\la/I_w$ or
$\ul{\Sub}\la/I_w$. We can also describe cluster tilting
subcategories with an infinite number of nonisomorphic indecomposable
objects in the categories $\fl\la$.

\begin{theorem}\label{teoII2.7}
With the above notation, assume that each $i$ occurs an
infinite number of times in $i_1,i_2,\cdots$. Then
  $\M=\add\{\la_m \mid 0\le m\}$  is a cluster tilting subcategory of $\fl\la$.
\end{theorem}


\begin{proof}
We already know that $\Ext_{\la}^1(\la_k,\la_m)= 0 $ for all $k$ and
$m$. Let now $X$ be indecomposable in $\fl\la$. Then $X$ is a
$\la/J^k$-module for some $k$. We have $J=I_1\cap\cdots\cap I_n\supseteq
I_1 \cdots I_n$, where $1,\cdots ,n$ are the vertices in the
quiver. By our assumptions we have $J^k\supseteq T_m$ for some $m$, so
that $X$ is a $\la_m$-module. Consider the exact sequence $0\to
C_1\to C_0\to X\to 0$ in $\mod\la_m$, with $C_1$ and $C_0$  in $\add
M_m$, obtained from Proposition \ref{propII2.5}. Assume that 
$\Ext_{\la}^1(X,\M)= 0 $. Since also $\Ext_{\la_m}^1(X,\M_m)= 0 $, the sequence splits, 
so that $X$
is in $\M_m$ and hence in $\M$.

It only remains to show that $\M$ is functorially finite. So
let $X$  be in $\fl\la$. Using the above exact sequence  $0\to
C_1\to C_0\to X\to 0$, we get the exact sequence
$$0\to (C,C_1)\to (C,C_0)\to (C,X)\to \Ext_{\la}^1(C,C_1)=0$$
for $C$ in $\M$. Hence $\M$ is contravariantly finite. 

For $X$ in $\fl\Lambda$, take the left $(\Sub\M)$-approximation $X\to Y$
and choose $m$ such that $Y\in \Sub\Lambda_m$. For any $Y\in \Sub\Lambda_m$, there exists
an exact sequence $0\to Y\to C_0\to C_1\to0$ with $C_i\in\M_m$ by Proposition \ref{pro-extra}(a).
Then the composition $X\to Y\to C_0$ is a left $\M$-approximation
since $\Ext^1_\la(C_1,\M)=0$, and hence $\M$ is also covariantly finite.
\end{proof}

Summarizing our results, we have the following.
\begin{theorem}\label{CY and ct}
\begin{itemize}
\item[(a)] For any $w\in W$, we have a stably 2-CY category $\C_w=\Sub\la/I_w$.
\item[(b)] For any reduced expression $w=s_{i_1}\cdots s_{i_m}$ of
  $w\in W$, we have a cluster tilting object
  $\bigoplus_{k=1}^m\la/I_{s_{i_1}\cdots s_{i_k}}$ in $\C_w$.
  In particular the number of non-isomorphic indecomposable summands
  in any cluster tilting object is $l(w)$.
\item[(c)] For any infinite reduced expression $s_{i_1}s_{i_2}\cdots$ such that each
  $i$ occurs an infinite number of times in $i_1,i_2,\cdots$, we have
  a cluster tilting subcategory $\add\{\Lambda/I_{s_{i_1}\cdots s_{i_k}} \mid 0\le k\}$ in $\fl\la$.
\end{itemize}
\end{theorem}

We end this section by showing that the subcategories $\Sub\la/I_w$ can be
characterized using torsionfree classes.

\begin{theorem}
Let $\la$ be the completed preprojective algebra of a connected non-Dynkin
quiver without loops. Let $\C$ be a torsionfree class in $\fl\la$ with
some cluster tilting object. Then we have $\C=\Sub\la/I_w$ for some element
$w$ in the Coxeter group associated with $\Lambda$.
\end{theorem}

\begin{proof}
We first prove that if $M$ is a cluster tilting object in $\C$, then
$\C=\Sub M$. We only have to show $\C\subset\Sub M$.
For any $X\in\C$, take a projective resolution
$\Hom_\la(M,N)\to\Ext^1_\la(M,X)\to 0$ ($*$) of $\End_\la(M)$-modules
with $N\in\add M$.
Replacing $M$ in ($*$) by $N$, we get an exact sequence $0 \to X \to
Y \to N \to 0$ as the image of the identity $1_N\in\End_\la(N)$ in $\Ext^1_\la(N,X)$.
Since $\C$ is extension closed, we have $Y\in \C$. 
We have an exact sequence
$\Hom_\la(M,N) \to \Ext^1_\la(M,X) \to \Ext^1_\la(M,Y) \to
\Ext^1_\la(M,N)=0$.
Since ($*$) is exact, we have $\Ext^1_\la(M,Y)=0$. 
Thus we have $Y\in\add M$ and $X\in\Sub M$.

Let now $I$ be the annihilator $\ann{}_\la M$ of $M$ in $\la$. 
Then $I$ is clearly a cofinite ideal in $\la$, and $\ann{}_{\la/I}M=0$.
Further $\Sub M$ is extension closed also in $\mod\la/I$. Hence the direct
sum $A$ of one copy of each of the non-isomorphic indecomposable
$\Ext$-injective $\la/I$-modules in $\Sub M$ is a cotilting
$\la/I$-module satisfying 
$\id{}_{\la/I}A\le 1$ and $\Sub M=\Sub A$ by \cite{Sm}.
Since $\Sub M$ is extension closed in the derived 2-CY
category $\fl\la$, the $\Ext$-injective $\la/I$-modules in $\Sub M$
coincide with the $\Ext$-projective ones, which are the projective
$\la/I$-modules. Hence we have that $A$ is a progenerator of $\la/I$
and $\Sub M=\Sub\la/I$. Since $\Sub\la/I$ is extension closed in $\fl\la$,
we have $\Ext^1_\la(\la/I,\la/I)=0$.

By Theorem \ref{teoII1.6}, we only have to show that $I$ is a partial
tilting left ideal. By Bongartz completion, we only have to show $\Ext^1_\la(I,I)=0$.
The natural surjection $\la\to\la/I$ clearly induces a surjection
$\Hom_\la(\la/I,\la/I)\to\Hom_\la(\la,\la/I)$. Since $\la$ is
derived 2-CY, we have injections
$\Ext^2_\la(\la/I,\la)\to\Ext^2_\la(\la/I,\la/I)$ and
$\Ext^1_\la(I,\la)\to\Ext^1_\la(I,\la/I)$.
Using the exact sequence $0\to I\to\la\to\la/I\to0$, we have a
commutative diagram
\[\begin{array}{cccccc}
&&\Ext^1_\la(\la/I,\la/I)=0\\
&&\uparrow\\
\Hom_\la(I,\la)&\to&\Hom_\la(I,\la/I)&\to&\Ext^1_\la(I,I)&\to
\Ext^1_\la(I,\la)\to\Ext^1_\la(I,\la/I)\\
\uparrow&&\uparrow&&\uparrow\\
\la&\to&\la/I&\to&0
\end{array}\]
of exact sequences. Thus we have $\Ext^1_\la(I,I)=0$.
\end{proof}

Note that we have proved that an extension closed subcategory of $\fl\la$
of the form $\Sub X$ for some $X$ in $\fl\la$ with some cluster tilting
object must be $\Sub\la/I_w$ for some element $w$ in the Coxeter group
associated with $\Lambda$.

We point out that there are other extension closed subcategories of $\fl\la$
with some cluster tilting object.
Let $Q$ be an extended Dynkin quiver and $Q'$ a Dynkin subquiver,
and $\la$ and $\la'$ the corresponding completed preprojective algebras.
Then clearly $\mod\la'=\fl\la'$ is an extension closed subcategory of
$\fl\la$. Hence any extension closed subcategory of $\mod\la'$ is extension
closed in $\fl\la$, so Example 1 in Section \ref{c1_sec3}
is an example of an extension closed stably 2-CY subcategory of $\fl\la$
with some cluster tilting object, but which is not closed under submodules.\\
\hspace{7mm}

\subsection[Realzations of categories]{Realization of cluster categories and stable categories for preprojective 
algebras of Dynkin type}\label{c2_sec3}
${}$ \\
In this section we show that for an appropriate choice of $T$ as a
product of tilting ideals $I_j=\la(1-e_j)\la$, any cluster category
is equivalent to some $\underline{\Sub}\la/T$. In particular, any cluster
category can be realized as the stable category of a Frobenius
category with finite dimensional homomorphism spaces. We also show that the stable 
categories for preprojective algebras of Dynkin type can be realized this way.

Let $Q$ be a finite connected quiver without loops, $kQ$ the associated 
path algebra, and $\la $ the completion of the preprojective algebra of $Q$. Choose a complete set of orthogonal primitive 
idempotents $e_1,...,e_n$ of $kQ$. We can assume that $e_i(kQ)e_j=0$
for any $i>j$. We regard $e_1,...,e_n$ as a complete set of orthogonal primitive 
idempotents of $\Lambda$, and let as before $I_i =\Lambda(1-e_i)\Lambda$.

Assume first that $Q$ is not Dynkin. 
We consider an exact stably 2-CY category associated to the square
$w^2$ of a Coxeter element $w=s_1s_2\cdots s_n\in W$.
Let $\la_i=\la/I_1I_2\cdots I_i$ and
$\la_{i+n}=\la/I_1I_2\cdots I_nI_1\cdots I_i$ for $1\le i\le n$.
We have seen in Section \ref{c2_sec2} that $\Sub\la_{2n}$, and also $\ul{\Sub}\la_{2n}$, has a cluster tilting object 
$M=\oplus_{i=1}^{2n}\la_i$. 

We shall need the following.

\begin{lemma}\label{lemII3.1}
Assume that $Q$ is not Dynkin.
Then $I_1\cdots I_nI_1\cdots I_nI_1\cdots$ gives rise to a
 strict descending chain of tilting ideals.
In particular, $s_1\cdots s_ns_1\cdots s_ns_1\cdots$ is reduced.
\end{lemma}

\begin{proof}
Assume to the contrary that the descending chain of ideals is not strict. Let 
$T_i=I_1 \cdots I_i$ and $U_i=I_1\cdots I_{i-1}I_{i+1}\cdots I_n$ for $i=1,\cdots, n$. 
Then we have $T_n^kT_{i-1}=T_n^kT_i$ for some $i=1,\cdots,n$ and $k\ge0$, where $T_0=\la$. 
Hence we obtain $T_n^{k+1}=T_n^kU_i$. 
Then we get $T_n^{k+m}=T_n^kU_i^m$ for any $m>0$. Since $U_i e_i =  \la e_i$ and $J\supseteq T_n$, 
we have $J^{m+k}e_i \supseteq T_n^{m+k}e_i=T_n^kU^m_ie_i = T_n^ke_i$. Since $(\la/T_n^k)e_i$ has finite 
length, we have $J^{m+k}e_i=T_n^ke_i$ for $m$ sufficiently large. Thus we have $T_n^ke_i=0$, which is a 
contradiction  since $\la e_i$ has infinite length.
We have the latter assertion from Proposition \ref{strict}.
\end{proof}


We have the following.

\begin{proposition}\label{teoII3.2}
Let $Q$ be a finite connected non-Dynkin quiver without oriented cycles
and with vertices $1, \cdots, n$ ordered as above.
Let $\la_{2n}=\la/{(I_1\cdots I_n)^2}$. Then $\la_n=\la/{I_1\cdots I_n}$ 
is a cluster tilting object in $\ul{\Sub}\la_{2n}$ with
$\End_{\ul{\Sub}\la_{2n}}(\la_n)\simeq kQ$.
\end{proposition}

\begin{proof}
Since the associated chain of ideals is strict descending by Lemma \ref{lemII3.1}, our general theory applies. 
We have a cluster tilting object $\bigoplus_{i=1}^{2n}\la_i$ in $\ul{\Sub}\la_{2n}$ by Theorem \ref{teoII2.6}.
We have $\add\bigoplus_{i=1}^{2n}\la_i=\la_n\oplus\la_{2n}$ in
$\fl\la$ by Proposition \ref{strict}.
Thus $\la_n$ is a cluster tilting object in $\ul{\Sub}\la_{2n}$.


Note that the path algebra $kQ$ is, in a natural way, a factor algebra of $\la$, and $kQ$ is hence a $\la$-module. 
We want to show that the $\la$-modules $\la_n$ and $kQ$ are isomorphic.

Let $P_j$ be the indecomposable projective $\la$-module corresponding to the vertex $j$. 
Then $I_{j+1}\cdots I_nP_j=P_j$ and $I_jP_j=JP_j$, the smallest submodule of $P_j$ such that 
the corresponding factor has only composition factors $S_j$. Further, $I_{j-1}I_jP_j$ 
is the smallest submodule of $I_jP_j=JP_i$ such that the factor has only composition 
factors $S_{j-1}$, etc. By our choice of ordering, this means that the paths starting at $j$, 
with decreasing indexing on the vertices, give a basis for $P_j/I_1\cdots I_nP_j$. 
In other 
words, we have $P_j/I_1\cdots I_nP_j\simeq(kQ)e_j$. Hence the $\la$-modules $\la_n=\la/I_1\cdots I_n$ 
and $kQ$ are isomorphic, so that $\End_{\la}(\la_n)\simeq kQ$.

It remains to show $\End_{\la_{2n}}(\la_n)\simeq\End_{\ul{\Sub}\la_{2n}}(\la_n)$.
By Lemma \ref{hom-calculation}, any morphism
from $\la_{n}$ to $\la_{2n}$ is given by a right multiplication of an
element in $(I_1\cdots I_n)/(I_1\cdots I_n)^2$.
This implies
$\Hom_\la(\la_n,\la_{2n})\Hom_\la(\la_{2n},\la_n)=0$.
Thus we have the assertion.
\end{proof} 


We now show that we have the same kind of result for Dynkin quivers.

\begin{proposition}\label{teoII3.3}
Let $Q'$ be a Dynkin quiver with vertices $1,\cdots,m$ contained in a finite
connected non-Dynkin quiver $Q$ without oriented cycles and with vertices $1,\cdots,n$ ordered as before. 
Let $\la$ be the preprojective algebra of $Q$ and 
$\la_{n+m}=\la/(I_1\cdots I_nI_1\cdots I_m)$. Then 
$\la_m=\la/ I_1\cdots I_m$ is a cluster tilting object in $\ul{\Sub}\la_{n+m}$
with $\End_{\ul{\Sub}\la_{n+m}}(\la_m)\simeq kQ'$.
\end{proposition}

\begin{proof}
Since we have seen in Lemma \ref{lemII3.1} that the product $(I_1\cdots I_n)^2$ gives
rise to a strict descending chain of ideals,
it follows that the same holds for $I_1\cdots I_nI_1\cdots I_m$.
The assertions follow as in the proof of Proposition \ref{teoII3.2}.
\end{proof}

Recall from \cite{kr2} that if a connected algebraic triangulated 2-CY category has a cluster tilting object $M$  
whose quiver $Q$ has no oriented cycles, then $\C$ is triangle equivalent to the cluster 
category $\C_{kQ}$. Then we get the following consequence of the last two results.

\begin{theorem}\label{teoII3.4}
Let $Q'$ be a finite connected quiver without oriented cycles. Let
$Q=Q'$ if $Q'$ is not Dynkin, and $Q$ as in
Proposition \ref{teoII3.3} if $Q'$ is Dynkin.
Let $\la$ be the preprojective algebra of $Q$.
Then there is a tilting ideal $I$ in $\la$ such that $\ul{\Sub}\la/I$
is triangle equivalent to the cluster category $\C_{kQ'}$ of $Q'$.
\end{theorem}

We finally show that also the categories $\ul{\mod}\la'$ where $\la'$ is the preprojective 
algebra of a Dynkin quiver $Q'$, can be realized this way.


\begin{theorem}\label{teoII3.5}
Let $Q'$ be a Dynkin quiver contained in a finite connected non-Dynkin
quiver $Q$ without loops. We denote by $\Lambda'$ the preprojective algebra of $Q'$,
by $W'$ the subgroup of $W$ generated by $\{s_i \mid i\in Q_0'\}$,
and by $w_0$ the longest element in $W'$.
Then $\Lambda'$ is isomorphic to $\Lambda/I_{w_0}$
and $\ul{\mod}\la'=\ul{\Sub}\Lambda/I_{w_0}$.
\end{theorem}

\begin{proof}
Let $I_{Q'}:=\la(\sum_{i\in Q_0\backslash Q_0'}e_i)\la$. 
Since we have $\la/I_{Q'}\simeq\la'$, we only have to show $I_{w_0}=I_{Q'}$.
We use the fact that $I_{Q'}$ is maximal amongst all two-sided ideals
$I$ of $\la$ such that any composition factor of $\la/I$ is $S_i$ for
some $i\in Q_0'$.

Since $w_0$ is a product of $s_i$ ($i\in Q_0'$), 
any composition factor of $\la/I_{w_0}$ is $S_i$ for some $i\in Q_0'$.
Thus we have $I_{w_0}\supseteq I_{Q'}$.
On the other hand, since $w_0$ is the longest element of $W'$,
we have $l(s_iw_0)<l(w_0)$ for any $i\in Q_0'$.
By Proposition \ref{length}, we have $I_iI_{w_0}=I_{w_0}$ for any $i\in Q_0'$.
This implies $I_{w_0}=I_{Q'}$.
\end{proof}

Using Theorem \ref{teoII3.5}, we see that our theory also applies to
preprojective algebras of Dynkin type. In particular, we can specialize
Theorem \ref{CY and ct} to recover the following result from \cite{gls1}.

\begin{corollary}
For a preprojective algebra $\la'$ of a Dynkin quiver the number of
non-isomorphic indecomposable summands in a cluster tilting object
is equal to the length $l(w_0)$ of the longest element in the associated
Weyl group, which is equal to the number of positive roots.
\end{corollary}

We also obtain a large class of cluster tilting objects associated with
the different reduced expressions of $w_0$.

Our results can also be viewed as giving an interpretation in terms of
tilting theory of some functors $\mathcal{E}_i$ used in \cite[5.1]{gls3}.\\
\hspace{7mm}

\subsection{Quivers of cluster tilting subcategories}
${}$ \\
In this section we show that the quivers of standard cluster tilting
subcategories associated with a reduced expression can be described
directly from the reduced expression.

Let $s_{i_1}s_{i_2}\cdots s_{i_k}\cdots$ be a (finite or infinite)
reduced expression associated with a graph $\Delta$ with vertices
$1,\cdots,n$.
We associate with this sequence a quiver $Q(i_1,i_2,\cdots)$ as follows,
where the vertices correspond to the $s_{i_k}$.
\begin{itemize}
\item[-] For two consecutive $i$ ($i\in\{1,\cdots,n\}$), draw an arrow
  from the second one to the first one.
\item[-] For each edge $i\stackrel{d_{ij}}{-}j$, pick out the expression
  consisting of the $i_k$ which are $i$ or $j$, so that we have
  $\cdots i i \cdots i j j \cdots j i i \cdots i \cdots$.
We draw $d_{ij}$ arrows from the last $i$ in a connected set of $i$'s
to the last $j$ in the next set of $j$'s, and the same from $j$ to
$i$. (Note that since by assumption both $i$ and $j$ occur an infinite
number of times if the expression is infinite, each connected set of $i$'s
or $j$'s is finite.)
\end{itemize}

Note that in the Dynkin case essentially the same quiver has been used in
\cite{fz3}.

For a finite reduced expression $s_{i_1}\cdots s_{i_k}$ we denote by
$\ul{Q}(i_1,\cdots i_k)$ the quiver obtained from $Q(i_1,\cdots i_k)$
by removing the last $i$ for each each $i$ in $Q_0$.

We denote by $\Lambda=T_0\supsetneq T_1\supsetneq\cdots$ the associated
strict descending chain of tilting ideals.
Then we have a cluster tilting
subcategory $\M(i_1,i_2,\cdots)=\add\{\Lambda_k \mid k>0\}$
for $\Lambda_k:=\Lambda/T_k$.

\begin{theorem}\label{quiver}
Let the notation be as above.
\begin{itemize}
\item[(a)]The quiver of the cluster tilting
subcategory $\M(i_1,i_2,\cdots)$ is $Q(i_1,i_2,\cdots)$.
\item[(b)]The quiver of $\ul{\End}_\la(\M(i_1,\cdots i_k))$ is
  $\ul{Q}(i_1,\cdots i_k)$.
\end{itemize}
\end{theorem}

Before we give the proof, we give some examples and consequences.

It follows from the definition that we get the same quiver if we
interchange two neighbors in the expression of $w$ which are not
connected with any edge in $\Delta$. But if we take two reduced expressions in
general, we may get different quivers, as the following examples show.

Let $\Delta$ be the graph
$\xymatrix@R0.1cm@C0.5cm{&2\ar@{-}[dr]\\ 1\ar@{-}[ur]&&3\ar@{-}[ll]}$
and $w=s_1s_2s_1s_3s_2=s_2s_1s_2s_3s_2=s_2s_1s_3s_2s_3$ expressions
which are clearly reduced.
The first expression gives the quiver
$$
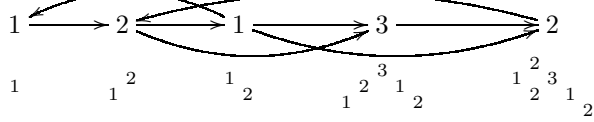
\xymatrix@R0.1cm{1\ar[r]&2\ar[r]\ar@/_1pc/[rr]&
1\ar@/_1pc/[ll]\ar@/_1pc/[rr]\ar[r]&3\ar[r]&2\ar@/_1pc/[lll]\\
{\begin{smallmatrix}1\end{smallmatrix}}&
{\begin{smallmatrix}&2\\ 1&\end{smallmatrix}}&
{\begin{smallmatrix}1&\\ &2\end{smallmatrix}}&
{\begin{smallmatrix}&&3&&\\ &2&&1&\\ 1&&&&2\end{smallmatrix}}&
{\begin{smallmatrix}&2&&&\\ 1&&3&&\\ &2&&1&\\ &&&&2\end{smallmatrix}}}$$
the second one gives the quiver
$$
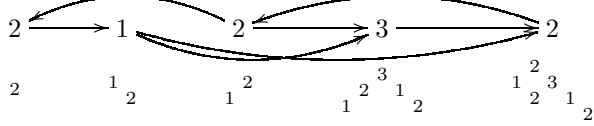
\xymatrix@R0.1cm{2\ar[r]&1\ar@/_1pc/[rrr]\ar@/_1pc/[rr]&
2\ar@/_1pc/[ll]\ar[r]&3\ar[r]&2\ar@/_1pc/[ll]\\
{\begin{smallmatrix}2\end{smallmatrix}}&
{\begin{smallmatrix}1&\\ &2\end{smallmatrix}}&
{\begin{smallmatrix}&2\\ 1&\end{smallmatrix}}&
{\begin{smallmatrix}&&3&&\\ &2&&1&\\ 1&&&&2\end{smallmatrix}}&
{\begin{smallmatrix}&2&&&\\ 1&&3&&\\ &2&&1&\\ &&&&2\end{smallmatrix}}}$$
and the third one gives the quiver
$$
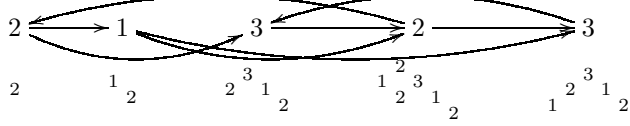
\xymatrix@R0.1cm{2\ar[r]\ar@/_1pc/[rr]&1\ar@/_1pc/[rr]\ar@/_1pc/[rrr]&
3\ar[r]&2\ar@/_1pc/[lll]\ar[r]&3\ar@/_1pc/[ll]\\
{\begin{smallmatrix}2\end{smallmatrix}}&
{\begin{smallmatrix}1&\\ &2\end{smallmatrix}}&
{\begin{smallmatrix}&3&&\\ 2&&1&\\ &&&2\end{smallmatrix}}&
{\begin{smallmatrix}&2&&&\\ 1&&3&&\\ &2&&1&\\ &&&&2\end{smallmatrix}}&
{\begin{smallmatrix}&&3&&\\ &2&&1&\\ 1&&&&2\end{smallmatrix}}}$$

We now investigate the relationship between the cluster tilting
objects given by different reduced expressions of the same element.

\begin{lemma}\label{replace}
Let $w=s_{i_1}\cdots s_{i_m}=s_{i'_1}\cdots s_{i'_m}$ be reduced
expressions and $\Lambda=T_0\supsetneq T_1\supsetneq\cdots$ and
$\Lambda=T_0'\supsetneq T_1'\supsetneq\cdots$ corresponding tilting ideals.
\begin{itemize}
\item[(a)] Assume that for some $k$ we have $i_k=i'_{k+1}$,
  $i'_k=i_{k+1}$ and $i_j=i'_j$ for any 
  $j\neq k,k+1$. Then the corresponding cluster tilting objects are isomorphic.
\item[(b)]  Assume that for some $k$ we have $i_{k-1}=i'_{k}=i_{k+1}$,
  $i'_{k-1}=i_{k}=i'_{k+1}$ and $i_j=i'_j$ for any 
  $j\neq k,k\pm1$. Then the corresponding cluster tilting objects are
  in the relationship of exchanges of $T_{k-1}e_{i_{k-1}}$ and $T'_{k-1}e_{i'_{k-1}}$.
\end{itemize}
\end{lemma}

\begin{proof}
(a) Obviously we have $T_j=T'_j$ for any $j<k$. Since
$s_{i_k}s_{i_{k+1}}=s_{i'_k}s_{i'_{k+1}}$, we have $I_{i_k}\
I_{i_{k+1}}=I_{i'_k}I_{i'_{k+1}}$.
Thus we have $T_j=T_j'$ for any $j>k+1$. In particular, we have
$T_je_{i_j}=T'_je_{i'_j}$ for any $j\neq k,k+1$. Since
$I_{i_k}e_{i_k}=I_{i'_k}I_{i'_{k+1}}e_{i'_{k+1}}$, we have
$T_ke_{i_k}=T'_{k+1}e_{i'_{k+1}}$. 
Similarly, we have $T_{k+1}e_{i_{k+1}}=T'_ke_{i'_k}$.
Thus the assertion follows.

(b) Since $I_{i_{k-1}}I_{i_{k}}I_{i_{k+1}}=I_{i'_{k-1}}I_{i'_{k}}I_{i'_{k+1}}$, we
have $T_je_{i_j}=T'_je_{i'_j}$ for any $j\neq k,k\pm1$. 
Since $I_{i_{k-1}}I_{i_{k}}e_{i_{k}}=I_{i'_{k-1}}I_{i'_{k}}I_{i'_{k+1}}e_{i_{k+1}}$,
we have $T_{k}e_{i_{k}}=T'_{k+1}e_{i'_{k+1}}$. Similarly we have
$T_{k+1}e_{i_{k+1}}=T'_{k}e_{i'_{k}}$. Thus we have the assertion.
\end{proof}

As an illustration, note that in the above example we obtain the
second quiver from the first by mutation at the left vertex.
Immediately we have the following conclusion.

\begin{proposition}\label{transitivity}
All cluster tilting objects in $\Sub(\Lambda/I_w)$ obtained from
reduced expressions of $w$ can be obtained from each other under
repeated exchanges.
\end{proposition}

\begin{proof}
This is immediate from Lemma \ref{replace} and \cite[3.3.1]{bb} since
we get from one reduced expression to another by applying the
operations described in Lemma \ref{replace}.
\end{proof}

Using Theorem \ref{teoII3.5} we see that for preprojective algebras of
Dynkin quivers we get the quivers of the endomorphism algebras associated
with reduced expressions of the longest element $w_0$.

For a stably 2-CY category or a triangulated 2-CY category $\C$ with
cluster tilting subcategories we have an associated {\em cluster tilting graph}
defined as follows. The vertices correspond to the non-isomorphic basic
cluster tilting objects, and two vertices are connected with an edge if
the corresponding cluster tilting objects differ in exactly one
indecomposable summand. For cluster categories this graph is known to be
connected \cite{bmrrt}, while this is an open problem in general.
For the categories $\Sub\la/I_w$ or $\ul{\Sub}\la/I_w$ it follows from
Proposition \ref{transitivity} that all standard cluster tilting objects
belong to the same component of the cluster tilting graph, and we call
this the {\em standard component}.

\medskip
We now illustrate with some classes of examples.
Let $Q$ be a connected non-Dynkin quiver without oriented cycles,
and vertices $1,\cdots,n$, where there is no arrow $i\to j$ for $j>i$.

(a) Let $w=s_1s_2\cdots s_ns_1s_2\cdots s_n$. The last $n$ vertices
correspond to projectives, so the quiver for the cluster tilting object
in the stable category $\ul{\Sub}\la/I_w$ is $Q$, which has no
oriented cycles. So we get an alternative proof of Proposition \ref{teoII3.2}.

(b) Choose $w=s_1s_2\cdots s_ns_1s_2\cdots s_n\cdots$.
Ordering the indecomposable preprojective modules as
$P_1,\cdots,P_n,\tau^{-1}P_1,\cdots,\tau^{-1}P_n,\cdots,\tau^{-i}P_1,\cdots,
\tau^{-i}P_n,\cdots$
where $P_i$ is the projective module associated with vertex $i$,
we have a bijection between the indecomposable preprojective modules
and the terms in the expression for $w$. Then the quiver of
the corresponding cluster tilting subcategory is the
preprojective component of the AR quiver of $kQ$, with an additional
arrow from $X$ to $\tau X$ for each indecomposable preprojective module
$X$. This is a direct consequence of our rule, since we know by Lemma
\ref{lemII3.1} that the expression for $w$ is reduced.


(c) Now take a part $\mathcal{P}$ of the AR quiver of the
preprojective component, closed under predecessors. Consider the
expression obtained from $s_1s_2\cdots s_ns_1s_2\cdots s_n\cdots$.
by only keeping the terms corresponding to the objects in
$\mathcal{P}$ under our given bijection. We show below that this new
expression is reduced. Then it follows directly by our rule that
when adding arrows $X\to\tau X$ when $X$ is nonprojective in $\mathcal{P}$, we get
the quiver of the cluster tilting object given by the above reduced
expression. That this quiver is the quiver of a cluster tilting object
was also shown in \cite{gls1b} for $\mathcal{P}$ being the AR quiver of a Dynkin
quiver, and in \cite{gls5} in the general case.

\begin{lemma}
The word associated with $\mathcal{P}$ obtained in this way is reduced.
\end{lemma}

\begin{proof}
The word satisfies the following conditions.
\begin{itemize}
\item[(a)] For each pair $(i,j)$ of vertices connected by some edge,
  $i$ and $j$ occur each second time after removing the other vertices. 
\item[(b)] $w=A_1A_2\cdots A_t$, where 
each $A_s$ is a strict increasing sequence of numbers in $\{1,\cdots,n\}$,
such that, if $j\notin A_s$, then $j\notin A_{s+1}$,
and if $i<j$ are connected with an edge and
$i\notin A_s$ (respectively, $j\notin A_s$), then $j\notin A_s$
(respectively, $i\notin A_{s+1}$).
\end{itemize}
The condition (a) is immediate from the construction.
For each pair $(i,j)$ connected with some edge, we
have in the AR quiver
$$
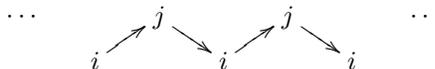
\xymatrix@R0.1cm@C0.5cm{\cdots&&j\ar[dr]&&j\ar[dr]&&\cdots\\
&i\ar[ur]&&i\ar[ur]&&i&}$$
the part involving $i$'s and $j$'s. Hence they must occur each second
time to give this quiver. Thus the condition (a) is satisfied.

If $A_1\neq(1,\cdots,n)$, then $w$ is a subsequence of
$(1,\cdots,n)$. So the word is clearly reduced in this case.
Thus we can assume $A_1=(1,\cdots,n)$.
We show that, for any word satisfying the conditions (a)(b) and
$A_1=(1,\cdots,n)$, the corresponding descending chain $\Lambda=T_0\supseteq
T_1\supseteq\cdots$ is strict. Then the word is reduced by Proposition \ref{strict}.

We assume $T_{k-1}=T_{k}$ for some $k$.
So $I_{i_1}\cdots I_{i_k}=I_{i_1}\cdots I_{i_{k-1}}$.
Take $s$ minimal such that $A_s\neq(1,\cdots,n)$, and take $i$ minimal
such that $i\notin A_s$. By assumption (b), all terms
appearing after the position of $i$ in $A_s$ are not connected with
$i$ by an edge in $\Delta$. In particular, the corresponding ideals
commute with $I_i$. By multiplying with $I_i$ from the right and using
commutative relations, we get an equality where $i$ is inserted in
$A_s$ after $1,\cdots,i-1$. Repeating this process, we get an equality
$$I_1\cdots I_nI_1\cdots I_n\cdots I_1\cdots I_{i_k}=I_1\cdots I_n
I_1\cdots I_n\cdots I_1\cdots I_{i_k-1}.$$
This contradicts Lemma \ref{lemII3.1}.
\end{proof}


\bigskip
In the rest of this section we give a proof of Theorem \ref{quiver}.
Note that (b) follows directly from (a) and Proposition \ref{strict}(b).

(a) Let $J$ be the Jacobson radical of $\Lambda$.
Let $\M=\M(i_1,i_2,\cdots)$ and $T_{l,k}=I_{i_l}\ I_{i_{l+1}}\cdots I_{i_k}$.
For $l>k$, this means $T_{l,k}=\Lambda$.
In the rest we often use the following equalities.
$$e_i\ J\ e_{i'}=\left\{\begin{array}{cc}\ e_i\ I_i\ e_{i'}&(i=i')\\
e_i\ \Lambda\ e_{i'}&(i\neq i'),\end{array}\right.\ \ 
I_i\ e_{i'}=\left\{\begin{array}{cc}J\ e_{i'}&(i=i')\\
\Lambda\ e_{i'}&(i\neq i'),\end{array}\right.\ \mbox{ and }\ 
e_{i'}\ I_i=\left\{\begin{array}{cc}e_{i'}\ J&(i=i')\\
e_{i'}\ \Lambda&(i\neq i').\end{array}\right.$$
We have
$$\Hom_\Lambda(\Lambda_l\ e_{i_l},\ \Lambda_k\ e_{i_k})=e_{i_l}\
(T_{l+1,k}/T_k)\ e_{i_k}$$
by Lemma \ref{hom-calculation}.
We have
$$\rad_{\M}(\Lambda_l\ e_{i_l},\ \Lambda_k\ e_{i_k})=
e_{i_l}\ (T_{(l+1-\delta_{l,k}),k}/T_k)\ e_{i_k}.$$ 
Moreover, we have
$$\rad^2_{\M}(\Lambda_l\ e_{i_l},\ \Lambda_k\ e_{i_k})=e_{i_l}\
((T_k+\sum_{j>0}T_{(l+1-\delta_{l,j}),j}\ e_{i_j}\
T_{(j+1-\delta_{j,k}),k})/T_k)\ e_{i_k}.$$ 
To get the quiver of $\M$, we have to compute 
$(\rad_{\M}/\rad_{\M}^2)(\Lambda_l\ e_{i_l},\ \Lambda_k\ e_{i_k})=E_{l,k}/D_{l,k}$ for
$$E_{l,k}=e_{i_l}\ T_{(l+1-\delta_{l,k}),k}\ e_{i_k}\supseteq
D_{l,k}=e_{i_l}\ T_k\ e_{i_k}+\sum_{j>0}e_{i_l}\ T_{(l+1-\delta_{l,j}),j}\
e_{i_j}\ T_{(j+1-\delta_{j,k}),k}\ e_{i_k}.$$ 
We denote by $k^+$ the minimal number satisfying $k<k^+$ and
$i_k=i_{k^+}$ if it exists.

(i) We consider the case when there are no arrows in $Q$ from $l$ to $k$.
We shall show $E_{l,k}=D_{l,k}$.

If $l>k$ and $i_l=i_k$, then we have $l>k^+>k$. Thus
$$E_{l,k}=e_{i_l}\ \Lambda\ e_{i_k}=e_{i_l}\ \Lambda\ e_{i_{k^+}}\
\Lambda\ e_{i_k}\subseteq D_{l,k}.$$

In the rest we assume that either $l\le k$ or $i_l\neq i_k$ holds.
First we show
$$E_{l,k}=e_{i_l}\ T_{l+1,k-1}\ (1-e_{i_k})\ \Lambda\ e_{i_k}=\sum_{a\neq
  i_k}e_{i_l}\ T_{l+1,k-1}\ e_a\ \Lambda\ e_{i_k}$$
by the following case by case study. 
\begin{itemize}
\item[-] If $l<k$, then
$E_{l,k}=e_{i_l}\ T_{l+1,k-1}I_{i_k}\ e_{i_k}=e_{i_l}\ T_{l+1,k-1}\
(1-e_{i_k})\ \Lambda\ e_{i_k}$.
\item[-] If $l=k$, then $E_{l,k}=e_{i_l}\ I_{i_k}\ e_{i_k}=e_{i_l}\
  \Lambda\ (1-e_{i_k})\ \Lambda\ e_{i_k}$.
\item[-] If $l>k$ and $i_l\neq i_k$, then
$E_{l,k}=e_{i_l}\ \Lambda\ e_{i_k}=e_{i_l}\ I_{i_k}\ e_{i_k}=e_{i_l}\
\Lambda\ (1-e_{i_k})\ \Lambda\ e_{i_k}$.
\end{itemize}

Thus we only have to show that $e_{i_l}\ T_{l+1,k-1}\ e_a\ \Lambda\
e_{i_k}\subseteq D_{l,k}$
for any $a\neq i_k$. We have the following three possibilities.
\begin{itemize}
\item[-] If $a\notin\{i_1,i_2,\cdots,i_{k-1}\}$, then we have
$T_{l+1,k-1}\ e_a=I_{i_{l+1}}\cdots I_{i_{k-1}}\ e_a=\Lambda
e_a=I_{i_1}\cdots I_{i_{k-1}}\ e_a=T_{k-1}e_a$. Thus
$$e_{i_l}\ T_{l+1,k-1}\ e_a\ \Lambda\ e_{i_k}=e_{i_l}\ T_{k-1}\ e_a\
\Lambda\ e_{i_k}\ \subseteq
e_{i_l}\ T_k\ e_{i_k}\subseteq D_{l,k}.$$
\item[-] If $a\notin\{i_{k+1},i_{k+2},\cdots,i_{k^+-1}\}$, then we have
$e_a\ \Lambda=e_a\ I_{i_k}\cdots I_{i_{k^+}}=e_a\ T_{k,k^+}$. Thus
$$e_{i_l}\ T_{l+1,k-1}\ e_a\ \Lambda\ e_{i_k}=
e_{i_l}\ T_{l+1,k-1}\ e_a\ T_{k,k^+}\ e_{i_k}\subseteq
e_{i_l}\ T_{l+1,k^+}\ e_{i_k}=e_{i_l}\ T_{l+1,k^+}\ e_{i_{k^+}}\
\Lambda\ e_{i_k}\subseteq D_{l,k}$$
since $l\neq k^+\neq k$.
\item[-] Otherwise, there is an arrow $j\to k$ with $i_j=a$.
Since $a\notin\{i_{j+1},i_{j+2},\cdots,i_{k-1}\}$,
we have $T_{j+1,k-1}\ e_a\ \Lambda=I_{i_{j+1}}\cdots I_{i_{k-1}}\ e_a\
\Lambda=\Lambda\ e_a\ \Lambda=\Lambda\ e_a\ I_{i_{j+1}}\cdots
I_{i_{k-1}}=\Lambda\ e_aT_{j+1,k-1}$. Thus 
$$e_{i_l}\ T_{l+1,k-1}\ e_a\ \Lambda\ e_{i_k}=e_{i_l}\ T_{l+1,j}\
T_{j+1,k-1}\ e_a\ \Lambda\ e_{i_k}=e_{i_l}\ T_{l+1,j}\ e_{i_j}\
T_{j+1,k}\ e_{i_k}\subseteq D_{l,k}$$ 
since $l\neq j\neq k$.

\end{itemize}
In each case we have $e_{i_l}\ T_{l+1,k-1}\ e_a\ \Lambda\
e_{i_k}\subseteq D_{l,k}$ for any $a\neq i_k$. 

(ii) We consider the case $l=k^+$.
We have $E_{k^+,k}=e_{i_{k^+}}\ \Lambda\ e_{i_k}$.
We shall show $D_{k^+,k}=e_{i_l}\ J\ e_{i_k}$.

Clearly we have $D_{k^+,k}\subseteq e_{i_{k^+}}\ J\ e_{i_k}$.
Conversely, we have
$$e_{i_{k^+}}\ J\ e_{i_k}=e_{i_{k^+}}\ I_{i_k}\ e_{i_k}=e_{i_{k^+}}\
\Lambda\ (1-e_{i_k})\ \Lambda\ e_{i_k}=\sum_{a\neq i_k}e_{i_{k^+}}\
\Lambda\ e_a\ \Lambda\ e_{i_k}.$$
Thus we only have to show $e_{i_{k^+}}\ \Lambda\ e_a\ \Lambda\
e_{i_k}\subseteq D_{k^+,k}$ 
for any $a\neq i_k$. We have the following two possibilities.
\begin{itemize}
\item[-] If $a\notin\{i_1,i_2,\cdots,i_{k-1}\}$, then we have
$\Lambda e_a=I_{i_1}\cdots I_{i_{k-1}}e_a=T_{k-1}\ e_a$. Thus
$$e_{i_{k^+}}\ \Lambda\ e_a\ \Lambda\ e_{i_k}=e_{i_{k^+}}\ T_{k-1}\
e_a\ I_{i_k}\ e_{i_k}\subseteq e_{i_{k^+}}\ T_k\ e_{i_k}\subseteq D_{k^+,k}.$$
\item[-] If $a\in\{i_1,i_2,\cdots,i_{k-1}\}$, then take the largest $j$
such that $i_j=a$. Then we have
$\Lambda=T_{k^++1,j}$ and $e_a\ \Lambda=e_a\ I_{i_{j+1}}\cdots
I_{i_k}=e_a\ T_{j+1,k}$. Thus
$$e_{i_{k^+}}\ \Lambda\ e_a\ \Lambda\ e_{i_k}=e_{i_{k^+}}\
T_{k^++1,j}\ e_{i_j}\ T_{j+1,k}\ e_{i_k}\subseteq
D_{k^+,k}$$
since $k^+\neq j\neq k$.
\end{itemize}
In each case we have $e_{i_{k^+}}\ \Lambda\ e_a\ \Lambda\
e_{i_k}\subseteq D_{k^+,k}$ for any $a\neq i_k$. 

(iii) Finally we consider the case when $l\neq k^+$ and there is an
arrow in $Q$ from $l$ to $k$. Then $l<k$.
We have $E_{l,k}=e_{i_l}\ J\ e_{i_k}$.
We shall show $D_{l,k}=e_{i_l}\ J^2\ e_{i_k}$.

First we show $D_{l,k}\subseteq e_{i_l}\ J^2\ e_{i_k}$.
We have $e_{i_l}\ T_k\ e_{i_k}\subseteq e_{i_l}\ I_{i_l}\ I_{i_k}\
e_{i_k}=e_{i_l}\ J^2\ e_{i_k}$.
We have the following three possibilities.
\begin{itemize}
\item[-] Assume $l\le j\le k$. Then
$e_{i_l}\ T_{(l+1-\delta_{l,j}),j}\ e_{i_j}\ T_{(j+1-\delta_{j,k}),k}\
e_{i_k}\subseteq
e_{i_l}\ J\ e_{i_j}\ J\ e_{i_k}\subseteq e_{i_l}\ J^2\ e_{i_k}$.

\item[-] Assume $k<j$. If $i_j\neq i_k$, then
$e_{i_l}\ T_{l+1,j}\ e_{i_j}\ T_{j+1,k}\ e_{i_k}\subseteq
e_{i_l}\ I_{i_j}\ e_{i_j}\ J\ e_{i_k}\subseteq e_{i_l}\ J^2\ e_{i_k}$.
If $i_j=i_k$, then $j\ge k^+$. Since there is an arrow $l\to k$, we
have $i_l\in\{i_{k+1},i_{k+2},\cdots,i_{k^+-1}\}\subseteq\{i_{l+1},i_{l+2},\cdots,i_{j-1}\}$. 
Thus $e_{i_l}\ T_{l+1,j}\ e_{i_j}\ T_{j+1,k}\ e_{i_k}\subseteq
e_{i_l}\ I_{i_l}\ I_{i_j}\ e_{i_j}\ \Lambda\ e_{i_k}\subseteq e_{i_l}\
J^2\ e_{i_k}$.

\item[-] Assume $l>j$. If $i_j\neq i_l$, then
$e_{i_l}\ T_{l+1,j}\ e_{i_j}\ T_{j+1,k}\ e_{i_k}\subseteq
e_{i_l}\ J\ e_{i_j}\ I_k\ e_{i_k}\subseteq e_{i_l}\ J^2\ e_{i_k}$.
If $i_j=i_l$, then $e_{i_l}\ T_{l+1,j}\ e_{i_j}\ T_{j+1,k}\ e_{i_k}\subseteq
e_{i_l}\ \Lambda\ e_{i_j}\ I_{i_l}\ I_{i_k}\ e_{i_k}\subseteq e_{i_l}\
J^2\ e_{i_k}$.
\end{itemize}

Next we show $e_{i_l}\ J^2\ e_{i_k}\subseteq D_{l,k}$. We have
$$e_{i_l}\ J^2\ e_{i_k}=e_{i_l}\ J\ (1-e_{i_k})\ \Lambda\ 
e_{i_k}=\sum_{a\neq i_k}e_{i_l}\ J\ e_a\ \Lambda\ e_{i_k}.$$
Thus we only have to show $e_{i_l}\ J\ e_a\ \Lambda\ e_{i_k}\subseteq D_{l,k}$
for any $a\neq i_k$. We have the following two possibilities.
\begin{itemize}
\item[-] If $a\notin\{i_1,i_2,\cdots,i_{k-1}\}$, then we have
$e_{i_l}\ J\ e_a=e_{i_l}\ \Lambda\ e_a=e_{i_l}\ I_{i_1}\cdots
I_{i_{k-1}}\ e_a=e_{i_l}\ T_{k-1}\ e_a$. Thus
$$e_{i_l}\ J\ e_a\ \Lambda\ e_{i_k}=e_{i_l}\ T_{k-1}\ e_a\ I_{i_k}\ 
e_{i_k}\subseteq e_{i_l}\ T_k\ e_{i_k}\subseteq D_{l,k}.$$
\item[-] If $a\in\{i_{1},i_{2},\cdots,i_{k-1}\}$, then take the largest $j$
  such that $i_j=a$. Then we have
$e_a\ \Lambda=e_a\ I_{i_{j+1}}\cdots I_{i_k}=e_a\ T_{j+1,k}$. 
Moreover if $j=l$, then $e_{i_l}\ J\ e_a=e_{i_l}\
T_{(l+1-\delta_{l,j}),j}\ e_a$.
If $j\neq l$, then $e_{i_l}\ J\ e_a\subseteq e_{i_l}\ I_{i_{l+1}}\cdots
I_{i_j}\ e_a=e_{i_l}\ T_{l+1,j}\ e_a$. Thus
$$e_{i_l}\ J\ e_a\ \Lambda\ e_{i_k}=e_{i_l}\ T_{(l+1-\delta_{l,j}),j}\
e_{i_j}\ T_{j+1,k}\ e_{i_k}\subseteq D_{l,k}$$
since $j\neq k$.
\end{itemize}
In each case we have $e_{i_l}\ J\ e_a\ \Lambda\ e_{i_k}\subseteq
D_{l,k}$ for any $a\neq i_k$. 
\qed \\
\hspace{7mm}

\subsection{Substructure}\label{c2_sec4}
${}$ \\
In this section we point out that the work in this chapter gives several illustrations of 
substructures of cluster structures. We also  give some concrete examples of 2-CY categories and 
their cluster tilting objects, to be applied in Chapter \ref{chap3}.

Let $s_{i_1}s_{i_2}\cdots s_{i_t}\cdots$ be an infinite reduced expression
which contains each $i\in\{1,\cdots,n\}$ an infinite number of times. 
Let $T_t=I_{i_1}\cdots I_{i_t}$, and $\la_t=\la/T_t$. Recall that for $t<m$, we have 
$\Sub\la_t \subseteq\Sub\la_m\subseteq \fl\la$. We then have the following.

\begin{theorem}\label{teoII4.1}
Let the notation be as above.
\begin{itemize}
\item[(a)]$\Sub\la_m$, $\ul{\Sub}\la_m$ and $\fl\la$ have a cluster structure using the 
cluster tilting subcategories with the indecomposable projectives as coefficients.
\item[(b)]{For $t<m$, the cluster tilting object
    $\la_1\oplus\cdots\oplus\la_t$ in $\Sub \la_t$ can be extended to
    a cluster tilting object
    $\la_1\oplus\cdots\oplus\la_t\oplus\cdots\oplus\la_m$ in $\Sub
    \la_m$, and determines a substructure of $\Sub\la_m$.} 
\item[(c)]{The cluster tilting object $\la_1\oplus\cdots\oplus\la_t$
    in $\Sub \la_t$ can be extended to the cluster tilting subcategory
    $\{\la_i \mid i\ge 0\}$ in $\fl\la$, and determines a substructure
    of $\fl\la$.} 
\end{itemize}
\end{theorem}
\begin{proof}
(a) Since $\Sub\la_m$ and $\ul{\Sub}\la_m$ are stably 2-CY and 2-CY respectively, they have a weak cluster 
structure. It follows from Proposition \ref{prop2.5} that we have no 
loops or 2-cycles, using the cluster tilting objects. Then it follows
from Theorem \ref{teoI1.6} that we have a 
cluster structure for $\Sub\la_m$ and $\ul{\Sub}\la_m$.

That also $\fl\la$ has a cluster structure follows by using that it is the case for all the $\Sub\la_m$.

\noindent
(b) and (c) follow directly from the definition of substructure and previous results.
\end{proof}

We now consider the Kronecker quiver 
$\xymatrix@R0.2cm{
  1\ar@<0.5ex>[r]\ar@<-0.5ex>[r]& 0
}$, 
and let $\la$ be the associated preprojective algebra. The only strict descending chains are 
$$I_0\supsetneq I_0I_1\supsetneq I_0I_1I_0\supsetneq \cdots(I_0I_1)^j\supsetneq(I_0I_1)^jI_0\supsetneq\cdots 
\text{   \hspace{0.5cm}and }$$
 $$I_1\supsetneq I_1I_0\supsetneq I_1I_0I_1\supsetneq \cdots(I_1I_0)^j\supsetneq(I_1I_0)^jI_1\supsetneq\cdots$$
We let $T_t$ be the product of the first $t$ ideals, and $\la_t=\la/T_{t}$. Both $I_0$ and $I_1$  
occur an infinite number of times in each chain. The indecomposable projective $\la$-modules $P_0$ and $P_1$ 
have the  following structure\\

$P_0=
\begin{smallmatrix}
  && 0 &&\\
  &1&&1&\\
  0&&0&&0\\
  && \cdot &&\\
  && \cdot &&\\
  && \cdot &&
\end{smallmatrix}
$\hspace{0.5cm}
$P_1=
\begin{smallmatrix}
  && 1 &&\\
  &0&&0&\\
  1&&1&&1\\
  && \cdot &&\\
  && \cdot &&\\
  && \cdot &&
\end{smallmatrix}
$ \hspace{0.5cm} \\\\where radical layer number $2i$ has $2i$ copies of 1 for $P_0$ and $2i$ copies of 0 for $P_1$, 
and radical layer number $2i+1$ has $2i+1$ copies of 0 for $P_0$ and $2i+1$ copies of 1 for $P_1$. We write 
$P_{0,t}=P_0/J^tP_0$ and $P_{1,t}=P_1/J^tP_1$. Then it is 
easy to see that for the chain $I_0\supsetneq I_0I_1\supsetneq\cdots $ we have $\la_1=\la/I_0=P_{0,1}=(0)$,  
$\la_2=\la/I_0I_1=P_{0,1}\oplus P_{1,2}=(0)\oplus\left(
\begin{smallmatrix}
  & 1&\\
  0&&0\\
\end{smallmatrix}
\right)$, $\la_3=\la/I_0I_1I_0=P_{0,3}\oplus P_{1,2}$,..., $\la_{2t}=P_{0,2t+1}
\oplus P_{1,2t}$, $\la_{2t+1}=P_{0,2t-1}\oplus P_{1,2t}$,..

Note that this calculation also shows that both our infinite chains are strict descending.

It follows from Section 4 (and is also easily seen directly) that the
quiver of the cluster tilting subcategory  $\{\la_i \mid i\ge 1\}$ is
the following:
$$
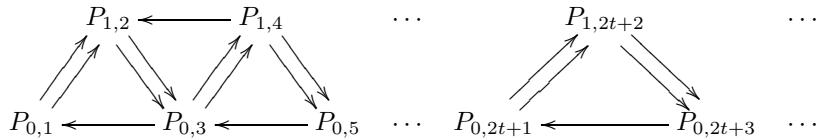
\xymatrix@C0.2cm{
  & P_{1,2}\ar@<0.1cm>[dr]\ar@<-0.1cm>[dr]&& P_{1,4}\ar[ll]\ar@<0.1cm>[dr]\ar@<-0.1cm>[dr]&
&\cdots&& P_{1,2t+2}\ar@<0.1cm>[dr]\ar@<-0.1cm>[dr]&& \cdots\\
  P_{0,1}\ar@<0.1cm>[ur]\ar@<-0.1cm>[ur]&& P_{0,3}\ar[ll]\ar@<0.1cm>[ur]\ar@<-0.1cm>[ur]&& P_{0,5}\ar[ll]&\cdots& 
P_{0,2t+1}\ar@<0.1cm>[ur]\ar@<-0.1cm>[ur]&& P_{0,2t+3}\ar[ll]& \cdots
}$$
In particular, we have the cluster tilting object $P_{0,1}\oplus P_{1,2}\oplus P_{0,3}$ for $\Sub\la_3$, where the 
last two summands are projective. Hence $P_{0,1}$ is a cluster tilting object in $\ul{\Sub}\la_3$. 
The quiver of the endomorphism 
algebra consists of one vertex and no arrows. Hence $\ul{\Sub}\la_3$ is equivalent to the cluster 
category $\C_k$, 
which has exactly two indecomposable objects. The other one is $JP_{0,3}$, obtained from the exchange sequence 
$0\to JP_{0,3}\to P_{0,3}\to P_{0,1}\to 0$. Note that it is also easy to see directly that there are no other 
indecomposable rigid nonprojective objects in $\Sub\la_3$.

For $\la_4$, we have the cluster tilting object $P_{0,1}\oplus P_{1,2}\oplus P_{0,3}\oplus P_{1,4}$ for $\Sub\la_4$. 
Again the last two $\la_4$-modules are projective, so $P_{0,1}\oplus P_{1,2}$ is a cluster 
tilting object in $\ul{\Sub}\la_4$. 
The quiver of the endomorphism algebra is 
$\xymatrix@C0.8cm{
  \cdot \ar@<0.5ex>[r]\ar@<-0.5ex>[r]&\cdot
}$
, which has no oriented cycles, and hence $\ul{\Sub}\la_4$ is triangle equivalent to the cluster category $\C_{k(
\xymatrix@C0.4cm{
   \cdot\ar@<0.4ex>[r]\ar@<-0.4ex>[r]&\cdot
})}$.
In particular the cluster tilting graph is connected.
We can use this to get a description of the rigid objects in $\Sub\la_4$.
\begin{proposition}
  Let $\la_4=\la/I_0I_1I_0I_1$ be the algebra defined above. Then the 
indecomposable rigid $\la_4$-modules in $\Sub\la_4$ are 
exactly the ones of the form $\Omega^i_{\la_4}(P_{0,1})$ and $\Omega^i_{\la_4}(P_{1,2})$ 
for $i\in \mathbb{Z}$
\end{proposition}
\begin{proof}
  For $\C_{k(
\xymatrix@C0.4cm{
   \cdot\ar@<0.4ex>[r]\ar@<-0.4ex>[r]&\cdot
})}$ the indecomposable rigid objects are the $\tau$-orbits of the objects induced by the 
indecomposable projective $k(
\xymatrix@C0.4cm{
   \cdot\ar@<0.4ex>[r]\ar@<-0.4ex>[r]&\cdot
})$-modules. Here $\tau=[1]$, and for $\ul{\Sub}\la_4$, $\Omega^{-1}=[1]$. This proves the claim.
\end{proof}

The cluster tilting graph for $\Sub\la_3$ and $\Sub\la_4$ are
$\xymatrix@C0.4cm@R0.4cm{
\cdot \ar@{-}[r] & \cdot}$ and 
$\xymatrix{\cdots\text{ }\cdot \ar@{-}[r] & \cdot \ar@{-}[r] & \cdot \ar@{-}[r] & \cdot \ar@{-}[r] &\cdot \text{ }\cdots}$.

We end with the following problem.

\begin{conjecture}
For any $w\in W$ the cluster tilting graph for $\Sub\la/I_w$ is connected.
\end{conjecture}


\section{Connections to cluster algebras}\label{chap3}

While the theory of 2-CY categories is interesting in itself, one of the motivations 
for investigating 2-CY categories comes from the theory of cluster algebras initiated by 
Fomin and Zelevinsky \cite{fz1}. In many situations the 2-CY categories can be used to construct 
new examples of cluster algebras, and also to give a new categorical model for
already known examples.
This has been done in for example \cite{ck1,ck2} and \cite{gls1}. In this chapter
we illustrate with some applications of the theory developed in the first two chapters. In
Section 2, we recall the definition and basic properties of 
a map $\varphi$ from the finite dimensional modules over the completed preprojective algebra
of a connected quiver with no loops
to the function field $\mathbb{C}(U)$ of the associated unipotent group
$U$. This was used in \cite{gls1} to model the cluster algebra structure of
the coordinate ring $\mathbb{C}[U]$ in the Dynkin case. It is what we call a (strong) cluster map. We also make explicit the
notion of subcluster algebra, and observe that 
a substructure of a (stably) 2-CY category together with a cluster map gives rise to a subcluster
algebra. This gives one way to construct new cluster algebras inside
$\mathbb{C}[U]$ in the Dynkin case, or model old ones, as we illustrate with examples in Section 2.
In Section 3 we deal with the non-Dynkin case. We conjecture that the stably 2-CY category
$\Sub \la/I_w$ discussed in Chapter \ref{chap2} gives a model for a cluster algebra structure on the
coordinate ring $\mathbb{C}[U^w]$ of the corresponding unipotent cell.
We give examples which support this conjecture.\\
\hspace{7mm}

\subsection{Cluster algebras, subcluster algebras and cluster maps}\label{c3_sec1}
${}$ \\
In this section we recall the notion of cluster algebras
\cite{fz1} and make explicit a notion of {\em subcluster algebras}.
Actually we extend the definition of cluster algebras to include the possibility of 
clusters with countably many elements. The coordinate rings of unipotent groups of non-Dynkin diagrams
are candidates for containing such cluster algebras.
We also introduce certain maps, called (strong) {\em cluster maps},
defined for categories with a cluster structure.
The image of a cluster map gives rise to a cluster algebra.
Cluster substructures on the category side give rise to subcluster algebras.

We first recall the definition of a cluster algebra, allowing countable clusters. 
Note that the setting used here is not the most general one.
Let $m \geq n$ be positive integers, or countable numbers.
Let $\F= \Q(u_1, \dots, u_m)$ be the field of rational functions over
$\Q$ in $m$ independent variables. A cluster algebra is a subring 
of $\F$, constructed in the following way.
A {\em seed} in $\F$ is a triple $(\ul{x}, \ul{c}, \widetilde{B})$,
where $\ul{x}$ and $\ul{c}$ are non-overlapping sets of elements in $\F$, where
we let $\widetilde{\ul{x}} = \ul{x} \cup \ul{c}$ and sometimes denote the seed
by the pair $(\widetilde{\ul{x}}, \widetilde{B})$. 
Here $\widetilde{\ul{x}} = \{x_1, \dots, x_m \}$ should be a transcendence basis for $\F$
and $\widetilde{B} = (b_{ij})$ is a locally finite $m \times n$-matrix 
with integer elements such that the submatrix $B$ of $\widetilde{B}$ consisting of the $n$ 
first rows is skew-symmetric. 

The set $\ul{x} = \{x_1, \dots, x_n \}$ is called the {\em cluster} of the seed,
and the set $\ul{c} = \{x_{n+1}, \dots, x_m \}$ 
is the {\em coefficient set} of the cluster algebra. The
set $\ul{\tilde{x}} = \ul{x} \cup \ul{c}$ is called an {\em extended cluster}.

For a seed $(\ul{\tilde{x}},\widetilde{B})$, with $\widetilde{B} =
(b_{ij})$, and for $k \in \{ 1, \dots,  n \}$,
a {\em seed mutation} in direction $k$ produces a new seed 
$(\ul{\tilde{x}}',\widetilde{B}')$. Here
$\ul{\tilde{x}}' = (\ul{\tilde{x}} \setminus \{x_k \}) \cup \{x_k' \}$,
where 
$$
x_k' = x_k^{-1} (\prod_{b_{ik} > 0} x_i^{b_{ik}} + \prod_{b_{ik} < 0} x_i^{-b_{ik}} )
$$
This is called an {\em exchange relation} and $\{ x_k, x_k'\}$ is called an {\em exchange pair}.
Furthermore
$$
b_{ij}' = \begin{cases} -b_{ij} & \text{if $i=k$ or $j=k$} \\ 
b_{ij} + \frac{\left| b_{ik} \right| b_{kj} + b_{ik} \left| b_{kj} \right|}{2} & \text{else.}
\end{cases}
$$

Fix an (initial) seed $(\ul{\tilde{x}}, \widetilde{B})$,
and consider the set $\Sc$ of all seeds obtained from $(\ul{\tilde{x}},\widetilde{B})$
by a sequence of seed mutations. The union $\X$ of all elements in the 
clusters in $\Sc$ are called the {\em cluster variables}, and for a fixed subset of coefficients 
$\ul{c}_0 \subseteq \ul{c}$,
the {\em cluster algebra} $\A(\Sc)$ with the coefficients $\ul{c}_0$ inverted, is the $\ZZ[\ul{c}, {\ul{c}_0}^{-1}]$-subalgebra 
of $\F$ 
generated by $\X$. Note that we, unlike in the original definition, do not
necessarily invert all coefficients. This is done in order to catch examples like
the coordinate ring of a maximal unipotent group in the Dynkin case 
and the homogeneous coordinate ring
of a Grassmannian.
Note that we often extend the scalars for cluster algebras to $\mathbb{C}$.

We now make explicit the notion of subcluster algebras.
Let $\A$ be a cluster algebra with cluster variables $\X$, coefficients $\ul{c}$,
and ambient field $\F = \Q(u_1, \dots , u_m)$. 
A {\em subcluster algebra} $\A'$ of $\A$ is a cluster algebra such that
there exists a seed $(\ul{x},\ul{c},Q)$ for $\A$ and a seed
$(\ul{x}',\ul{c}',Q')$ for $\A'$ such that 
\begin{itemize}
\item[(S1)]{$\ul{x}'\subseteq \ul{x}$ and $\ul{c}' \subseteq \ul{x} \cup \ul{c}$.}
\item[(S2)]{For each cluster variable $x_i\in\ul{x}'$, the set of
    arrows entering and leaving $i$ in $Q$ lie in $Q'$.}
\item[(S3)]{The invertible coefficients $\ul{c}_0' \subseteq \ul{c}'$ satisfy 
$\underline{c}_0 \cap \underline{c}' \subseteq \underline{c}_0'$.}
\end{itemize}

Note that a subcluster algebra is not necessarily a subalgebra since
some coefficients may be inverted. Also note that $\A'$ is determined by
the seed $(\ul{x}',\ul{c}',Q')$ and the set $\ul{c}_0'$ of invertible
coefficients.

The definition implies that clusters in the subcluster algebra can be uniformly extended.

\begin{proposition}\label{lemIII1.1} 
\begin{itemize}
\item[(a)]Seed mutation in $\A'$ is compatible with seed mutation in $\A$.
\item[(b)]There is a set $\ul{v}$ consisting of cluster variables and coefficients in $\A$, such that
for any extended cluster $\ul{x}'$ in $\A'$, $\ul{x}' \cup \ul{v}$ is
an extended cluster in $\A$.
\end{itemize}
\end{proposition}

\begin{proof}
This follows directly from the definition.
\end{proof}

%
%

Inspired by \cite{gls1,gls2} and \cite{cc, ck1,ck2} we introduce certain maps, which we call 
(strong) cluster maps,
defined for a 2-CY category with a (weak) cluster structure, 
and such that the image gives rise to a cluster algebra. We show that such
maps preserve substructures, as defined above and in Section \ref{c1_sec2}.
 
Recall that a category $\C$ is stably $2$-CY if it is either an exact Frobenius category
where $\ul{\C}$ is $2$-CY or a functorially finite extension closed subcategory
$\B$ of a triangulated $2$-CY category $\C$.

Let $\C$ be a stably 2-CY category with a 
cluster structure defined by cluster tilting objects, where projectives are 
coefficients.
We assume that the cluster tilting objects have $n$ cluster variables and $c$ coefficients,
where $1 \leq n \leq \infty$ and $0 \leq c \leq \infty$. 
For a cluster tilting object $T$, we denote by $ B_{\End_{\C}(T)}$ the $m \times n$-matrix
obtained by removing the last $m-n$ columns of
the skew-symmetric $m \times m$ matrix corresponding
to the quiver of the endomorphism algebra $\End_{\C}(T)$, where
the columns are ordered such that those corresponding to projective summands of $T$ come last. 
We can also think of this as dropping the arrows between vertices in $\End_{\C}(T)$ corresponding
to indecomposable projective summands of $T$ from the quiver of $\End_{\C}(T)$.

Let $\F = \Q(u_1, \dots, u_m)$.
Given a connected component $\Delta$ of the cluster tilting graph of $\C$,
a {\em cluster map} (respectively, {\em strong cluster map}) for
$\Delta$ is a map $\varphi \colon 
\E=\add\{T\mid T\in\Delta\} \to \F$ (respectively, $\varphi \colon 
\C \to \F$) where isomorphic objects have the same image, satisfying
the following three conditions.
\begin{itemize}
\item[(M1)]{For a cluster tilting object $T$ in $\Delta$, $\varphi(T)$
    is a transcendence basis for $\F$.}
\item[(M2)]{(respectively, (M2$'$)) For all indecomposable objects $M$
    and $N$ in $\E$ (respectively, $\C$) with 
    $\dim_k\Ext^1(M,N)=1$, we 
    have $\varphi(M)\varphi(N) = \varphi(V) + \varphi(V')$ where $V$
    and $V'$ are the middle of the non-split triangles/short exact
    sequences $N\to V\to M$ and $M\to V'\to N$.}
\item[(M3)]{(respectively, (M3$'$)) $\varphi(A \oplus A')= \varphi(A)
    \varphi(A')$ for all $A,A'$ in $\E$ (respectively, $\C$).}
\end{itemize}

Note that a pair $(M,N)$ of indecomposable objects in $\E$ is an
exchange pair if and only if $\Ext^1(M,N) \simeq k$ (see \cite{bmrrt}).
Note that a map $\varphi\colon \C\to \F$ satisfying
(M2$'$) and (M3$'$) is called a {\em cluster character} in \cite{p}.
Important examples of (strong) cluster maps appear in \cite{ck1,ck2,
  gls1}, and more recently in \cite{gls5,p}.


\begin{theorem}\label{propIII2.2} 
With $\C$ and $\E$ as above, let $\varphi \colon \E \to \F$ be a cluster map. 
Then the following hold.
\begin{itemize}
\item[(a)]{Let $\A$ be the subalgebra of $\F$ generated by
    $\varphi(X)$ for $X\in\E$. Then $\A$ is a cluster algebra and
    $(\varphi(T),B_{\End_{\C}(T)})$ is a seed for $\A$ for any cluster
    tilting object $T$ in $\Delta$.}
\item[(b)]{Let $\B$ be a subcategory of $\C$ with a substructure, and
    $\E'$ a subcategory of $\B$ defined by a connected component of
    the cluster tilting graph of $\B$. Then $\varphi(X)$ for $X\in\E'$
    generates a subcluster algebra of $\A$.}
\end{itemize}
\end{theorem}

\begin{proof}
(a) follows from the fact that $\C$ has a cluster structure. For (b), 
let $T'$ be a cluster tilting object in the subcategory that extends to a cluster tilting object
$T$ for $\E$. Then $\varphi(T')$ gives a transcendence basis for a subfield $\F'$ of $\F$,
and $T'$ together with its matrix $B_{\End_{\C}(T')}$ gives a seed for a subcluster algebra.
\end{proof}

For any subset $\ul{c}_0$ of coefficients of $\A$, we have a cluster
algebra $\A[\ul{c}_0^{-1}]$. We say that {\em the cluster algebra
  $\A[\ul{c}_0^{-1}]$ is modelled by the cluster map $\varphi:\E\to\F$}.\\
\hspace{7mm}

\subsection{The GLS $\varphi$-map with applications to the Dynkin case}\label{c3_sec2}
${}$ \\
Let $Q$ be a finite connected quiver without loops, and let $\la$ be
the associated completed preprojective algebra over an algebraically
closed field $k$, and $W$ the associated Coxeter group. For the non-Dynkin
case, we have
the derived 2-CY category $\fl \Lambda$ of finite dimensional left
$\Lambda$-modules, and for each $w$ in $W$ the stably 2-CY category
$\Sub \la/I_w$ as investigated in Chapter \ref{chap2}.

On the other hand, associated with the underlying graph of $Q$, is a Kac-Moody
group $G$ \cite{kp} with a maximal unipotent subgroup $U$.
Let $H$ be the torus. Recall that the Weyl group
$\operatorname{Norm}(H)/H$ is isomorphic to the Coxeter group $W$
associated with $Q$.
For an element $w \in W=\operatorname{Norm}(H)/H$ and 
any lifting $\widetilde{w}$ of $w$ in $G$, define the unipotent
cell $U^w$ \cite{BZ} to be the intersection
\[U^w = U \cap B_- \widetilde{w} B_-,\]
where $B_-$ is the opposite Borel subgroup corresponding to $U$. Then
$U^w$ is independent of the choice of lifting of $w$. It is a
quasi-affine algebraic variety of dimension
$l(w)$ and we have $U = \bigsqcup_{w \in W} U^w$.

Let $U(\mathfrak{n})$ be enveloping algebra
of the maximal nilpotent subalgebra
of the Kac-Moody Lie algebra $\mathfrak{g}$
associated to $Q$ and
let $U(\mathfrak{n})^*$ and $U(\mathfrak{n})^*_{\text{gr}}$
be the dual and graded dual of
$U(\mathfrak{n})$ respectively.
Note that both $U(\mathfrak{n})^*$ and
$U(\mathfrak{n})^*_{\text{gr}}$ become
algebras with respect to the $\circ$-product which is dual
to the coproduct on $U(\mathfrak{n})$.

Recall that the matrix coefficient function $f^{\tau}_{v,\zeta}$
associated to a representation $\tau:
U(\mathfrak{n}) \longrightarrow \mathfrak{gl}(V)$
and vectors $v \in V$ and $\zeta \in V^*$ is the linear
form in $U(\mathfrak{n})^*$
defined by $f^{\tau}_{v,\zeta}(x) = \zeta
\Big( \tau(x) \cdot v \Big)$.

Define the restricted dual $U(\mathfrak{n})^*_{\text{res}}$ of
the enveloping algebra $U(\mathfrak{n})$
to be the span of the unit element $1 \in U(\mathfrak{n})^*$
together with all matrix coefficient functions
$f^{\tau}_{v,v^*}$ for integrable lowest weight
representations
$\tau \colon U(\mathfrak{n}) \longrightarrow \mathfrak{gl}(V)$
with $v \in V$ and $\zeta \in V^*$.
The restricted dual $U(\mathfrak{n})^*_{\text{res}}$ is a
subalgebra of $U(\mathfrak{n})^*$ due to the fact
that $f^{\tau}_{v,\zeta} \circ f^{\tau'}_{v', \zeta'}
\ = \ f^{\tau \otimes \tau'}_{ v \otimes v', \zeta \otimes \zeta'}$.
One can check that the graded dual $U(\mathfrak{n})^*_{\text{gr}}$
is a subalgebra of the restricted dual $U(\mathfrak{n})^*_{\text{res}}$.

Let $\mathbb{C}[U]$ denote the ring of (strongly)
regular functions on the unipotent group $U$ as defined
in \cite{kp2}. 
Given a representation $\varrho \colon U \longrightarrow \GL(V)$
of the unipotent group $U$ we may differentiate it to obtain
a $U(\mathfrak{n})$-representation $d \varrho \colon U(\mathfrak{n})
\longrightarrow \mathfrak{gl}(V)$ uniquely determined by the
formula

\[ d {\varrho}(x) \cdot v = \ { {\displaystyle \partial} \over
{\displaystyle \partial t} } \Big|_{t=0} \, \exp(tx) \cdot v \]

\bigskip
\noindent
for all $x \in \mathfrak{n}$. Accordingly a representation
$\varrho \colon U \longrightarrow \GL(V)$
will be called an integrable lowest weight representation
of the unipotent group $U$ if the differentiated
representation $d \varrho \colon U(\mathfrak{n}) \longrightarrow
\mathfrak{gl}(V)$ is integrable and lowest weight.

Theorem 1 of \cite{kp2} implies that
$\mathbb{C}[U]$ is spanned by the unit element together with
all matrix coefficient functions
$F^{\varrho}_{v,\zeta}$ for
integrable lowest weight representations $\varrho: U \longrightarrow
\GL(V)$ with $v \in V$ and $\zeta \in V^*$. The function
$F^{\varrho}_{v,\zeta}: U \longrightarrow
\Bbb{C}$ is strongly regular and is defined for $g \in U$
by $F^{\varrho}_{v,\zeta}(g) := \zeta \Big( \varrho(g) \cdot v \Big)$.
One easily checks that, with regard to the usual product of functions,
the identity $F^{\varrho}_{v,\zeta} F^{\varrho'}_{v', \zeta'}
\ = \ F^{\varrho \otimes \varrho'}_{ v \otimes v', \zeta \otimes \zeta'}$
holds.

In view of these facts it follows that the mapping $\iota$ given by
\[ f^{d\varrho}_{v, \zeta} \longmapsto F^{\varrho}_{v,\zeta} \]
defines an algebra isomorphism between the restricted dual
$U(\mathfrak{n})^*_{\text{res}}$ and the coordinate ring
$\mathbb{C}[U]$. In particular $\iota$ restricts to an
embedding of the graded dual $U(\mathfrak{n})^*_{\text{gr}}$
into $\mathbb{C}[U]$.

\bigskip
We now apply the construction of the GLS $\varphi$-map \cite{gls1} to the non-Dynkin case.
For $M\in\fl\Lambda$, Geiss-Leclerc-Schr\"oer \cite{gls5} constructed
$\delta_M\in U(\mathfrak{n})^*_{\text{gr}}$ by using Lusztig's Lagrangian
construction of $U(\mathfrak{n})$ \cite{Lu1,Lu2}.
We denote by $\varphi(M)\in\mathbb{C}[U]$ the image of $\delta_M$
under the above map $\iota:U(\mathfrak{n})^*_{\text{gr}}\to \mathbb{C}[U]$.
Then $\varphi(M)$ has the following property.

Let 
$x_i(a) = \operatorname{exp}(ae_i)$
for $a \in \mathbb{C}, i \in Q_0$ be
elements in $U$, for a Chevalley generator $e_i$. 
For ${\bf i} = \big(i_1, \dots, i_k \big) \in \mathbb{Z}_{\geq 0}^k$
the symbol {\LARGE $\chi_{\bf \s i}$}$(M)$ denotes the 
Euler characteristic of the variety $\mathcal{F}_{\bf i}(M)$ of 
all flags in $M=M_0 \supset M_1 \supset \cdots \supset M_k = 0$
with $M_{l-1}/M_l$ isomorphic to the simple module $S_{i_l}$. 
Then for any ${\bf i} = \big(i_1, \dots, i_k \big) \in \mathbb{Z}_{\geq 0}^k$
and $a_1,\cdots,a_k\in\mathbb{C}$, we have
\begin{equation}\label{eq-gls}
\varphi(M)(x_{i_1}(a_1) \cdots x_{i_k}(a_k)) \ = \ \sum_{ {\bf j} \in \mathbb{Z}_{\geq 0}^k} \ 
\text{\huge $\chi$}_{\rev {\bf i}^{\bf j} } (M) \ \frac{{a_1^{j_1} \cdots a_k^{j_k}}}{
{j_1! \cdots j_k!}}               
\end{equation}
where $ \rev {\bf i}^{\bf j}$ is the
$j_1 + \cdots + j_k$ tuple which, when read from left to right, starts
with $j_k$ occurrences of $i_k$, followed by $j_{k-1}$ occurrences of $i_{k-1}$,
and ultimately $j_1$ occurrences $i_1$.

Notice that one can prove that property \eqref{eq-gls} uniquely
determines $\varphi(M)$.
The following result (a) was shown in \cite{gls2, gls5} and (b) was
shown in \cite{gls1,gls1a}. 

\begin{theorem}
\begin{itemize}
\item[(a)] $\varphi:\fl\la\to\mathbb{C}[U]\subset\mathbb{C}(U)$
  satisfies (M2$'$) and (M3$'$). 
\item[(b)] If $Q$ is Dynkin, then $\mathbb{C}[U]$
(respectively, $\mathbb{C}[U^{w_0}]$ for
the longest element $w_0$) is a cluster algebra
  modelled by a strong cluster map
  $\varphi:\mod\la\to\mathbb{C}(U)$ for the
  standard component of the cluster tilting graph of $\mod\la$ with no (respectively, all) coefficients inverted.
\end{itemize}
\end{theorem}


The image of a substructure $\B$ of $\mod \la$
gives a subcluster algebra of  $\mathbb{C}[U]$ for the Dynkin case, and we
illustrate this with the examples 
from \ref{c1_sec3}. We omit the calculation involved in proving
the isomorphisms between the subcluster algebras arising from $\B$ and 
the coordinate rings of the varieties under consideration.
See \cite{bl} for general background on Schubert varieties and (isotropic) Grassmannians.
For a
subset $J$ of size $k$ in $[1 \dots n]$ the symbol $[J]$ will denote
the $k \times k$ matrix minor of an $n \times n$ matrix with row set
$[1 \dots k]$ and column set $J$.
In the first two examples $G$ is $\SL_n(\mathbb{C})$ and $U$ is the subgroup of all upper triangular 
$n \times n$ unipotent matrices.

\bigskip
\noindent
\textbf{Example 1} ($\Gr_{2,5}$-Schubert variety)

Let $\la$ be of type $A_4$, and let $\B$ be the full additive subcategory
of $\mod \la$ from Example 1 in \ref{c1_sec3}. The associated algebraic group is
then $\SL_5(\mathbb{C})$. Consider the Grassmannian $\Gr_{2,5}$, and the Schubert
variety $X_{3,5}$ associated with the subset $\{3,5 \}$ of
$\{1,2,3,4,5\}$. Let $w_{3,5} = \begin{pmatrix} 1 & 2 & 3 & 4 & 5 \\ 3 &4 & 1 & 5 & 2
\end{pmatrix}$ be the associated Grassmann permutation 
in $S_5$, and $U^{w_{3,5}}$ the unipotent cell in $U$ associated to 
$w_{3,5}$. 
Note that the
Schubert variety $X_{3,5}$ is birationally isomorphic to the unipotent
cell $U^{w_{3,5}}$ \cite{BZ}.
Then $\mathbb{C}[U^{w_{3,5}}]$ is known to be a subcluster algebra
of $\mathbb{C}[U]$ \cite{fz3}.

Under the GLS-map $\varphi$ from $\mod \la$ to $\mathbb{C}[U]$ one can
check that $\varphi(M_x) = [x]$, with $M_x$ as defined in Example 1 in
\ref{c1_sec3}. Since $\B$ has a cluster substructure of
$\mod \la$, we know 
that the image gives rise to a subcluster algebra of $\mathbb{C}[U]$.
Then the image of $\mathcal{B}$ under the strong cluster map $\varphi$
is precisely $\mathbb{C}\big[U^{w_{3,5}}\big]$.
To see this, we mutate a seed from \cite{fz3} which generates the cluster algebra
structure for $\mathbb{C}[U]$, to get a new seed which contains $\varphi(T)$ for
the cluster tilting object $T$ in $\B$ in Example 1 in
\ref{c1_sec3}. Then one proves that the image is
$\mathbb{C}[U^{w_{3,5}}]$ after a proper choice of which coefficients to invert.

\bigskip
\noindent
\textbf{Example 2} (Unipotent Cell in
$\SL_4\big(\mathbb{C}\big)$)

Let $\la$ be the preprojective algebra of type $A_3$, and $\B$ the subcategory
of $\mod \la$ from Example 2 in \ref{c1_sec3}. The associated algebraic group is 
$\SL_4(\mathbb{C})$. Let $U^w$ be the unipotent cell of the unipotent subgroup $U$,
associated with the permutation $w = \begin{pmatrix} 1 & 2 & 3 & 4 \\ 4 & 3 & 1 & 2
\end{pmatrix}$. It is shown in \cite{fz3} that $\mathbb{C}[U^{w}]$ is a cluster algebra
of type $A_2$, and implicitly that it is a subcluster algebra of $\mathbb{C}[U]$.
In view of Section \ref{c3_sec1}, the image $\varphi(\B)$ has a 
subcluster algebra structure modelled by $\B$. 
As in Example 1, one begins
by mutating a seed from
\cite{fz3} which generates the cluster algebra structure for $\mathbb{C}[U]$, such
that the new seed contains $\varphi(T)$. Then one proves that
the image is $\mathbb{C}[U^{w}]$, after a proper choice of which
coefficients to invert.

\bigskip
\noindent
\textbf{Example 3} (The $\SO_{8}(\mathbb{C})$-Isotropic Grassmannians (cf. \cite[10.4.3]{gls3}))

Let $\la$ be the preprojective algebra of the Dynkin quiver $D_4$. 
Let $\varrho$ be the $4 \times 4$ anti-diagonal matrix 
whose $i,j$ entry is $(-1)^{i}\delta_{i,5-j}$ and let
$J$ be the $8 \times 8$ anti-diagonal matrix, written in block form as
\[ \begin{pmatrix} 0 & \ \varrho \\ \\ \varrho^T & \ 0 \end{pmatrix}. \]
The {\em even special orthogonal group} $\SO_{8}(\mathbb{C})$ is the
group of $8 \times 8$ matrices
\[ \left\{ g \in \SL_{8}(\mathbb{C}) \ \Big| \ g^T J g = J \right\}. \]
The {\em maximal unipotent subgroup} $U$ of $\SO_{8}(\mathbb{C})$ consists
of all $8 \times 8$ matrices in $\SO_{8}(\mathbb{C})$ which are upper triangular
and unipotent, i.e. having all diagonal entries equal to $1$. A more
explicit description in terms of matrices in block form is
\begin{equation}\label{formel1}
U \ = \ \left\{
\begin{pmatrix} u & u \varrho v \\ \\ 0 & \varrho^T \big(u^{-1}\big)^T \varrho
\end{pmatrix} \ \Bigg| \begin{array}{l} \text{$u$ is upper triangular
    unipotent in } \SL_{4}(\mathbb{C})
\\ \text{$v$ is skew-symmetric in } {\rm M}_{4}(\mathbb{C})\end{array} \right\} \end{equation}

The isotropic 
Grassmannian $\Gr_{2,8}^{\iso}$ is the closed subvariety of the classical Grassmannian $\Gr_{2,8}$
consisting of all isotropic $2$-dimensional subspaces of $\mathbb{C}^8$.
Let $\widehat{\Gr}^{\iso}_{2,8}$ be the corresponding affine cone. Let $q \colon U \to
\widehat{\Gr}^{\iso}_{2,8}$ denote the map given by 
$q(u)= u_1 \wedge
u_2$, where $u_1, u_2$ are the first two rows
of $u$ in $U$, and let 
$q^{\ast} \colon \mathbb{C}\big[ \widehat{\Gr}_{2,8}^{\iso} \big]
\longrightarrow \mathbb{C}\big[U \big]$ be the associated
homomorphism of coordinate rings.

Let $\varphi \colon \mod \la \to \mathbb{C}[U]$ be the GLS $\varphi$-map. Then one can show that
\[ \begin{array}{lllllllll} 
 \varphi\big(M_{16}\big) &=[16] &\varphi\big(M_{24}\big) &=[24]
&\varphi\big(M_{25}\big) &=[25] &\varphi\big(M_{26}\big) &=[26] \\
 \varphi\big(M_{68}\big) &=[68] &\varphi\big(M_{18}\big) &=[18] 
&\varphi\big(M_-\big)    &= \psi_{\s -} &\varphi\big(M_+\big) &= \psi_{\s +} \\ 
\varphi\big(P_1\big)     &= [8]   &\varphi\big(P_2\big)  &= [78]
&\varphi\big(P_3\big)    &=  {\displaystyle [678] \over {\displaystyle \Pfaff_{[1234]}} }&\varphi\big(P_4\big) 
&= \Pfaff_{[1234]} \end{array} \]

Here $\Pfaff_{[1234]}$ 
denotes the Pfaffian of the $4 \times 4$ skew-symmetric part $v$ appearing in (\ref{formel1}) of the unipotent element 
and $\psi_{\pm} = {1 \over 2}\big( [18] - [27] + [36] \pm [45]\big)$.
The functions  ${\displaystyle [678] \over {\displaystyle \Pfaff_{[1234]}}}$ and $\Pfaff_{[1234]}$
are examples of generalized minors of type $D$ \cite{fz-d}.

In the notation of Example 3 in Section \ref{c1_sec3}, we have seen that
$$T= 
M_{16} \oplus M_{24} \oplus M_{25} \oplus M_{26} \oplus M_{68} \oplus M_{18}
\oplus M_{-} \oplus M_{+} \oplus P_2 $$ is a cluster tilting object in
$\B = \Sub P_2$,
which can be extended to a cluster tilting object 
$$\widetilde{T}= 
T \oplus P_1 \oplus 
P_3 \oplus P_4 $$
for $\mod \la$. One shows that the initial seed used in \cite{fz3},
which determines a cluster algebra structure for $\mathbb{C}[U]$,
is mutation equivalent to the initial seed determined by $\widetilde{T}$, which 
hence generates the same cluster algebra. Since the subcategory 
$\B$ of $\mod \la$ in \ref{c1_sec3}, Example 3, has a substructure, as defined in Chapter \ref{chap1}, 
the connected component of the cluster tilting graph of $\B$ containing $T$
determines a subcluster algebra $\A'$ of $\mathbb{C}[U]$
(where $[18], \psi_{\pm},[78]$ are taken as noninverted
coefficients). Then we can prove that $\A'$ coincides with 
$\Im q^{\ast}$ (which we conjecture to be true more generally).

Notice that we have the cluster algebra structure for $\widehat{\Gr}^{\iso}_{2,8}$ by adjoining the coefficient $[12]$ to $\Im q^{\ast}$.
\\
\hspace{7mm}

\subsection{Cluster Structure of the loop group $\SL_2( \mathcal{L} )$}\label{c3_sec3}
${}$ \\
Let $Q$ be a finite connected non-Dynkin quiver without loops, and
let $\la$ be the associated completed preprojective algebra and
$W$ the associated Coxeter group. Let $G$ be the associated Kac-Moody
group $G$ with a maximal unipotent subgroup $U$, and let $U^w$ be the
unipotent cell associated with $w\in W$.
%
Using the GLS-map $\varphi:\fl\la\to\mathbb{C}[U]$ and the
restriction map $\mathbb{C}[U]\to\mathbb{C}[U^w]$, define the induced
map
\[\varphi_w:\Sub \la/I_w\subset\fl\la\stackrel{\varphi}{\to}\mathbb{C}[U]\to\mathbb{C}[U^w].\]
Using our results from Chapter \ref{chap2} we know that the transcendence degree
$l(w)$ of $\mathbb{C}(U^w)$ is equal to the number of non-isomorphic
summands of a cluster tilting object in $\Sub \la/I_w$.
It is then natural to pose the following.

\begin{conjecture}\label{first conjecture}
For any $w \in W$, the coordinate ring $\mathbb{C}[U^w]$ is a cluster
algebra modelled by a strong cluster map $\varphi_w:\Sub
\la/I_w\to\mathbb{C}[U^w]$ for the standard component of the cluster
tilting graph of $\Sub \la/I_w$ with all coefficients inverted.
\end{conjecture}

Recall that any infinite reduced expression where all generators occur
an infinite number of times gives rise to a cluster tilting
subcategory with an infinite number of non-isomorphic indecomposable
objects. Since the GLS-map $\varphi:\fl\la\to\mathbb{C}[U]$ satisfies
(M2$'$) and (M3$'$), it is natural to ask the following.

\begin{question}
The coordinate ring $\mathbb{C}[U]$ contains a cluster algebra
modelled by $\varphi:\fl\la\to\mathbb{C}[U]$ for any
connected component of the cluster tilting graph of $\fl\la$.
\end{question}

As a support for Conjecture \ref{first conjecture}, we show that this
is the case when $Q$ is the Kronecker quiver 
$\xymatrix@R0.2cm{
  1\ar@<0.5ex>[r]\ar@<-0.5ex>[r]& 0
}$, 
and the length of $w$ is at most 4. Without loss of generality,
we only have to consider the case $w=w_i$ for $w_1=s_0$, $w_2=s_0
s_1$, $w_3=s_0 s_1 s_0$ and $w_4 = s_0 s_1 s_0 s_1$.
The cases $w_i$ for $i=1$ or $2$ are clear since $\Lambda/I_{w_i}$ is a
cluster tilting object in $\Sub\la/I_{w_i}$ and $\mathbb{C}[U^{w_i}]$
is generated by invertible coefficients.
For the case $w_i$ for $i=3$ or $4$, we have cluster tilting objects
$T_3 = P_{0,1} \oplus P_{1,2} \oplus P_{0,3}$ in $\Sub \la/I_{w_3}$
and $T_4 = P_{0,1} \oplus P_{1,2} \oplus P_{0,3} \oplus P_{1,4}$ in
$\Sub \la/I_{w_4}$, where $P_{i,k}=P_i/J^kP_i$ for $i=0,1$ and $k>0$. 

The Kac-Moody group $\widehat{\SL_2}(\mathcal{L})$ associated with the
Kronecker quiver is defined as the unique non-trivial central
extension
\[ 1 \longrightarrow \Bbb{C}^* \longrightarrow \widehat{\SL_2}(\mathcal{L})
\overset{\pi}{\longrightarrow} \SL_2(\mathcal{L}) \longrightarrow
1 \]
of the algebraic loop group $\SL_2(\mathcal{L})$ by $\Bbb{C}^*$ (see
\cite{kp,ps} for details.)
The group $\SL_2(\mathcal{L})$ consists of all $\mathcal{L}$-valued $2
\times 2$ matrices $g =(g_{ij})$ with determinant 1, where
$\mathcal{L}$ is the Laurent polynomial ring $\mathbb{C}[t,t^{-1}]$.
The maximal unipotent subgroup (respectively, unipotent cells) of
$\widehat{\SL_2}(\mathcal{L})$ is mapped isomorphically onto the
maximal unipotent subgroup $$ U  \ = \ \Bigg\{ g \in \begin{pmatrix}
1 + t\mathbb{C}[t] & \mathbb{C}[t] \\ \\ t\mathbb{C}[t] & 1 + t\mathbb{C}[t] 
\end{pmatrix} \ \Bigg| \ \det(g) = 1 \ \Bigg\} $$
(respectively, corresponding unipotent cells) of $\SL_2(\mathcal{L})$.
Hence we deal with $\SL_2(\mathcal{L})$ instead of
$\widehat{\SL_2}(\mathcal{L})$.
The torus $H$ is the subgroup of $\SL_2(\mathcal{L})$ where the
elements are the diagonal matrices of the form $\begin{pmatrix} a & 0
  \\ 0 & a^{-1} \end{pmatrix}$.
The Weyl group $W=\operatorname{Norm}(H)/H$ is generated by the two
non-commuting involutions $s_0$ and $s_1$. 

Each $\varphi(P_{i,k})$ for $i= 0,1$ and $k>0$ is a regular function 
by computing explicit determinental formulas. For $g = \begin{pmatrix}
  g_{11} & g_{12} \\ g_{21} & g_{22}\end{pmatrix} \in U$, let $T_g$ be
the $\mathbb{Z}\times \mathbb{Z}$ matrix whose $(M,N)$ entry is given
by the residue formula
\[ \big(T_g\big)_{M,N} \ = \ \Res {\Di {g_{rs} \over {t^{n-m+1} } }} \]
with $N = 2n + r$ and $M = 2m + s$ and with $n, m \in \mathbb{Z}$
and $r,s \in \{1, 2\}$.
Let $\Delta^\sigma_{k;i}(g)$ be the determinant of the $k \times k$
submatrix of $T_g$ whose row and column sets are given respectively by
\begin{eqnarray*}
 \text{rows}  &=&  
\Big\{ k-i,\ k-i-2, \ k-i-4, \ k - i - 6, \ \dots \ , \ 2 - k - i \Big\},\\
\text{columns} & = &
\Big\{k -i + 1, \ k-i, \ k-i-1, \ k-i-2, \ \dots \ , \ 2-i  \Big\}. 
\end{eqnarray*}
The following theorem is a special case of \cite[Theorem 2]{scott}.

\begin{theorem}
$\varphi(P_{0,k})=\Delta^\sigma_{k;1}$ and
$\varphi(P_{1,k})=\Delta^\sigma_{k;0}$ for any $k>0$.
\end{theorem}

In the present case of the Weyl group elements $w_3= s_0 s_1 s_0$ and $w_4= s_0 s_1 s_0 s_1$, the 
corresponding affine unipotent varieties $U^w$
are given by: 
\begin{equation}\label{unipotent cell}
\begin{array}{ll}
U^{w_3} \ = \ \left\{ \begin{pmatrix} 1 + At & B  \\ 
Dt + Et^2 & 1 + Ft  \end{pmatrix} \ \Bigg| \
\begin{array}{ll} A+F = BD &  \\ AF = BE & E \ne 0
\end{array} \right\} \\
U^{w_4} \ = \ \left\{ \begin{pmatrix} 1 + At & B + Ct \\ 
Dt + Et^2 & 1 + Ft + Gt^2 \end{pmatrix} \ \left| \
\begin{array}{ll} A+F = BD & \ \ AG = CE \\ AF - CD = BE - G 
& \ \ G \ne 0 \end{array} \right. \right\} \end{array} 
\end{equation}
In terms of the complex parameters $A, \dots, G$ we have  
\begin{equation}\label{eqn3}
\begin{array}{ll} \varphi\big(P_{0,1}\big) = \Delta^\sigma_{1;1} = D
& \varphi\big(P_{1,2}\big) = \Delta^\sigma_{2;0} = DF-E \\
\varphi\big(P_{0,3}\big) = \Delta^\sigma_{3;1} = DEF - D^2 G - E^2 
& \varphi\big( P_{1,4}\big) = \Delta^\sigma_{4;0} = G(DEF - D^2 G - E^2) 
\end{array} 
\end{equation}
We are now ready to prove the crucial result on transcendence bases.

\begin{proposition}\label{propIV2.4}
The collections $\big\{ \Delta^\sigma_{1;1}, \Delta^\sigma_{2;0},
\Delta^\sigma_{3;1} \big\} \text{ and } \big\{ \Delta^\sigma_{1;1}, \Delta^\sigma_{2;0},
\Delta^\sigma_{3;1}, \Delta^\sigma_{4;0} \big\}$ are respectively
transcendence bases for the rational function fields
$\mathbb{C}\big(U^{w_3}\big)$ and $\mathbb{C}\big(U^{w_4}\big)$.
\end{proposition}

\begin{proof} Consider first the case of $\big\{ \Delta^\sigma_{1;1}, \Delta^\sigma_{2;0},
\Delta^\sigma_{3;1}, \Delta^\sigma_{4;0} \big\}$ within $\mathbb{C}\big(U^{w_4}\big)$.
The transcendence degree of $\mathbb{C}\big(U^{w_4}\big)$ over $\mathbb{C}$ is 4 since  
$\dim_{\mathbb{C}} \, U^{w_4}$ is 4. In this case any
4 rational functions which generate the field
$\mathbb{C}\big(U^{w_4}\big)$ must constitute a transcendence basis for 
$\mathbb{C}\big(U^{w_4}\big)$. The global coordinates $A,B,C,D,E,F, G^{\pm 1}$
clearly generate $\mathbb{C}\big(U^{w_4}\big)$, and we will be finished if we
can express each global coordinate as a rational function
in $\Delta^\sigma_{1;1}$, $\Delta^\sigma_{2;0}$,
$\Delta^\sigma_{3;1}$, and $\Delta^\sigma_{4;0}$. To do this we
introduce five auxiliary functions $\widetilde{\Delta}^\sigma_{1;1}$,
$\widetilde{\Delta}^\sigma_{2;0}$, $\Psi$, $\Omega$, and $\Sigma$
defined implicitly by the formulas
\begin{equation}\label{eqn4}
\begin{array}{lllllll}
& \widetilde{\Delta}^\sigma_{1;1} \ \Delta^\sigma_{1;1} &= \
\Big( \Delta^\sigma_{2;0} \Big)^2 + \ \Delta^\sigma_{3;1} &
& \widetilde{\Delta}^\sigma_{2;0} \ \Delta^\sigma_{2;0} &= \
\Big( \Delta^\sigma_{1;1} \Big)^2 \Delta^\sigma_{4;0} + \ 
\Big(\Delta^\sigma_{3;1} \Big)^2 \\
& \Psi \ \Delta^\sigma_{2;0}   &= \ \Big(\widetilde{\Delta}^\sigma_{1;1}\Big)^2 
+ \ \Delta^\sigma_{4;0} &
& \Omega \ \Delta^\sigma_{1;1} &= \ \Big(\widetilde{\Delta}^\sigma_{2;0}\Big)^2          
+ \ \Big(\Delta^\sigma_{3;1} \Big)^3 \\
& \Sigma \ \widetilde{\Delta}^\sigma_{2;0} &= \ \Big(\Delta^\sigma_{3;1} \Big)^4 
\Delta^\sigma_{4;0}  + \ \Omega^2 \end{array} \end{equation}
\noindent
These are exactly the cluster variables obtained from the initial cluster from at most two mutations.
Evidently these five functions can be rationally expressed in terms
of $\Delta^\sigma_{1;1}$, $\Delta^\sigma_{2;0}$,
$\Delta^\sigma_{3;1}$, and $\Delta^\sigma_{4;0}$. 
Moreover, after using the equations \eqref{eqn3} and carrying out the divisions arising in 
solving the above system, each function 
can be written as a polynomial in the
global coordinates $A, \dots, G^{\pm 1}$:
\begin{equation}\label{eqn5}
\begin{array}{lll} 
&\widetilde{\Delta}^\sigma_{1;1} &= \ DF^2 - EF - DG \\
&\widetilde{\Delta}^\sigma_{2;0} &= \ E \big( DEF - D^2G -E^2 \big) \\
&\Psi                            &= \ \big(F^2 - G \big)\big(DF - E\big) - DFG \\
&\Omega                          &= \ \big(EF - DG\big)\big(DEF - D^2G - E^2\big)^2 \\
&\Sigma                          &= \ \big(EF^2 - DFG - EG\big)\big(DEF - D^2G - E^2\big)^3 
\end{array} 
\end{equation}
\noindent
The following identities can be easily checked by using the relations
\eqref{unipotent cell}.
\begin{equation}\label{eqn6}
\begin{array}{lllll}  
&A = 
{\Di {  \Omega \widetilde{\Delta}^\sigma_{2;0} \over { \Big(\Delta^\sigma_{3;1}\Big)^{4} } }} 
&B= {\Di { \Sigma \over { \Big(\Delta^\sigma_{3;1}\Big)^{4}}} }
&C ={\Di { \Omega \Delta^\sigma_{4;0} \over { \Big(\Delta^\sigma_{3;1}\Big)^{4}}} } 
&D = \Delta^\sigma_{1;1} 
\\
&E= {\Di { \widetilde{\Delta}^\sigma_{2;0} \over { \Delta^\sigma_{3;1} }} }
&F = { \Di {\Sigma \Delta^\sigma_{1;1}  + \ \Omega \widetilde{\Delta}^\sigma_{2;0} \over
{ \Big(\Delta^\sigma_{3;1}  \Big)^{4}}} } 
&G^{\pm 1} = 
\Bigg( {\Di { \Delta^\sigma_{4;0} \over { \Delta^\sigma_{3;1} } } } \Bigg)^{\pm 1}
\end{array} 
\end{equation}
\noindent
Thus $\big\{ \Delta^\sigma_{1;1}, \Delta^\sigma_{2;0},
\Delta^\sigma_{3;1}, \Delta^\sigma_{4;0} \big\}$ is a transcendence basis.

\medskip
\noindent
The case of $\big\{ \Delta^\sigma_{1;1}, \Delta^\sigma_{2;0},
\Delta^\sigma_{3;1} \big\}$ is handled similarly. Over 
the unipotent cell $U^{w_3}$ the global
coordinates $A,B,D,E^{\pm 1},F$ are rationally expressed as
\[ A= {\Di {\widetilde{\Delta}^\sigma_{1;1} \over 
{\Big( \Delta^\sigma_{2;0} \Big)^2 }} }  \quad
B = {\Di {\Big(\widetilde{\Delta}^\sigma_{1;1} \Big)^2 \over 
{\Big( \Delta^\sigma_{2;0} \Big)^3 }} } \quad 
D = \Delta^\sigma_{1;1} \quad 
E^{\pm 1} = \Bigg( {\Di {\Delta^\sigma_{3;1} \over
{\Delta^\sigma_{2;0} }} }  \Bigg)^{\pm 1} \quad
F = {\Di { \widetilde{\Delta}^\sigma_{1;1} \over 
{\Delta^\sigma_{2;0} }} }  \]
\noindent
from which it follows that $\big\{ \Delta^\sigma_{1;1}, \Delta^\sigma_{2;0},
\Delta^\sigma_{3;1} \big\}$ generates the rational function field
$\mathbb{C}\big(U^{w_3}\big)$. The transcendence degree of 
$\mathbb{C}\big(U^{w_3}\big)$ is 3 and consequently
the collection $\big\{ \Delta^\sigma_{1;1}, \Delta^\sigma_{2;0},
\Delta^\sigma_{3;1} \big\}$ is a transcendence basis.
\end{proof}

\begin{theorem}\label{}
For $w=w_3=s_0s_1s_0$ or $w=w_4=s_0s_1s_0s_1$, the coordinate ring
$\mathbb{C}[U^w]$ is a cluster algebra modelled by a strong cluster
map $\varphi_w\colon\Sub\la/I_w\to\mathbb{C}[U^w]$ for the standard
component of the cluster tilting graph of $\Sub \la/I_w$ with all
coefficients inverted.
\end{theorem}

\begin{proof}
We treat the case $w=w_4$, and leave the other easier case to the reader.
Since $\varphi_{w_4}:\Sub\la/I_{w_4}\to\mathbb{C}[U^{w_4}]$ is a strong
cluster map by Proposition \ref{propIV2.4}, we have a cluster algebra
$\A^{w_4} \subset\mathbb{C}[U^{w_4}]$ where we invert all
coefficients. 

We start with proving the inclusion $\A^{w_4} \subset
\mathbb{C}\big[U^{w_4} \big]$. By construction the functions 
$\Delta^\sigma_{1;1}$, $\Delta^\sigma_{2;0}$,
$\Delta^\sigma_{3;1}$, and $\Delta^\sigma_{4;0}$ are
regular. The
coefficients $\Delta^\sigma_{3;1}$ and $\Delta^\sigma_{4;0}$
are invertible, so we must verify that their
inverses are in $\mathbb{C}\big[U^{w_4} \big]$. 
Put
\[\Delta_{3;1} = \ BCF - B^2G - C^2,\ \ \ \Delta_{4;0}= \ G\big(BCF -B^2G - C^2 \big).\]
Using relations \eqref{unipotent cell} one has
$\Delta^\sigma_{3;1} \ \Delta_{3;1} \ = \ \Delta^\sigma_{4;0} \ \Delta_{4;0} \ = \ G^4$.
The function $G \ne 0$ over $U^{w_4}$, and consequently $\Delta^\sigma_{3;1}$
and $\Delta^\sigma_{4;0}$, are invertible with inverses
given by 
\[ \Big( \Delta^\sigma_{3;1} \Big)^{-1} = { \Di \Delta_{3;1} \over
{\Di G^4}} \qquad \Big(  \Delta^\sigma_{4;0} \Big)^{-1} =
{\Di \Delta_{4;0} \over {\Di G^4 } }. \]
\noindent
The two exchange relations for the initial seed are precisely the
first two relations given  in \eqref{eqn4}, namely 
\[\begin{array}{llllll}
& \widetilde{\Delta}^\sigma_{1;1} \ \Delta^\sigma_{1;1} &= \
\Big( \Delta^\sigma_{2;0} \Big)^2 + \ \Delta^\sigma_{3;1},
&& \widetilde{\Delta}^\sigma_{2;0} \ \Delta^\sigma_{2;0} &= \
\Big( \Delta^\sigma_{1;1} \Big)^2 \Delta^\sigma_{4;0} + \ 
\Big(\Delta^\sigma_{3;1} \Big)^2
\end{array} \]
\noindent
and we have seen in \eqref{eqn5}  that we have the following expressions for the cluster variables
 $\widetilde{\Delta}^\sigma_{1;1}$ and 
$\widetilde{\Delta}^{\sigma}_{2;0}$:
\begin{equation}
\begin{array}{llllll} 
&\widetilde{\Delta}^\sigma_{1;1} &= \ DF^2 - EF - DG,
&&\widetilde{\Delta}^\sigma_{2;0} &= \ E \big( DEF - D^2G -E^2 \big) 
\end{array} 
\end{equation}
which are clearly regular functions on $U^{w_4}$.
Define $x_1 = \Delta^\sigma_{1;1}$,
$x_2 = \Delta^\sigma_{2;0}$, $x_3 = \Delta^\sigma_{3;1}$, $x_4 =
\Delta^\sigma_{4;0}$, and define also $\widetilde{x}_1 =
\widetilde{\Delta}^\sigma_{1;1}$, and $\widetilde{x}_2
= \widetilde{\Delta}^\sigma_{2;0}$. Select an arbitrary
cluster variable $x$ in $\A^{w_4}$. 
By the Laurent phenomenon we know that $x$ can
be expressed as a polynomial in 
$x_1^{\pm 1}$, $x_2^{\pm 1}$, $x_3$, and $x_4$. 
>From this it follows that 
$x$ is regular on the Zariski open set $\mathcal{U}$ of
$U^{w_4}$ defined by
\[ \mathcal{U} = \ \Big\{ g \in U^{w_4} \ \Big| \ x_i(g) \ne 0 
\ \text{for $i = 1, 2$ }  \Big\} \]
\noindent
Again, by the Laurent phenomenon, $x$ can also be expressed, simultaneously,
as a polynomial in both $\Big\{ \widetilde{x}_1^{\pm 1}, \ x_2^{\pm 1}
, \ x_3, \ x_4 \Big\}$ and 
$\Big\{ x_1^{\pm 1}, \ \widetilde{x}_2^{\pm 1}, \ x_3,
\ x_4 \Big\}$. Consequently $x$ is simultaneously regular on both of the
Zariski open subsets
\[ \begin{array}{lll}
&\mathcal{U}^{(1)} &= \ \Big\{ g \in U^{w_4} \ \Big| \ \widetilde{x}_1(g) \ne 0
\ \text{and} \ x_2(g) \ne 0 \Big\} \\
&\mathcal{U}^{(2)} &= \ \Big\{ g \in U^{w_4} \ \Big| \ x_1(g) \ne 0
\ \text{and} \ \widetilde{x}_2(g) \ne 0 \Big\}
\end{array} \]
\noindent
Taken all together we conclude that $x$ is regular over the union
$\mathcal{U} \cup \ \mathcal{U}^{(1)} \cup \ \mathcal{U}^{(2)}$.
The complement $V$ of this union inside $U^{w_4}$ consists of those 
points $g$ for which $x_1(g) = \widetilde{x}_1(g) = 0$
and/or $x_2(g) = \widetilde{x}_2(g) = 0$. Using the 
expressions for $\widetilde{x}_1$ and $\widetilde{x}_2$ given (in \eqref{eqn3}) we see that 
\[ V \ = \ \left\{ \begin{pmatrix} 1 & B \\ Et^2 & 1 + Gt^2 \end{pmatrix}
\ \Bigg| \begin{array}{ll} &BE = G \\ &G \ne 0 \end{array}
\right\} \bigcup 
\left\{ \begin{pmatrix} 1 & Ct \\ Dt & 1 + Gt^2 \end{pmatrix}
\ \Bigg| \begin{array}{ll} &CD = G \\  &G \ne 0 \end{array}
\right\}
\]
\noindent
The only hypothetical singularities of $x$ must lie in $V$,
and since the codimension of $V$ inside $U^{w_4}$ is clearly
2, it follows that $x$ has no singularities and is thus
regular on all of $U^{w_4}$. The cluster algebra
$\A^{w_4}$ is generated by the coefficients and
cluster variables and so we conclude that
$\A^{w_4} \subset \mathbb{C}\big[ U^{w_4}\big]$.

\medskip
\noindent
The reverse inclusion $\mathbb{C}\big[ U^{w_4}\big] \subset
\A^{w_4}$ follows from the fact 
the functions $x$, $\Omega$, and $\Sigma$ 
defined previously are also cluster
variables of $\A^{w_4}$ 
and the fact that the global coordinates 
$A, \dots, G^{\pm 1}$ are generated inside of $\A^{w_4}$
by using the expressions in \eqref{eqn6}.
\end{proof}

\medskip
Finally we show the cluster graphs for $\mathbb{C}[U^{w_3}]$ and
$\mathbb{C}[U^{w_4}]$.
Then $\mathbb{C}[U^{w_3}]$ has a seed consisting of a cluster $\{ \Delta^\sigma_{1;1}
\}$, coefficients $\{ \Delta^\sigma_{2;0},\Delta^\sigma_{3;1} \}$ and the quiver
$$\xymatrix@C0.4cm@R0.4cm{
&\Delta^\sigma_{2;0} &\\
\Delta^\sigma_{1;1}\ar@<0.5ex>[ru]\ar@<-0.5ex>[ru]&& \Delta^\sigma_{3;1}\ar[ll]
}$$
\noindent
which is the quiver of $\End_\la(T_3)$, where we drop the arrows between coefficients. The cluster graph is 
$$\xymatrix@C0.4cm@R0.4cm{
\{ \Delta^\sigma_{1;1} \} \ar@{-}[r] & \{ \widetilde{\Delta}^\sigma_{1;1} \}}
$$
\noindent
where the only other cluster variable $\widetilde{\Delta}^\sigma_{1;1}$ is determined by 
$\Delta^\sigma_{1;1} \widetilde{\Delta}^\sigma_{1;1} = (\Delta^\sigma_{2;0})^2 + \Delta^\sigma_{3;1}$.
The cluster type of $\mathbb{C}[U^{w_3}]$ is $A_1$.

Similarly, $\mathbb{C}[U^{w_4}]$ has a seed consisting of a cluster
$\{ \Delta^\sigma_{1;1},\Delta^\sigma_{2;0}\}$, coefficients $\{ \Delta^\sigma_{3;1},\Delta^\sigma_{4;0} \}$ and the quiver
$$\xymatrix@C0.4cm@R0.4cm{
&\Delta^\sigma_{2;0}\ar@<0.5ex>[rd]\ar@<-0.5ex>[rd]& & \Delta^\sigma_{4;0}\ar[ll] \\
\Delta^\sigma_{1;1}\ar@<0.5ex>[ru]\ar@<-0.5ex>[ru]&& \Delta^\sigma_{3;1}\ar[ll] & 
}$$
The cluster graph is 
$$\xymatrix{\cdots\text{ }\cdot \ar@{-}[r] & \cdot \ar@{-}[r] & \{ \Delta^\sigma_{1;1}, \Delta^\sigma_{2;0} \} \ar@{-}[r] & \cdot \ar@{-}[r] &\cdot \text{ }\cdots}$$
The cluster type of $\mathbb{C}[U^{w_4}]$ is $\widehat{A_1}$.

Note that this gives an example of a substructure of a cluster 
structure coming from the inclusion $\Sub\la/I_{w_3}\subset\Sub\la/I_{w_4}$, 
and a cluster map such that we get a subcluster algebra of a cluster algebra, 
namely $\mathbb{C}[U^{w_3}]$ as a subcluster algebra of $\mathbb{C}[U^{w_4}]$.


\end{document}